\documentclass[preprint]{imsart}
\RequirePackage[OT1]{fontenc}
\RequirePackage{amsthm,amsmath,amsfonts,amssymb}
\RequirePackage[numbers]{natbib}
\RequirePackage[colorlinks,citecolor=blue,urlcolor=blue]{hyperref}
\RequirePackage{graphicx}

\usepackage{verbatim}
\usepackage{multirow}
\usepackage{color, subfigure}

\theoremstyle{plain}

\newtheorem{theorem}{Theorem}[section]
\newtheorem{lemma}[theorem]{Lemma}

\theoremstyle{remark}

\newtheorem{prop}[theorem]{Proposition}
\newtheorem{cor}[theorem]{Corollary}

\newtheorem{thm}[theorem]{Theorem}

\newtheorem*{lemma*}{Lemma}

\newcommand{\beq}{\begin{equation}}
\newcommand{\eeq}{\end{equation}}

\newcommand{\var}{\mathrm{Var}}

\newcommand{\eps}{\epsilon}

\newcommand{\bitem}{\begin{itemize}}
\newcommand{\eitem}{\end{itemize}}

\newcommand{\goto}{\rightarrow}

\newcommand{\beqn}{\begin{equation}}
\newcommand{\eeqn}{\end{equation}}
\newcommand{\balign}{\begin{align}}
\newcommand{\ealign}{\end{align}}

\newcommand{\be}{\begin{equation}}
\newcommand{\ee}{\end{equation}}
\newcommand{\ba}{\begin{array}}
\newcommand{\ea}{\end{array}}
\def\ben{\begin{equation*}}
\def\een{\end{equation*}}
\def\bea{\begin{eqnarray}}
\def\eea{\end{eqnarray}}
\usepackage{amssymb}

\newcommand{\tf}{\tilde{f}}

\begin{document}


\begin{frontmatter}
\title{High Dimensional Quadratic Discriminant Analysis: Optimality and Phase Transitions}
\runtitle{QDA: Optimality and Phase Transitions}
\thankstext{T1}{Equal contribution.}

\begin{aug}
\author[A]{\fnms{Wanjie} \snm{Wang}\thanksref{t1}\ead[label=e1,mark]{wanjie.wang@nus.edu.sg}},
\author[B]{\fnms{Jingjing} \snm{Wu}\thanksref{t2}\ead[label=e2]{jinwu@ucalgary.ca}}
\and 
\author[A]{\fnms{Zhigang} \snm{Yao}\thanksref{t3}\ead[label=e3,mark]{zhigang.yao@nus.edu.sg}}

\thankstext{t1}{Supported by MOE Start-up grant 155-000-173-133 and Tier 1 grant 0004813-00-00.}
\thankstext{t2}{Supported by NSERC Discovery Grants RGPIN-2018-04328.}
\thankstext{t3}{Supported by MOE Tier 2 grant 155-000-184-114 and Tier 1 grant 0004813-00-00.}


\address[A]{Department of Statistics and Data Science, 
National University of Singapore, 
\printead{e1,e3}}
\address[B]{Department of Mathematics and Statistics,
University of Calgary,
\printead{e2}}
\end{aug}

\begin{abstract}

Consider a two-class classification problem where we observe samples $(X_i, Y_i)$ for $i = 1, \cdots, n$, $X_i \in \mathcal{R}^p$ and $Y_i \in \{0, 1\}$. 
Given $Y_i = k$, $X_i$ is assumed to follow a multivariate normal distribution with mean $\mu_k \in \mathcal{R}^k$ and covariance matrix $\Sigma_k$, $k=0,1$. Supposing a new sample $X$ from the same mixture is observed, our goal is to estimate its class label $Y$. Such a high-dimensional classification problem has been studied thoroughly when $\Sigma_0 = \Sigma_1$. However, the discussions over the case $\Sigma_0 \neq \Sigma_1$ are much less over the years. 

This paper presents the quadratic discriminant analysis (QDA) for the weak signals (QDAw) algorithm, and the QDA with feature selection (QDAfs) algorithm.
QDAfs applies Partial Correlation Screening in \cite{PCS} to estimate $\hat{\Omega}_0$ and $\hat{\Omega}_1$,  and then applies a hard-thresholding on the diagonals of $\hat{\Omega}_0 - \hat\Omega_1$. QDAfs further includes the linear term $d^\top X$, where $d$ is achieved by a hard-thresholding on $\hat{\Omega}_1\hat{\mu}_1 - \hat{\Omega}_0\hat{\mu}_0$. QDAfs achieves theoretical optimality and outperforms recent works on the linear discriminant analysis of high-dimensional data on a real data set. 

We further propose the rare and weak model to model the signals in $\Omega_0 - \Omega_1$ and $\mu_0 - \mu_1$. Based on the signal weakness and sparsity in $\mu_0 - \mu_1$, we propose two ways to estimate labels: 1) QDAw for weak but dense signals; 2) QDAfs for relatively strong but sparse signals. 
We figure out the classification boundary on the 4-dim parameter space: 1) Region of possibility, where either QDAw or QDAfs will achieve a mis-classification error rate of 0; 2) Region of impossibility, where all classifiers will have a constant error rate. 
 The numerical results from real datasets support our theories and demonstrate the necessity and superiority of using QDA over LDA for classification.

\end{abstract}

\begin{keyword}[class=MSC]
\kwd[Primary ]{62H30}
\kwd[; secondary ]{62F05}
\end{keyword}

\begin{keyword}
\kwd{High-dimensional classification}
\kwd{quadratic discriminant analysis}
\kwd{phase transitions}
\kwd{rare and weak signals}
\end{keyword}

\end{frontmatter}


\section{Introduction}
Consider a two-class classification problem, where we have $n$ labeled training samples $(X_i,Y_i)$, $i=1, \dots, n$. Here, $X_i$'s are $p$-dimensional feature vectors and $Y_i\in\{0,1\}$ are the corresponding class labels. $X_i$ is assumed to have mean $\mu_k$ and covariance matrix $\Sigma_k$ where $k = Y_i \in \{0, 1\}$. The goal is to estimate the label of a new observation $X$. A significant amount of work has been done in this field; see \citep{Anderson, RegularizedQDA, lachenbruch1979discriminant}.

Fisher's Linear Determinant Analysis (LDA) in \cite{Fisher} utilizes a weighted average of the features of the test sample to make a prediction. The optimal weight vector for LDA where the two classes are assumed to share the same correlation structure $\Omega = \Sigma_0^{-1}  = \Sigma_1^{-1}$ satisfies
\begin{equation} \label{Fisherw}
d \propto \Omega  (\mu_1 - \mu_0).
\end{equation}
When $n \gg p$, the mean vectors $\mu_0$ and $\mu_1$, and the precision matrix $\Omega$ can be easily estimated. Therefore Fisher's LDA is approachable. 

In modern analytical approaches, high-dimensional data have flooded, where a prominent number of measurements of features, often in the millions, are gathered for a single subject (\citep{large-infer}).
Although the number of features is huge, usually only a small portion of them are regarded as relevant to the classification decision, but these are not known in advance. 
In this sense, the traditional classification methods will lose power because of the large amount of noise. 
Methods have been proposed to reduce the noises in methods; see \citep{FAIR, ROAD, FJY, Ingster, Wasserman}. 
For such high dimensional classification problems, recent developments provide methods to estimate the precision matrix $\Omega$ in LDA; see \cite{CLIME, PCS}. 

However, LDA still faces two problems:
\begin{itemize}
\item It does not account for the information from $\Omega_0$ and $\Omega_1$. The features in two classes may not share the same conditional independence structure, however, LDA doesn't take this information into consideration. Such difference will impact the distribution of the new variable $X$ and the estimation of $d$. 
\item The theoretical analysis for the case $\Omega_0 \neq \Omega_1$ is short of discussion. Actually, 
the analysis of error rates is quite complicated for the high dimensional model, even in the simplest case where $\Omega_0= \Omega_1 = I_p$, as the signals are rare and weak. An extensive discussion may be found in  \cite{DJ08, Jin-survey1, FJY, JinPNAS, Jin-survey2, JKW}.
\end{itemize}
These two problems motivates us to derive an algorithm and relative theoretical analysis for the case $\Omega_0 \neq \Omega_1$.

Inherited from the previous studies, many applications in high-dimensional classifier problems share the same aspects: (1) The signals in the mean vector are comparatively rare and weak; and (2) the precision matrices $\Omega_0$ and $\Omega_1$ are sparse. 
Such sparsity properties guide us to propose rare and weak model about them and solve the problem. 

In this paper, we propose a rare and weak model for the high-dimensional classification problem. In the new model, we allow $\Omega_0 \neq \Omega_1$ under both the diagonals and off-diagonals with different parameterizations. It is a more reasonable fit with real data than the LDA model. 
To parameterize $(\mu_0, \mu_1)$ and $(\Omega_0, \Omega_1)$, we normalize the data by centering the mean vector and scaling features to have unit variance. These are standard processing steps in the real data analysis, which will be discussed later with more details. In theoretical analysis, we will see that using $\Omega_0 - \Omega_1$ improves the classification accuracy.

This paper further explores the high-dimensional classification problem in the following aspects:  
\begin{itemize}
    \item We propose a rare and weak model to account for both $(\mu_0, \mu_1)$ and $(\Omega_0, \Omega_1)$. It models the weakness and sparsity of $\mu_0 - \mu_1$ and the diagonals and off-diagonals of $\Omega_0 - \Omega_1$.
    We tie all the parameters to $p$, which allows us to explore the phase diagram of the classification problem. 
    \item We propose several algorithms: the Quadratic Discriminant Analysis with Feature Selection (QDAfs) algorithm that works when the signals in $\mu_0 - \mu_1$ are relatively strong and sparse; and the Quadratic Discriminant Analysis for Weak signals (QDAw) algorithm that works when the signals in $\mu_0 - \mu_1$ are weak and relatively sparse. 
    \item We derive the phase diagram for the high-dimensional classification problem when $\Omega_0 \neq \Omega_1$. We find that the region that QDAfs and QDAw will give satisfactory classification results. We also find out the region that no classifier will succeed, i.e. region of impossibility. Under the rare and weak model, the success region of QDAw/QDAfs and the region of impossibility can form the whole phase diagram, which proves the optimality of QDAw/QDAfs. 
\end{itemize}
Our method performs better than the optimal high-dimensional LDA method on a real data set in the numerical analysis section, which suggests that the second-order information should be incorporated to improve the classification results.

\subsection{Quadratic Discriminant Analysis on high-dimensional data}\label{sec:qda}
Consider the two-class classification problem mentioned above, where we observe $n$ training samples $(X_i,Y_i)$, $i=1, \dots, n$. 
Given $Y_i=k$, we assume the feature vector $X_i \in \mathcal{R}^p$ follows a multivariate normal distribution with mean $\mu_k$ and covariance matrix $\Sigma_k = \Omega_k^{-1}$. 
Let $X$ denote an independent test sample from the same population; then, 
\be \label{T0}
X|Y \sim (1 - Y) N(\mu_0, \Omega_0^{-1})+ Y N(\mu_1,\Omega_1^{-1}).
\ee
We would like to classify $X$ as being from either $Y=0$ or $Y=1$. 

For two-population classification problems, the QDA method is commonly used to exploit both the mean and covariance information; see \cite{RegularizedQDA, mclachlan2004}. Consider the ideal case that both $\mu_k$ and $\Omega_k$ are known, $k = 0, 1$. The ratio of the likelihood functions in two classes gives the optimal classifier, which results in the QDA method, that
\be \label{CR2}
\ba{l}
\displaystyle \hat{Y} = I \left\{X^\top\left (\Omega_0 -\Omega_1\right )X-2\bigl(\mu_0^\top\Omega_0 -\mu_1^\top\Omega_1 \bigr)X\right. \\
\displaystyle \left. \ \ \ + \bigl(\mu_0^\top\Omega_0\mu_0- \mu_1^\top\Omega_1\mu_1+\ln|\Omega_1|- \ln|\Omega_0|\bigr) \ > \ 0\right\} ,
\ea
\ee 
where $I(A)$ is the indicator function of event $A$. If $P(Y_i = 1) \neq 0.5$, an additional term $2\ln \tfrac{P(Y_i = 1) }{1-P(Y_i = 1) }$ on the right-hand side of the inequality will improve the accuracy. However, as such term does not have effects on the possibility and impossibility regions, our analysis applied to the balanced case (i.e. $P(Y_i = 1) =0.5$) will suffice.


For real data, the parameters $\mu_0$, $\mu_1$, $\Omega_0$ and $\Omega_1$ are all unknown. Estimation of them is challenging in the high-dimensional setting with $p \gg n$. 
There are various extensions of QDA for the high-dimensional data; see \cite{aoshima2019high, wu2019quadratic, RFQD}. When the signals are sparse, it is modified accordingly with sparsity assumptions on $\Sigma_1$, $\Sigma_0$ and $\mu_1 - \mu_0$; see \cite{fan2015innovated, jiang2018direct, li2015sparse}. 
Since sparsity assumptions on the precision matrices $\Omega_0$ and $\Omega_1$ are more commonly seen in applications (\cite{lauritzen1996graphical, yu2017learning, yuan2007model}), which means sparse conditional dependency between features, we propose an approach based on the sparsity of $\Omega_0$, $\Omega_1$ and $\mu$. 

We begin with estimating these parameters in the high-dimensional setting. 
\begin{itemize}
    \item Updated Partial Correlation Screening (PCS) approach on precision matrices $\Omega_0$ and $\Omega_1$.
    \begin{itemize}
    \item[Step 1.] Estimate $\Omega_0$ and $\Omega_1$ by PCS in \cite{PCS}, denoted as $\hat{\Omega}_0$ and $\hat{\Omega}_1$. 
    \item[Step 2.] Let $\hat\Omega_{\rm diff} = \hat{\Omega}_0 - \hat{\Omega}_1$. For each diagonal $\hat\Omega_{\rm diff}(i,i)$, update it as $\hat\Omega_{\rm diff}(i,i) = \hat\Omega_{\rm diff}(i,i)1\{|\hat\Omega_{\rm diff}(i,i)| \leq 2\sqrt{2\ln p/n}\}$. 
    \end{itemize}
    \item Estimation of the linear component $(\Omega_0\mu_0 - \Omega_1\mu_1)^\top X$. 
       \begin{itemize}
    \item[Weak signals:] Let $\hat{\mu} = a*{\bf 1}$, where $\bf 1$ is a vector with all ones and $a$ is a constant. Estimate it by $g_w(X; \hat{\mu}, \hat{\Omega}_0, \hat{\Omega}_1) = \hat{\mu}^\top (\hat{\Omega}_0 + \hat{\Omega}_1)X$. 
    \item[Strong signals:] Let $d = \hat\Omega_1\hat{\mu}_1 - \hat{\Omega}_0 \hat{\mu}_0$, where $\hat{\mu}_1$ and $\hat{\mu}_0$ are the sample mean vectors of class 1 and 0. Let $d^{(t)}(j) = 1\{|d(j) \geq t|\}$. 
    Estimate it by $g_s(X; \hat{\mu}, \hat{\Omega}_0, \hat{\Omega}_1) = (d\circ d^{(t)})^{\top} X$. 
    \end{itemize}
\end{itemize}

With the estimates $\hat\Omega_0$, $\hat\Omega_1$ and $g(X; \hat{\mu}, \hat{\Omega}_0, \hat{\Omega}_1)$, the high-dimensional QDA is proposed in Table \ref{tab:alg1}. We call the algorithm with $g_w(X; \hat{\mu}, \hat{\Omega}_0, \hat{\Omega}_1)$ as QDA with weak signals (QDAw), and the algorithm with $g_s(X; \hat{\mu}, \hat{\Omega}_0, \hat{\Omega}_1)$ as QDA with a feature-selection step (QDAfs). 
\begin{table}[ht!]
\caption{Algorithm 1: Pseudocode for QDA on high-dimensional data} 
\begin{tabular}{ll}\\ \hline
& \underline{Input}: data points $(X_i, Y_i)$, $1 \leq i \leq n$;  threshold $t > 0$; new data point $X$. \\  
& \underline{Output}: label $\hat{Y}$. \\
1. &Parameter Estimation: Estimate $\hat{\Omega}_0$, $\hat{\Omega}_1$, $\hat\Omega_{\rm diff}$, $\hat{\mu}$ and $d$ according to the procedure as above.  \\ 
2.  & Define a constant $C$ according to the parameters; details in the algorithms in later sections. \\
3. & Let $g(X) = g_s(X)$ if $\mu_1 - \mu_0$ has strong signals or $g(X) = g_w(X)$ if $\mu_1 - \mu_0$ has only weak signals.\\
4.  & \underline{QDA Score}: Calculate the QDA score $Q(X)=X^\top \Omega_{\rm diff} X + 2g(X) + C$. \\
5. &\underline{Prediction}: Predict $\hat Y = I\{Q(X) > 0\}$.  \\
 \hline
\end{tabular}
\label{tab:alg1} 
\end{table}

There are multiple high-dimensional precision matrix estimation methods; see \cite{CLIME, fan2016overview, PCS}. We estimate $\Omega_0$ and $\Omega_1$ with the Partial Correlation Screening (PCS) approach in \cite{PCS} for this algorithm. PCS has good control on the Frobenius norm of $\hat{\Omega}_k - \Omega_k$, which is the main factor in the error analysis of QDA. 
The thresholding step on the diagonals of $\hat{\Omega}_0 - \hat{\Omega}_1$ is as an adjustment on the element-wise error of PCS, at the order of $\sqrt{\ln p/n}$. Without the thresholding step, the random error in $X^\top \hat\Omega_{\rm diff} X$ is large enough to cover the truth if signals in $\Omega_{\rm diff}$ and $\mu_0 - \mu_1$ are weak. The theoretical limit of QDA with PCS can be found in Proposition \ref{prop:unknownomega}, Theorems \ref{thm:PCS} and \ref{thm:real}.

When the signals are strong, we propose a feature selection step on $d =  \hat\Omega_1\hat{\mu}_1 - \hat{\Omega}_0 \hat{\mu}_0$ instead of $\hat{\mu}_1 - \hat{\mu}_0$. The inclusion of precision matrix in the feature selection step has been shown optimality in the linear classifier case where $\Omega_0 = \Omega_1 = \Omega$. In \cite{FJY}, it has been proved such innovated thresholding is better than the thresholding on $\hat{\mu}$ or $\Omega^{-1/2}\hat{\mu}$. When it comes to quadratic forms, we borrowed this idea. 

When the signals are weak, we suggest estimating $\hat{\mu}$ as a constant vector. It can be seen as a simple aggregation of all the features. The constants of it reduce the random error and hence achieve the optimal boundary; see Theorem \ref{thm:unknownmu}. 
In the supplementary material \cite[Section A]{QDAsupp}, we have explored the performance of the original QDA in the region where the signals in $\mu$ are weak. We have proved two theorems about the phase transition phenomenon of QDA/QDAfs when $\Omega_1$ is known or unknown. By QDA/QDAfs, there must be $\max\{\|\mu_0 - \mu_1\|^2, \|\Omega_1 - I\|^2\} \gg \sqrt{p/n}$ for successful classification when the signals in $\mu_0 - \mu_1$ are weak.
There is a gap between this upper bound and the statistical lower bound in Theorem \ref{thm:unknownmulower}. By QDAw, this gap will be overcome.

Finally, in Step 2 we find the constant $C$ by minimizing the training error for real data sets. Actually, in our detailed algorithms in Sections \ref{sec:main}, we define $C$ clearly based on the scenarios in concern. 


\subsection{Asymptotic rare and weak signal model}\label{sec:arw}
We propose a rare and weak signal model for both  mean vectors and precision matrices. 

The two classes have mean vectors $\mu_0 \in \mathcal{R}^p$ and $\mu_1 \in \mathcal{R}^p$. With a location shifting of the distance $\frac{1}{2}|\mu_1-\mu_0|$, we can take $\mu_1 = \mu$ and $\mu_0 = -\mu$. 
Hence, the signals in $\mu_i$ are the non-zeros in $\mu$.
We model $\mu$ as
\be\label{Param1}
\mu_i \ \overset{i.i.d.}\sim \ (1-\epsilon_p)\mathcal M_0+\epsilon_p \mathcal H, \ \ \ \ i=1,\dots,p,
\ee
where $\mathcal M_0$ is the point mass at 0 and $\mathcal H$ is a distribution that concentrates at $\tau_p$ and has no point mass at 0. 
Hence, the density of the signals can be captured by $\epsilon_p$ and the strength can be captured by $\tau_p$. To model the sparsity and weakness, we assume that when $p \goto \infty$, 
\be
\epsilon_p \ \rightarrow \ 0, \ \ \ \ \tau_p \ \rightarrow \ 0. 
\ee

When it comes to delicate theoretical analysis, we assume the signals have the same signs and strengths, which means $\mathcal{H} = \mathcal{M}_\tau$, the point mass at $\tau$. Such assumption is generally used in high-dimensional applications to facilitate the theoretical analysis; see \cite{DJ08, FJY, JKW}. 

The covariance matrix of class $k$ is $\Sigma_k$. Say for each feature $j$, its conditional variances given $Y = 0$ and $Y=1$ share the same main term; otherwise the signal in the variances are strong enough. Let $D_0$ be the diagonal matrix where the diagonals are the variances of features in Class 0. We normalize $\Sigma_k$ by $D_0^{-1/2} \Sigma_k D_0^{-1/2}$ so that $\Sigma_0$ have diagonals as 1 and $\Sigma_1$ have diagonals as $1+o(1)$. So we suppose the diagonals of $\Omega_0$ and $\Omega_1$ are around 1 without loss of generality. Let $D^{(k)}_{\Omega} = Diag(\Omega_k)$. We model it as follows:
\be\label{Omegadiag}
D^{(k)}_{\Omega}(i,i)  \stackrel{i.i.d.}{\sim} 1+\mathcal{D}_p, \quad 
1\leq i \leq p, \ \ k = 0, 1.
\ee
Here, $\mathcal{D}_p$ is a distribution with the magnitude concentrating at $\xi_p$. 

In many applications, $\Omega_k$'s, instead of $\Sigma_k$'s, have comparatively small number of non-zero entries in each row. We model the non-zeros on the off-diagonals as $V$, where 
\be\label{Omega0}
V^{(k)}_{ij} = V^{(k)}_{ij} \ \overset{i.i.d.}\sim \ (1-\nu_p)\mathcal M_0+\frac{\nu_p}{2} \mathcal M_{\eta_p} + \frac{\nu}{2} \mathcal M_{-\eta_p}, \ \ k = 0, 1,\ \  1 \leq i < j \leq p,
\ee
where $\mathcal M_{\eta_p}$ and $\mathcal M_{-\eta_p}$ are the point mass at $\eta_p$ and $-\eta_p$ respectively. The analysis still holds when $M_{\eta_p}$ and $M_{-\eta_p}$ are replaced by a symmetric distribution that concentrates on $\eta_p$ and $-\eta_p$ with no point mass at $0$. The information on diagonals and off-diagonals are modelled by two separate parameters because they have different effects on the clustering results; see Theorems \ref{thm:PCS} and \ref{thm:real}. 

Combine the modelling on the diagonals and off-diagonals, the precision matrices follow
\be\label{Omega1}
\Omega_k \ = \ \Sigma_k^{-1} \ = \ D^{(k)}_{\Omega} + V^{(k)}, \quad 
k = 0, 1.
\ee
To model the sparsity and weakness of signals in $\Omega_0$ and $\Omega_1$, we assume 
\be\label{Omegapara}
\nu_p \rightarrow 0, \qquad \eta_p \rightarrow 0, \qquad\xi_p \goto 0.
\ee

In our analysis, we consider the case when $\Omega_0$ is known and $\Omega_0$ is unknown. In the former case, we can set $X = \Omega_0^{1/2}X$ and update $\Omega_1$ and $\mu$ accordingly. Because of the sparsity in $\Omega_0$, the updated $\Omega_1$ and $\mu$ are still sparse. 
Hence, without loss of generality, we assume $\Omega_0=I$ and assume $\Omega_1$ follows (\ref{Omega0}) and (\ref{Omega1}). The model that satisfies (\ref{Param1}) -- (\ref{Omegapara}) is called rare and weak model.

One aim of this paper is to derive the statistical limits for the classification problem on the phase diagram. To derive it, we should tie all the parameters to $p$ by some constant parameters. The sample size $n$ goes to infinity at a slower rate than $p$, so we tie the sample size $n$ to $p$ by
\be\label{Param3}
n = n_p = p^{\delta}, \qquad 0 < \delta < 1. 
\ee
For $\mu$, we define the signal sparsity parameter $\epsilon$ and the weakness parameter $\tau$ as
\be\label{Param2}
\epsilon_p \ = \ p^{-\zeta}, \ \ \ \ \tau_p \ = p^{-\theta}, \ \ \ \ 0<\zeta, \theta <1.
\ee
For the precision matrices $\Omega_i$, $i = 0,1$, we similarly define the parameters as 
\be\label{Param}
\eta_p= p^{-\alpha}, \ \ \ \nu_p=p^{-\beta}, \ \ \ \xi_p=p^{-\gamma}, \ \ \ \ \ 0<\alpha, \gamma<1,\ \  0< \beta <2.
\ee
Here, $\alpha$, $\beta$, $\gamma$, $\zeta$, $\theta$, and $\delta$ are all constants. The regions of interest can be interpreted as the regions on the space formed by these parameters. 

Finally, to guarantee that $\Omega_k = D_{\Omega}^{(k)} + V^{(k)}$ is a positive definite matrix, $V^{(k)}$ must be weak enough so that $\Omega_k$ is positive definite. According to Lemma \ref{lemma1}, this requirement is satisfied with high probability under the condition
\be\label{con1}
\beta > 1 - 2\alpha. 
\ee
Hence, we discuss the regions under this condition only. 
The model that satisfies (\ref{Param1}) -- (\ref{con1}) is called asymptotic rare and weak model for classification.

\subsection{Phase transitions}\label{sec:mean} 
Under the ARW model, we analyze the regions of possibility and impossibility for any classifiers. 
In detail, we calibrate the impact of quadratic terms on classification in the following terms:
\begin{itemize}
\item the possibility and impossibility regions for the classification problem under the ideal case and the case that $\Omega_0$ and $\Omega_1$ are known. 
\item the possibility and impossibility regions for the classification problem when $\Omega_0$ and $\Omega_1$ are unknown but some sparsity conditions of them are satisfied. 
\end{itemize}
We summarize our results about the first part here. When all parameters are unknown, we present the conditions and results by Theorems \ref{thm:PCS} and \ref{thm:real} in Section \ref{sec:main}. 

Define a function 
\be\label{eqn:rho}
\rho_{\delta}(\zeta) = \left\{\begin{array}{ll}
1/2 - \zeta, & 0 < \zeta \leq  (1 - \delta)/2, \\
\delta/2, & (1-\delta)/2 < \zeta \leq  1 - \delta, \\
(1 - \zeta)/2, & 1 - \delta < \zeta < 1. 
\end{array}
\right.
\ee
Such a function can be found in multiple works about high-dimensional problems in the analysis of lower bounds; see \cite{ImpossibilityClass, JKW}. In different settings, the meaning of this function is different. 

\begin{thm}\label{thm:unknownmulower}[Lower bound]
Under the ARW model with $\Omega_0 = I$, if $\|\Omega_1 - I\|_F^2 \goto 0$ and $\theta > \rho_{\delta}(\zeta)$, then for any classifier $L$, when $p \to \infty$, there is 
\[
\mbox{Mis-classification Rate of $L$} \geq 1/2.
\]
\end{thm}

Let QDAw be the QDA classifier with $\hat{\mu} = a*{\bf 1}$ and $g(X) = g_w(X)$ as the estimation in weak signal case; and let QDAfs be the QDA classifier with $g(X) = g_s(X; \hat{\mu}, \Omega_1, I)$ as the strong signal case. Details of the two algorithms in Algorithms \ref{tab:algweak} and \ref{tab:algstrong}. 
They give the matching upper bound as the following theorem. 
\begin{thm}\label{thm:unknownmu}[Upper Bound]
Under the ARW model with $\Omega_0 = I$,
\begin{itemize}
\item [(i)] the mis-classification rate of $QDAw$ with $a = p^c$ for an arbitrary constant $0 < c < 1$ goes to 0 when $p \goto \infty$, if $\theta \geq \delta/2$ and one of the following conditions hold 
\begin{itemize}
    \item[(a)] $\|\Omega_1 - I\|_F^2 \gg p^c$; or, 
    \item[(b)] $\theta < \rho_{\delta}(\zeta)$;
\end{itemize}
\item[(ii)] the mis-classification rate of $QDAfs$ goes to 0 when $p \goto \infty$, if $\theta < \delta/2$ and one of the following conditions hold 
\begin{itemize}
    \item[(a)] $\|\Omega_1 - I\|_F^2 \to \infty$; or, 
    \item[(b)] $\theta < \rho_{\delta}(\zeta)$.
\end{itemize}
\end{itemize}
\end{thm}
For QDAw, the constant $c$ can be chosen arbitrarily. When $\|\Omega_1 - I\| \goto \infty$, i.e. $\max\{1 - 2\gamma, 2 - 2\alpha - \beta\} > 0$, we can always choose $c = \frac{1}{2}\max\{1 - 2\gamma, 2 - 2\alpha - \beta\}$ so that the inequality holds. Hence, the two boundaries match.

Theorems \ref{thm:unknownmulower}--\ref{thm:unknownmu} show the determining factor of the classification problem contains two parts, $\|\Omega_1 - I\|_F^2$ and $\mu$. Since we discuss the case $\Omega_0 = I$ here, $\|\Omega_1 - I\|_F^2$ can be regarded as $\|\Omega_1 - \Omega_0\|_F^2$, which is the effect of quadratic term. For $\mu$, we have to consider two cases: 
\begin{itemize}
\item[(i)] When $\theta < \delta/2$, the sample size is large enough so that the signals in the mean vector can be almost perfectly recovered. With the feature selection step, the QDAfs achieves an asymptotic mis-classification rate of 0 when $\max\{\|\Omega_1 - I\|_F^2, \|\mu\|^2\} \goto \infty$, and $1/2$ otherwise. In addition, the latter region is proven to be a failure region for all classifiers, which is referred to as the region of impossibility. 

\item[(ii)] When $\theta > \delta/2$, the sample size is insufficient for the signal recovery, and the feature selection step is ineffective. 
We use $g_w(X)$ to aggregate all the information in $\mu$ to do estimation, which is $O(\|\mu\|_1)$. 
The QDAw mis-classification rate converges to 0 when $\max\{\|\Omega_1 - I\|_F^2, \|\mu\|_1\} \goto \infty$, and $1/2$ otherwise. Again, the latter region is proven to be a failure region for all classifiers. 
\end{itemize}

\begin{figure}[htb!]
    \centering
    \subfigure[Phase transition on $\Omega_1$]{\includegraphics[width=0.48\textwidth]{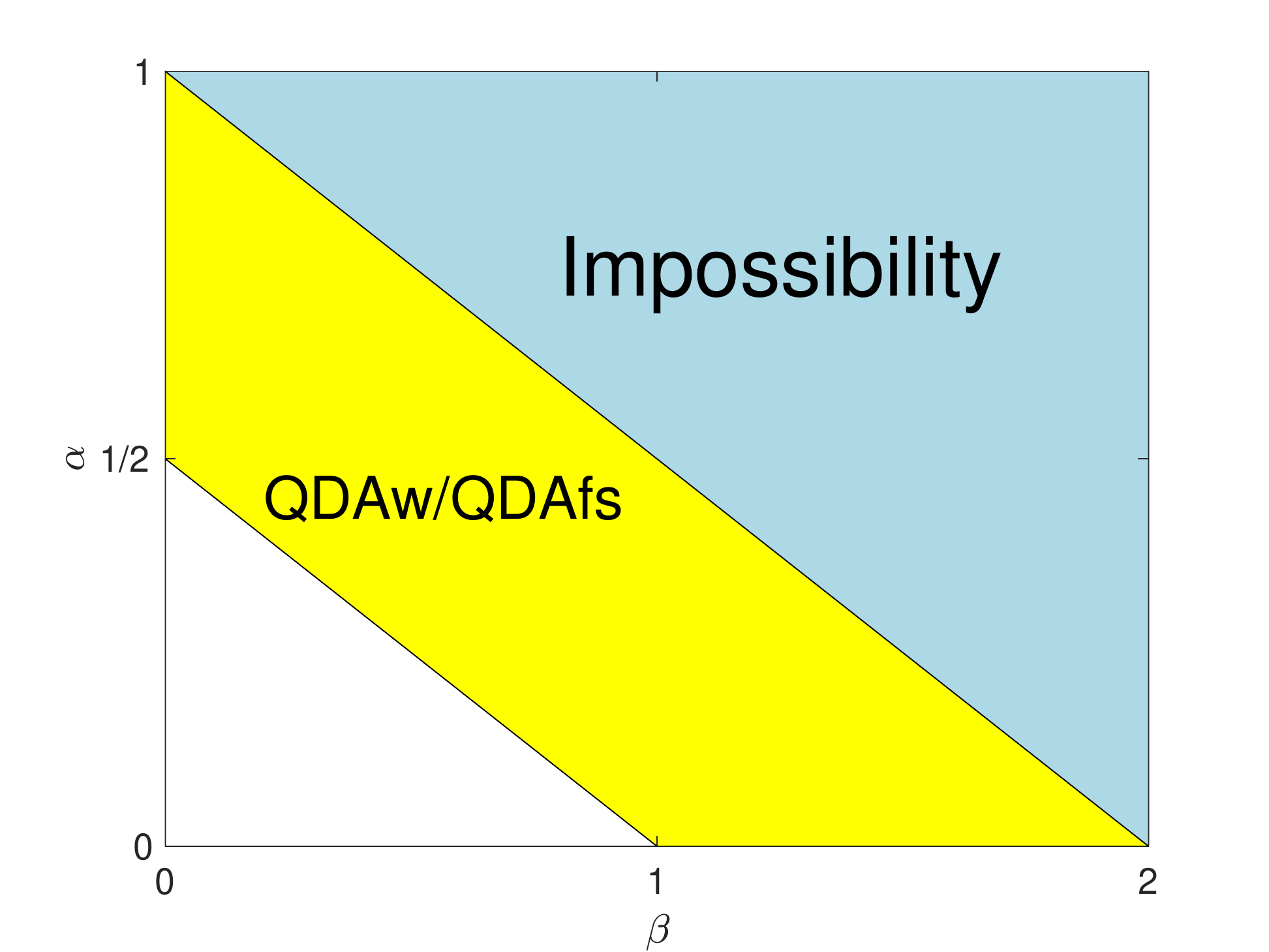}} 
               \subfigure[Phase transition on $\mu$]{\includegraphics[width=0.48\textwidth]{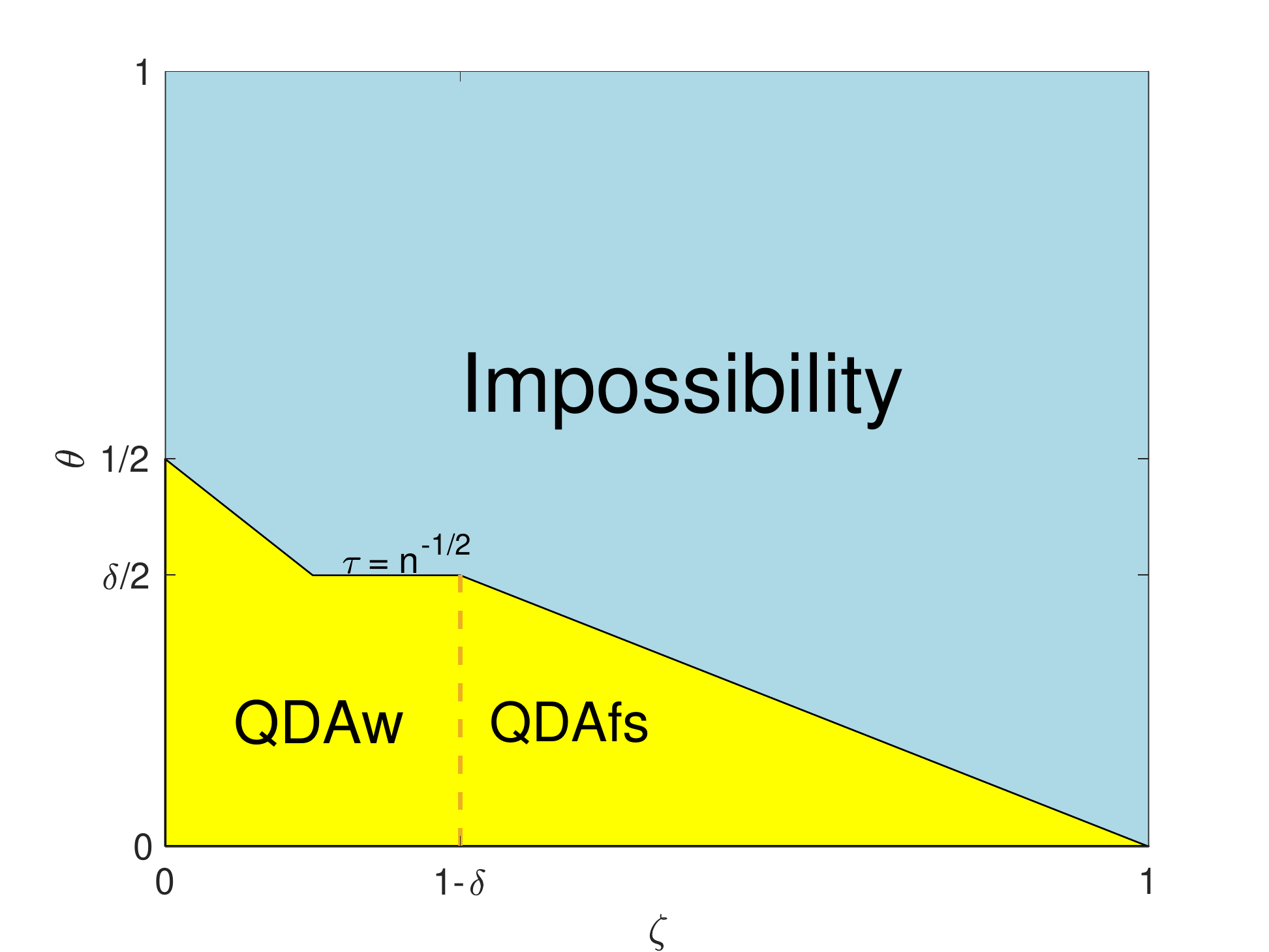}} 
      \caption{The possibility/impossibility regions derived in Theorems \ref{thm:unknownmulower} and \ref{thm:unknownmu} when $\delta$ and part of the rest parameters are fixed: (a) $\delta$, $\theta$ and $\zeta$ are fixed, $\gamma > 1/2$ and $\theta < \rho_{\delta}(\zeta)$; (b) $\delta$, $\alpha$ and $\beta$ are fixed, $\gamma > 1/2$ and $2 - 2\alpha - \beta < 0$.}
    \label{figure: parameter}
\end{figure}

Figure \ref{figure: parameter} provides a sense of the relationship between the sparsity and weakness parameters of the mean and covariance matrix. 
Subfigure (a) is on the $\alpha$-$\beta$ plane to present results about $\|\Omega_1 - I\|_F^2$. Subfigure (b) is on the $\theta$-$\zeta$ plane to present the regions on $\|\mu\|^2$. To see the effects clearly, we assume the information from the other part is insufficient for each subfigure. We can see QDAw and QDAfs are the optimal methods. 

In Subfigure (a), we suggest QDAw and QDAfs in the region of possibility, instead of only one method. The reason is although the contribution from $\Omega_1$ and $\mu$ seems independent of each other, the performance of QDAfs still relies on the signal strength in $\mu$. When the signals in $\mu$ cannot be successfully recovered, QDAfs requires $\|\Omega_1 - I\|_F^2 \gg \sqrt{p/n}$ to success, which cannot achieve the bound $\|\Omega_1 - I\|_F^2 \gg 0$ by QDAw. It comes from the effects of $\mu$ on $X^\top \Omega_{\rm diff} X$. This effect is rarely discussed in previous literature.

The methods employed and results obtained in this work are unique compared with other literature on QDA methods for high dimensional data with sparse signals (\cite{li2015sparse, fan2015innovated, wu2019quadratic}). We propose QDAw and QDAfs for different types of mean vectors and show they match the statistical lower bound, which is rarely discussed in other works.  

\subsection{A Real Data Example}\label{sec:quickdata}
We use a quick example to demonstrate how this works on the real data. We consider the rats dataset with summaries given in Table \ref{table:data}. This dataset consists of $181$ samples measured on the same set of $8491$ genes, with $61$ samples labeled by \cite{liver} as toxicants and the other $120$ as  drugs.  
The original rats dataset was collected in a study of gene expressions of live rats in response to different drugs and a toxicant; we use the cleaned version by \cite{liver}. 
\begin{table}[ht]
\caption{A gene-expression microarray rats dataset.}
\begin{tabular}{|l|l||c|c|}
\hline
Data Name & Source  &    $n$ ($\#$ of subjects)  & $p$ ($\#$ of genes)   \\
\hline
Rats & Yousefi et al. (2010)    &  181 &  8491  \\
\hline
\end{tabular}
\label{table:data}
\end{table}

This dataset has been carefully studied in \cite{PCS}, with the performance of the two-class classification compared among a sequence of popular classifiers, including SVM in \cite{SVM}, Random Forest in \cite{RF}, and HCT-PCS. The HCT-PCS, which achieves optimal classification when it adapts  LDA \citep{DJ08, FJY} in the rare and weak signal setting, was shown to have very promising classification results with this data. 

That said, in HCT-PCS, all samples of the two classes are assumed to share the same precision matrix, leaving room for improvement. We now apply QDAfs with data normalization (details in Table \ref{tab:alg2}) to this data set and compare the results with those from the LDA with HCT-PCS approach. Here, we leave out all the implementation details, which will be introduced in Section \ref{sec:rats}, and only highlight our findings for this rats data:
\begin{itemize}
\item QDA further outperforms LDA with HCT-PCS, and produces better results than those other methods in \cite{PCS}, including SVM and Random Forest, suggesting that QDA gives a better separation by taking into account the second-order difference between the two classes.

\end{itemize}

\begin{figure}[htb!]
    \centering
    \subfigure[]{\includegraphics[scale=0.34]{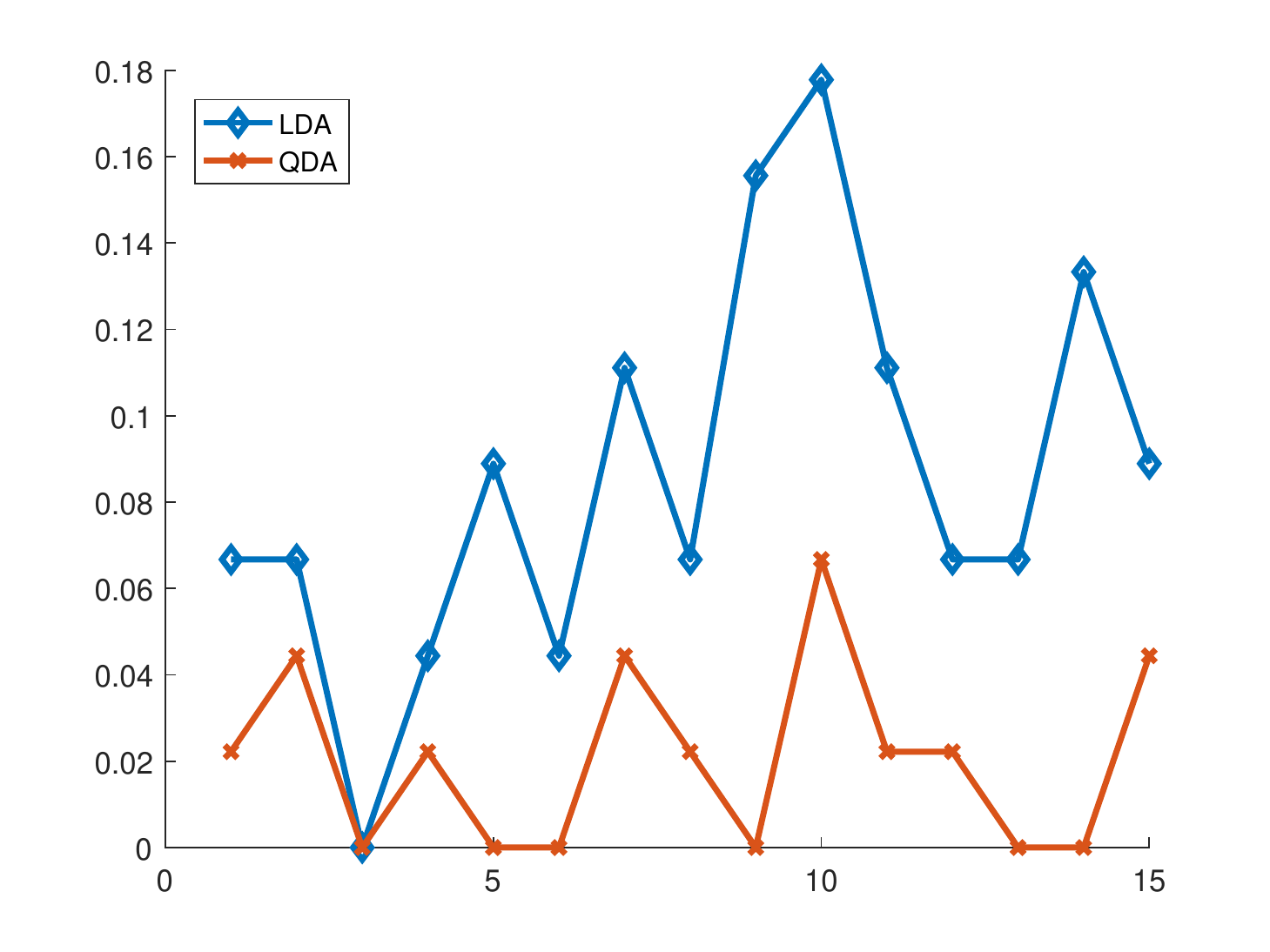}} \hspace{1cm}
	   \subfigure[]{\includegraphics[scale=0.45]{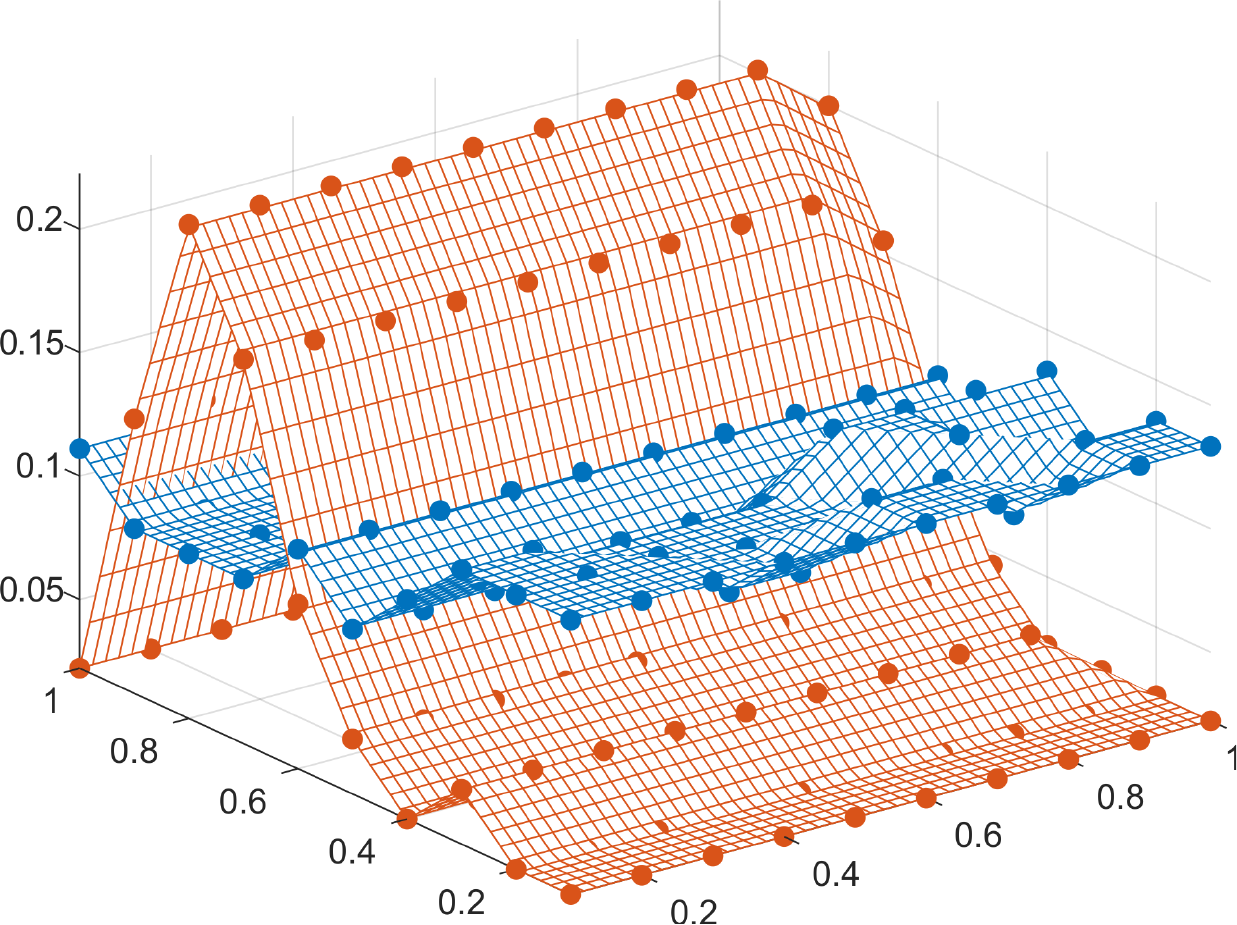}}
	   	   \caption{Comparison of testing errors for the rats data: (a)  error rates (y-axis) of LDA (blue) and QDA (red) for 15 data splittings (x-axis) at a certain sparsity of the precision matrix; (b) Zoom-in errors for the rats data for varying choices of parameters in estimating the precision matrices for one splitting of 15 splittings.}
    \label{figure:errorrats-k30-intro}
\end{figure}

We record 15 random splits of the rats data for the training data and test data. The test error is illustrated in Figure \ref{figure:errorrats-k30-intro} (below, left). We can see that the test error of LDA are all above those of QDA at every data splitting, given that all the tuning parameters are selected in the same way. Figure \ref{figure:errorrats-k30-intro} (below, right) demonstrates the surface of the test error between LDA and QDA, by varying the tuning parameters in the precision-matrix estimation. This Zoom-in plot shows that the QDA does bring necessary improvement over LDA when the precision matrices are appropriately estimated. 

\subsection{Content and Notations}
The main results for the phase transitions under various scenarios are discussed in Section \ref{sec:main}. 
Proofs of lower bounds are given in Section \ref{sec:lowerproof} and that of upper bounds are given in Section \ref{sec:proof}. In Section \ref{sec:rats}, we present numerical results of the proposed methods and algorithms on real data. In Section \ref{sec:dis}, some concluding remarks and potential directions of future work are discussed. The details of the proofs are provided in the supplementary materials. 

Here we list the notations used throughout the paper. Let the eigenvalues of $W$ be denoted by $\lambda_1 \geq \lambda_2 \geq \cdots \geq \lambda_p$. For a matrix $M$, we use $\|M\|$ and $\|M\|_F$ to denote its spectral normal and Frobenius norm, respectively, $det(M)$ to denote the determinant of $M$ and $Tr(M)$ to denote its trace which equals the summation of the eigenvalues of $M$. We use $diag(c_1,\dots, c_p)$ to denote a diagonal matrix with diagonal elements $c_1,\dots, c_p$, and use $I(A)$ to denote an indicator function over event $A$. For two vectors or matrices $a$ and $b$ of same dimension, $a\circ b$ denotes the Hadamard (entrywise) product.  


\section{Phase transition for the classification problem}\label{sec:main}
Throughout the whole paper, we consider the mixture model 
\be \label{M12}
X|Y \sim (1 - Y) N(-\mu, \Omega_0^{-1})+ Y N(\mu,\Omega_1^{-1}),
\ee
and the mis-classification rate as 
\be\label{eqn:error}
MR = [P_{\epsilon, \tau, \eta, \nu, \xi}(\hat Y = 0|Y = 1) + P_{\epsilon, \tau, \eta, \nu, \xi}(\hat Y = 1|Y = 0)]/2,
\ee
where $X$ is a fresh data vector with $Y$ as the true label and $\hat{Y}$ being the estimated label. We refer to the error rate  by a classifier $L$ as  $MR(L)$. 
Since we consider $P(Y = 0) = q = 1/2$, so $MR$ is the average of two types of errors. For a general $q$, we should update $MR$ as $MR = (1-q)P(\hat Y = 0|Y = 1) + q P(\hat Y = 1|Y = 0)$ and the results still hold. 


\subsection{New QDA approaches}\label{sec:newqda}
When all the parameters are known, the classifier is in (\ref{CR2}), where we calculate 
\be\label{eqn:classqda}
Q(X) = X^\top(\Omega_0 -\Omega_1)X+2\mu^\top(\Omega_0 + \Omega_1 )X + \mu^\top(\Omega_0-\Omega_1)\mu+\ln|\Omega_1|- \ln|\Omega_0|,
\ee
and estimate $\hat{Y} = I(Q(X) > 0)$. 

We have presented the estimation of $X^\top (\Omega_0 - \Omega_1)X$ and $2\mu^\top(\Omega_0 - \Omega_1 )X$ for the unknown parameter case in Section \ref{sec:qda}. Here we define two functions to better present the estimates in the algorithms. 
For a symmetric matrix $M$ and threshold $t$, define  
\be\label{eqn:trunm}
T(M;t) = M - diag(diag(M)\circ I\{|diag(M)| <= t\}).
\ee
By $T(M;t)$, the diagonals of $M$ that are smaller than $t$ are truncated to be 0 and the off-diagonals do not change. 
Given matrices $A$, $B$, vector $d$ and threshold $t$, define a vector $d^{(t)}$ with $d^{(t)}(j) = I\{|d(j)| \geq t\}$ and the functions
\be\label{eqn:linear}
g_w(X; \hat{\mu}, A, B) = 
\hat{\mu}^\top(A + B)X, 
\quad
g_s(X; d, t) = (d\circ d^{(t)})^\top X.
\ee

With all the preparations, we present the algorithms for various scenarios. Tables \ref{tab:algweak} and \ref{tab:algstrong} are the algorithms employed in Theorem \ref{thm:unknownmu}, for the case that both $\Omega_i$'s are known. The algorithm for that $\Omega_0$ is known and $\Omega_1$ is unknown is in Table \ref{tab:QDAomega1}, and the most general case is in Table \ref{tab:QDAgeneral}. 
For previous cases, the constants are clearly stated. 

\begin{table}[ht!]
\caption{Algorithm QDAw: weak and relatively dense signals in $\mu$, $\Omega_1$ is known, $\Omega_0 = I$.} 
\begin{tabular}{ll}\\ \hline
& \underline{Input}: data points $(X_i, Y_i)$, $1 \leq i \leq n$;  constant $0 < c < 1$; new data point $X$; true precision matrix $\Omega_1$. \\
& \underline{Output}: label $\hat{Y}$. \\
1. & \underline{Parameter Estimation}:  Let $\hat{\mu}_0 = \frac{1}{n_0} \sum_{i: Y_i = 0} X_i$, where $n_0 = n - n_1$.  Let $\hat{\mu} = p^{(c-1)/2}*{\bf 1}$, \\
&where ${\bf 1} \in \mathcal{R}^{p}$ is the vector of ones.\\
2.  & Let $C = \hat{\mu}_0^\top (\Omega_1 - I)\hat{\mu}_0 + \ln |{\Omega_1}| + \frac{1}{n_0}Tr({\Omega_1} - I)$. \\
3.  & \underline{QDA Score}: Calculate the QDA score $Q=X^\top(I - {\Omega_1}) X + 2g_w(X; \hat{\mu}, I, \Omega_1) + C$. \\
4. &\underline{Prediction}: Predict $\hat Y = I\{Q > 0\}$.  \\
 \hline
\end{tabular}
\label{tab:algweak} 
\end{table}

\begin{table}[ht!]
\caption{Algorithm QDAfs: moderately strong and sparse signals in $\mu$, $\Omega_1$ is known, $\Omega_0 = I$.} 
\scalebox{1}{ 
{\begin{tabular}{ll}\\ \hline
& \underline{Input}: data points $(X_i, Y_i)$, $1 \leq i \leq n$;  threshold $t > 0$; new data point $X$. true precision matrix $\Omega_1$. \\
& \underline{Output}: label $\hat{Y}$. \\
1. &\underline{Parameter Estimation}: Let $\hat{\mu}_0 = \frac{1}{n_0} \sum_{i: Y_i = 0} X_i$ and $\hat{\mu}_1 = \frac{1}{n_1} \sum_{i: Y_i = 1} X_i$, where \\
&$n_1 = \sum_{i =1}^n Y_i$ and $n_0 = n - n_1$.   \\ 
2.  & Let $d_0=\hat{\mu}_0$, $d_1= {\Omega_1} \hat{\mu}_1$ and $d = d_1 - d_0=(d(1),\dots, d(p))^\top$.  \\
3.  & \underline{Thresholding}: Let $d^{(t)}$ denote the indicator vector of feature selection, i.e. $d^{(t)}(j) = 1\{|d(j)| \geq t\}$, \\
&for $j=1,\dots,p$.\\
&Let $\Omega_1^{(d)}$ be the sub-matrix of $\Omega_1$ constrained on rows and columns that $d^{(t)} = 1$.\\
4.  & Let $C = (\hat{\mu}_0 \circ d^{(t)})^\top (I - {\Omega_1})(\hat{\mu}_0 \circ d^{(t)}) + \ln |{\Omega_1}| + \frac{1}{n_0}Tr(\Omega_1^{(d)} - I)$. \\
5.  & \underline{QDA Score}: Calculate the QDA score $Q=X^\top(I - {\Omega_1}) X + 2g_s(X; d, t) + C$. \\
6. &\underline{Prediction}: Predict $\hat Y = I\{Q > 0\}$.  \\
 \hline
\end{tabular}}
} 
\label{tab:algstrong} 
\end{table}

\begin{table}[ht!]
\caption{Algorithm QDAfs/QDAw-PCS: $\Omega_0 = I$.} 
\scalebox{1}{ 
{\begin{tabular}{ll}\\ \hline
& \underline{Input}: data points $(X_i, Y_i)$, $1 \leq i \leq n$;  threshold $t > 0$; new data point $X$. \\  
& \underline{Output}: label $\hat{Y}$. \\
1. &\underline{Parameter Estimation}: Let $\hat{\Omega}_1$ be the estimation from PCS. \\
2. & \underline{Thresholding}: Update $\hat{\Omega}_1$: $\hat{\Omega}_1 = T(\hat{\Omega}_1 - I; \sqrt{2\ln p/n}) + I$.\\
3. & Apply QDAw in Table \ref{tab:algweak} or QDAfs in Table \ref{tab:algstrong} with the input precision matrix as $\hat{\Omega}_1$. \\
 \hline
\end{tabular}}
} 
\label{tab:QDAomega1} 
\end{table}

\begin{table}[ht!]
\caption{Algorithm QDAfs-PCS: all unknown.} 
\scalebox{1}{ 
{\begin{tabular}{ll}\\ \hline
& \underline{Input}: data points $(X_i, Y_i)$, $1 \leq i \leq n$;  threshold $t > 0$; new data point $X$. \\  
& \underline{Output}: label $\hat{Y}$. \\
1. &\underline{Parameter Estimation}: Let $\hat{\Omega}_0$ and $\hat{\Omega}_1$ be the estimation from PCS. \\
2. & \underline{Thresholding}: Let $\hat{\Omega}_{\rm diff} = T(\hat{\Omega}_0 - \hat{\Omega}_1; 2\sqrt{2\ln p/n})$.\\
3. &\underline{Parameter Estimation}: Let $\hat{\mu}_0 = \frac{1}{n_0} \sum_{i: Y_i = 0} X_i$ and $\hat{\mu}_1 = \frac{1}{n_1} \sum_{i: Y_i = 1} X_i$, where \\
&$n_1 = \sum_{i =1}^n Y_i$ and $n_0 = n - n_1$.   \\ 
4.  & Let $d_0=\hat{\Omega}_0 \hat{\mu}_0$, $d_1= {\hat\Omega_1} \hat{\mu}_1$ and $d = d_1 - d_0=(d(1),\dots, d(p))^\top$.  \\
5. & Let $C = (\hat{\mu}_0 \circ d^{(t)})^\top \hat{\Omega}_{\rm diff} (\hat{\mu}_0 \circ d^{(t)}) + \ln |\hat{\Omega}_0 - diag(\hat{\Omega}_0) + I| -  \ln |\hat{\Omega}_1 - diag(\hat{\Omega}_1) + I|$.\\
6.  & \underline{QDA Score}: Calculate the QDA score $Q(X)=X^\top\hat{\Omega}_{\rm diff} X + 2g_s(X; d, t) + C$. \\
7. &\underline{Prediction}: Predict $\hat Y = I\{Q(X) > 0\}$.  \\
 \hline
\end{tabular}}
} 
\label{tab:QDAgeneral} 
\end{table}

\subsection{Ideal case}\label{sec:ideal}
When all the parameters are known, the classical QDA classifier provides the optimal results; see Proposition \ref{thm:idealunequal}. 

\begin{prop}\label{thm:idealunequal}[Phase transition for the ideal case.]
Consider the rare and weak signal model (\ref{Param1}) -- (\ref{Omegapara}). 
\begin{itemize}
\item[(i)] The QDA classifier (\ref{eqn:classqda}) has a mis-classification rate $MR(QDA) \goto 0$ as $p\to\infty$, if $\|\Omega_1 - \Omega_0\|_F^2 + 8\|\mu\|^2 \goto \infty$. 

If further (\ref{Param3})--(\ref{con1}) are satisfied, then $MR(QDA) \goto 0$ if one of the following conditions is satisfied, 
\begin{itemize}
\item [(1)]  $\beta<2-2\alpha$; or 
\item [(2)] $\gamma < 1/2$; or
\item [(3)] $\zeta<1-2\theta$.
\end{itemize}
\item[(ii)] The mis-classification rate $MR(L)$ of any classifier $L$ converges to $1/2$ when $p\to\infty$, 
if $\|\Omega_1 - \Omega_0\|_F^2 + 8\|\mu\|^2 \goto 0$. 
\end{itemize}
\end{prop}

{\bf Remark 1}. Proposition \ref{thm:idealunequal} describes an exact phase diagram of the classification problem. When $\|\Omega_1 - I\|_F^2 + 8\|\mu\|^2 \goto \infty$, the QDA method achieves a mis-classification rate of 0 asymptotically. On the complement region, all classifiers fail. 
In this sense, Proposition \ref{thm:idealunequal} demonstrates QDA succeeds in the whole possibility region and thus is optimal.

{\bf Remark 2}. It can be observed the contribution of $\mu$ is $\|\mu\|^2$ and the contribution of $\Omega_1$ is $\|\Omega_1 - I\|_F^2$. 
There is no intersection between them because all the parameters are known. 
When we consider data-trained classifiers, the interaction may happen to be the choice of algorithms (see Theorem \ref{thm:unknownmu}), or 
the two possibility regions will depend on the parameter from the other part (see Theorems \ref{thm:PCS} and \ref{thm:real}).

{\bf Remark 3}. Since the contribution of $\Omega_1$ is $\|\Omega_1 - I\|_F^2$, the diagonals and off-diagonals perform independently. 
The off-diagonals of $\Omega - I$ is modelled in (\ref{Omega0}) and (\ref{Param}) by $\alpha$ and $\beta$ to measure the signal strength and sparsity. The diagonals only have one signal strength parameter $\xi = p^{-\gamma}$. Hence, condition (1) is an inequality between $\alpha$ and $\beta$ while condition (2) is about $\gamma$ solely.  


We do not consider the sparsity in the diagonals of $\Omega - I$ due to model complexity. In that sense, the model thus raised will have 6 sparsity and weakness indices in total: 2 for means, 2 for diagonals, and 2 for off-diagonals of $\Omega_1$. We can readily obtain the possibility region for this model, but it will be difficult to visualize and is thus omitted here.

\subsection{Phase transitions with partial information}\label{sec:unknownpre}
With Theorems \ref{thm:unknownmulower} and \ref{thm:unknownmu} in Section \ref{sec:mean}, we have discussed the phase transition phenomenon when $\mu$ is unknown but $\Omega_0$ and $\Omega_1$ are known. 
Here we further show the results when $\Omega_0 = I$ without loss of generality, and $\Omega_1$ is unknown. 

When $\Omega_1$ is unknown, the first problem is to estimate it.  
It has been discussed in numerous publications in the literature, such as \cite{CLIME, FJY, glasso, PCS}. However, estimation of high-dimensional precision matrix is restricted to the case that $\Omega_1$ is sparse.
Here, we consider the case that the signals in $\Omega_1$ are sparse and strong, that  
\begin{equation}\label{eqn:Param4}
    \alpha < \delta/2, \quad 1-\delta/2 < \beta < 2.
\end{equation}
Under (\ref{eqn:Param4}), with high probability, the number of non-zero entries in each row of $\Omega_1$ is $o(n)$ and the signal strength $\eta \gg 1/\sqrt{n}$. 

Under (\ref{eqn:Param4}), we suggest to estimate $\Omega_1$ by the Partial Correlation Screening (PCS) approach in \cite{PCS}. The PCS approach has a good control on $\|\hat{\Omega}_1 - \Omega_1\|_{\max}$, and hence $\|\hat{\Omega}_1 - \Omega_1\|_F^2$ is under control. 
With the PCS estimation, we further apply a truncation on the diagonals of $\hat{\Omega}_1$, where we assign $\hat{\Omega}_1(i,i)$ to be 1 if it is close to 1, i.e. $|\hat{\Omega}_1(i,i) - 1| \leq \sqrt{2\ln p/n}$). 
The truncation step helps to remove the noise on diagonals, but suffers a loss on the diagram of $\gamma$. 

With $\hat{\Omega}_1$, we apply QDAw and QDAfs to estimate $Y$; details in Table \ref{tab:QDAomega1}. We call it as QDAw-PCS or QDAfs-PCS. 
To find the optimality of it, we first find the phase diagram of the classification problem when $\mu$ is given; see the following proposition. 
When $\mu$ is given, we do not need QDAw or QDAfs. The classical QDA in (\ref{eqn:classqda}) with $\hat{\Omega}_1$, called QDA-PCS, will work. 

\begin{prop}\label{prop:unknownomega}
Consider model (\ref{M12}) with the parameterizations (\ref{Param1})--(\ref{con1}) and (\ref{eqn:Param4}). Suppose $\mu$ is given and $\Omega_0 = I$.

\begin{itemize}
\item [(i)] the QDA classification rule (\ref{eqn:classqda}) with $\Omega_1 = \hat{\Omega}_1$ as the truncated PCS estimate of $\Omega_1$ has a mis-classification rate $MR(QDA$-$PCS) \goto 0$ as $p\to\infty$, if 
\begin{itemize}
\item [(i)]  $\gamma < \delta/2$; or
\item [(ii)]  $\gamma > (2\alpha + \beta - 1)/2$ and $2\alpha + \beta - 2 < 0$; or 
\item [(iii)] $\gamma > \frac{1}{2}\min\{1, 2\alpha + \beta - 1\}$ and $2\theta + \zeta - 1 < 0$.
\end{itemize}
\item[(ii)]  $MR(L) \geq 1/2$ for any classifier $L$ when $p \goto \infty$, if $\gamma > 1/2$, $ 2\alpha + \beta - 2 > 0$ and $2\theta + \zeta - 1> 0$.
\end{itemize}
\end{prop}

When $\gamma > 1/2$, then QDA-PCS achieves the optimal boundary. When $\gamma < \delta/2$, then QDA-PCS also provides satisfactory classification results. However, when $2\alpha + \beta - 1 > \delta$, then on $\delta/2 < \gamma < (2\alpha + \beta - 1)/2$, QDA-PCS suffers a power loss because of the truncation on the diagonals. 

Now we consider the case $\mu$ and $\Omega_1$ are both unknown. Similar as QDA-PCS in Proposition \ref{prop:unknownomega}, here we find QDAw-PCS and QDAfs-PCS can match the lower bound.  
\begin{thm}\label{thm:PCS}
Consider model (\ref{M12}) with the parameterizations (\ref{Param1})--(\ref{con1}) and (\ref{eqn:Param4}). Suppose $\Omega_0 = I$.

\begin{itemize}
\item [(i)] MR(QDAw-PCS) with arbitrary constant $c$ goes to 0 when $p \goto \infty$, if $\theta \geq \delta/2$ and one of the following conditions hold 
\begin{itemize}
\item [(a)]  $\gamma < \delta/2$; or
\item [(b)]  $\gamma > (2\alpha + \beta - 1)/2$ and $2\alpha + \beta - 2 < c$; or 
\item [(c)] $\gamma > \frac{1}{2}\min\{1, 2\alpha + \beta - 1\}$ and $\theta < \rho_{\delta}(\zeta)$.
\end{itemize}
\item[(ii)] MR(QDAfs-PCS) goes to 0 when $p \goto \infty$, if $\theta < \delta/2$ and one of the following conditions hold 
\begin{itemize}
\item [(a)]  $\gamma < \delta/2$; or
\item [(b)]  $\gamma > (2\alpha + \beta - 1)/2$ and $2\alpha + \beta - 2 < 0$; or 
\item [(c)] $\gamma > \frac{1}{2}\min\{1, 2\alpha + \beta - 1\}$ and $\theta < \rho_{\delta}(\zeta)$.
\end{itemize}
\item[(iii)] $MR(L) \geq 1/2$ for any classifier $L$ when $p \goto \infty$, if $\gamma > 1/2$, $2\alpha + \beta -2 > 0$, and $\theta > \rho_{\delta}(\zeta)$ in (\ref{eqn:rho}).
\end{itemize}

\end{thm}

{\bf Remark 1}. The upper bound of QDAfs-PCS or QDAw-PCS matches the lower bound when $\gamma > \min\{1, 2\alpha + \beta - 1\}/2$ or $\gamma < \delta/2$. When the diagonals have parameter $\delta/2 < \gamma < \min\{1, 2\alpha + \beta - 1\}/2$, the random error is too large to recognize the truth and do successful classification. 

{\bf Remark 2}. Even for the case that $\theta < \rho_\delta(\zeta)$ and $\mu$ performs the main role in classification, we still need the condition that $\gamma > \frac{1}{2}\min\{1, 2\alpha + \beta - 1\}$ so that $\|\hat{\Omega}_1 - I\|_F^2$ is under control. This can be seen as the intervene between the quadratic and the linear terms. 

\subsection{Phase transitions with unknown parameters}\label{sec:real}
The most generalized case is that all the parameters are unknown. According to Theorem \ref{thm:PCS}, we consider both $\Omega_1$ and $\Omega_0$ have sparse and strong off-diagonal signals as (\ref{eqn:Param4}) and very weak diagonal signals that 
\be\label{eqn:gamma}
\gamma > 1/2. 
\ee
Under this condition, the diagonals are too weak to do successful classification. 

\begin{thm}\label{thm:real}
Consider the ARW model (\ref{M12}) with the parameterizations (\ref{Param1})--(\ref{con1}), (\ref{eqn:Param4}) and (\ref{eqn:gamma}). 
\begin{itemize}
\item[(i)] MR(QDAfs-PCS) in Table \ref{tab:QDAgeneral} goes to 0 when $p \goto \infty$, if $\theta < \delta/2$ and one of the following conditions hold 
\begin{itemize}
\item [(a)]  $2 - 2\alpha - \beta > 0$; or \item [(b)] $1 - 2\theta - \zeta > 0$.
\end{itemize}
\item[(ii)] $MR(L) \geq 1/2$ for any classifier $L$ when $p \goto \infty$, if $\theta < \delta/2$, $2 - 2\alpha - \beta < 0$ and $1-2\theta-\zeta < 0$.
\end{itemize}
\end{thm}

The theorem suggests QDAfs-PCS is optimal if $\Omega_i$ and $\mu$ have strong and sparse signals.


\section{Proof of lower bounds}
\label{sec:lowerproof}
We present the lower bound when $\mu$ and $\Omega_i$'s are all known in Proposition \ref{thm:idealunequal}, only $\mu$ and $\Omega_0$ are known in Theorem \ref{thm:unknownmulower}, only $\Omega_0$ is known in Theorem \ref{thm:PCS} and all are unknown in Theorem \ref{thm:real}. 
When the signal strength and sparsity parameters falls below the lower bound, any classifier $L$ will fail.
In this section we will prove these results.

\subsection{Proof of lower bound in Proposition \ref{thm:idealunequal}}\label{subsec:lowerideal}
In this ideal case both $\mu$ and $\Omega_i$'s are known. Let $f$ be the density function of $X \sim N(-\mu, \Omega_0^{-1})$ and $g$ be the density function of $X \sim N(\mu, \Omega_1^{-1})$. The Hellinger affinity between $f$ and $g$ is defined as $H(f, g) = \int \sqrt{f(x) g(x)} dx$.
\begin{lemma}\label{lemma:H} 
For any classifier $L = L(X|\mu, \Omega_0, \Omega_1)$, 
\[
|MR(L) - 1/2 | \leq C(1 - H(f, g))^{1/2}. 
\]
\end{lemma}

This lemma is well known, and so we omit the proof. According to this lemma, $H(f, g) = 1 + o(1)$ suffices to prove the impossibility. Introduce the normal density into $H(f, g)$, with basic calculations we have 
\be\label{eqn:Hfg}
H(f, g) = \exp\{-\frac{1}{2}\bigl[ \|\mu\|^2(1+o(1)) +  \|\Omega_0 - \Omega_1\|^2_F/8\bigr]\}. 
\ee
Therefore, when $\|\mu\|^2 + \|\Omega_0 - \Omega_1\|^2_F/8 = o(1)$, $H(f, g) = 1 + o(1)$ and the mis-classification error from any classifier will be close to 1/2. 

Under (\ref{Param1})--(\ref{Omegapara}), with high probability, $ \|\Omega_0 - \Omega_1\|_F^2 \leq 4(p\xi_p^2 + \eta_p^2p^2\nu_p)(1+o(1))$ and $\|\mu\|^2 = \tau_p^2 p \epsilon_p^2(1+o(1))$. 
 In the region of impossibility, both terms converge to 0, and then $H(f, g) = 1 + o(1)$. As a result, $MR(L) \goto 1/2$ for any classifier $L$. The lower bound in Theorem \ref{thm:idealunequal} is proved. 

\subsection{Proof of Theorem \ref{thm:unknownmulower}}\label{subsec:lowermu}
 Here we consider the case $\mu$ and $\Omega_1$ are unknown, and $\Omega_0 = I$ without loss of generality. Since we have to use training data to estimate $\mu$, the density functions are updated to be $f = f(X, X_1, \cdots, X_n; \Omega_1)$ and $g = g(X, X_1, \cdots, X_n; \Omega_1)$, where $X_i \sim N((2Y_i-1) \mu, I + Y_i(\Omega_1^{-1} - I))$ with known $Y_i$'s for both cases. The new data point is assumed to be $X \sim N(-\mu, I)$ for $f$ and $X\sim N(\mu, \Omega_1^{-1})$ for $g$. We want to prove $H(f,g) = 1 + o(1)$. 

Here, $f$ and $g$ differ at both mean and covariance matrix of $X$. We define $\tf$ to be a middle state, that $\tf = \tf(X, X_1, \cdots, X_n; \Omega_1)$, where $X \sim N(\mu, I)$ and others are the same. Hence, $\tf$ differs with $f$ only on the mean vector of $X$, and differs with $g$ only on the covariance matrix of $X$. 
When both both $\|f - \tf\|_1 = o(1)$ and $\|g - \tf\|_1 = o(1)$, there is $\|f - g\|_1 = o(1)$ and hence $H(f,g) = 1+o(1)$.

Consider $\|f - \tf\|_1$ first. It comes to the classification problem with an identity covariance matrix. In \cite{JKW}, it has been proved that $\|f - \tf\|_1 = C p\eps_p^2(e^{\tau_p^2} - 1)(1+o(1)) e^{n\tau_p^2}$ when $\zeta < 1 - \delta$, and $\|f - \tf\|_1 = C \sqrt{(e^{2p\epsilon_p \tau_p^2} - e^{-2p\epsilon_p \tau_p^2})/2}(1+o(1))$ when $\zeta > 1 - \delta$. Introducing (\ref{Param2}) that models $\epsilon_p$ and $\tau_p$, $\|f - \tf\|_1 = o(1)$ when one of the following can be satisfied: 
\begin{itemize}
    \item[(a)] $\theta \geq \delta/2$, $\zeta < 1 - \delta$, $\zeta + \theta < 1/2$; or 
    \item[(b)] $\zeta > 1 - \delta$, $ 1 - \zeta - 2\theta < 0$. 
\end{itemize}

Consider $\|g - \tf\|_1$ where $\mu$ is unknown and (\ref{Param1}) holds. With some calculations, 
\bea
\|g - \tf\|_1 & = & \int \int \frac{1}{(2\pi)^{p/2}} e^{-\frac{1}{2}(X - \mu)^\top (X - \mu)}|1 - det(\Omega_1)^{1/2} e^{-\frac{1}{2}(X - \mu)^\top (\Omega_1 - I)(X - \mu)}| dX dF(\mu)\nonumber\\
& = & \int  \frac{1}{(2\pi)^{p/2}} e^{-\frac{1}{2}X^\top X}|1 - det(\Omega_1)^{1/2} e^{-\frac{1}{2}X^\top (\Omega_1 - I)X}| dX. 
\eea
It equals to the $L_1$ distance between $f_s \sim N(0, I)$ and $g_s \sim N(0, \Omega_1^{-1})$. Therefore, to show $\|g - \tf\|_1 = o(1)$, it is to prove $\|f_s - g_s\| = o(1)$, which is equivalent with $H(f_s, g_s) = 1 + o(1)$. For $H(f_s, g_s)$, there is no training data and we can calculate the Hellinger distance directly, which is 
\be\label{eqn:prooflower2omega}
H(f_s, g_s) = \frac{det(\Omega_1)^{1/4}}{det((\Omega_1 + I)/2)^{1/2}} = \|\Omega_1 - I\|^2_F/8(1 + o(1)).
\ee
As a conclusion, $\|g - \tf\|_1 = o(1)$ when $\|\Omega_1 - I\|^2_F \goto 0$. 

Recall that $H(f,g) = 1+o(1)$ when both $\|f - \tf\|_1$ and $\|g - \tf\|_1$ are $o(1)$. 
Combine it with the results for $\|f - \tf\|_1$ and $\|g - \tf\|_1$. Therefore, $H(f,g) = 1+o(1)$ when $\|\Omega - I\|_F^2 \goto 0$ and one of the following conditions can be satisfied: 
\begin{itemize}
    \item[(a)] $\theta \geq \delta/2$, $\zeta < 1 - \delta$, $1 - 2\zeta - 2\theta < 0$; or 
    \item[(b)] $\zeta > 1 - \delta$, $ 1 - \zeta - 2\theta < 0$. 
\end{itemize}
Consider condition (a), $\zeta < 1-\delta$ always holds when $\theta \geq \delta/2$ and $\zeta + \theta < 1/2$, so the condition $\zeta < 1-\delta$ can be removed. Actually, when $\theta \geq \delta/2$, the region of impossibility will be decided by condition (a) because $1 - \zeta - 2\theta < 0$ in (b) always indicate $1 - 2\zeta - 2\theta < 0$ in (a). So we only need to consider the case $\theta < \delta/2$ for condition (b). When $\theta < \delta/2$, $1 - \zeta < 2\theta < \delta$, so the condition $\zeta > 1 - \delta$ always hold. 
 Hence, the conditions can be simplified as $\|\Omega - I\|_F^2 \goto 0$ and one of the following conditions can be satisfied:
\begin{itemize}
    \item[(a)] $\theta > \delta/2$, $1 - 2\zeta - 2\theta > 0$; or 
    \item[(b)] $\theta < \delta/2$, $ 1 - \zeta - 2\theta < 0$. 
\end{itemize}
Theorem \ref{thm:unknownmulower} is proved. 

\subsection{Proof of lower bound in Theorem \ref{thm:PCS}}\label{subsec:lowerpcs}
Consider the case $\Omega_0 = I$ without loss of generality. Because loss of information about $\Omega_1$, the region of impossibility cannot be larger than that in the case $\Omega_1$ is known in Theorem \ref{thm:unknownmulower}. Hence, $MR(L) \geq 1/2+o(1)$ when $2 - 2\alpha - \beta < 0$, $1 - 2\gamma < 0$ and one of the following conditions are satisfied: 
\begin{itemize}
    \item[(a)] $\theta > \delta/2$, $1 - 2\zeta - 2\theta > 0$; or 
    \item[(b)] $\theta < \delta/2$, $ 1 - \zeta - 2\theta < 0$. 
\end{itemize}
The region of impossibility in Theorem \ref{thm:PCS} is proved. 

\subsection{Proof of lower bound in Theorem \ref{thm:real}}\label{subsec:lowerreal}
When both $\mu$ and $\Omega_i$'s are unknown, with the same analysis in Section \ref{subsec:lowerpcs}, we have the region of impossibility in Theorem \ref{thm:real}.

\section{Proof of upper bounds}\label{sec:proof}
In this section, we present the proof of upper bounds in Theorem \ref{thm:unknownmu}, Proposition \ref{prop:unknownomega} and Theorem \ref{thm:PCS}. 
This section is structured as follows. In Section \ref{sec:lemmas}, we present some mathematical results as the preparations.  In Section \ref{sec:unknownmu}, we present the upper bounds of QDAw and QDAfs in Theorem \ref{thm:unknownmu}. We prove the case that $\Omega_1$ is unknown in Section \ref{sec:proofpcs}. All the proofs of the lemmas in this section can be found in the supplementary material \cite{QDAsupp}. 
In this section, we always use $\Omega = \Omega_1$ for simplification without confusion. 

We begin with the expression of the mis-classification rate $MR$ in terms of QDA. 
Given $\mu$ and $\Omega$, the two types of mis-classification rates are defined as 
\be\label{eqn:p}
p_{0,\mu,\Omega} = P_{Y = 0}(Q>0|\mu, \Omega), \quad
p_{1, \mu, \Omega} = P_{Y = 1}(Q < 0|\mu, \Omega).
\ee
Then, the population mis-classification rate ($MR$) of QDA is 
\be\label{eqn:mr}
MR(QDA) \ = \ [E[p_{0,\mu,\Omega}]+  E[p_{1,\mu,\Omega}]/2.
\ee
Given a parameter set $(\alpha, \beta, \gamma, \zeta, \theta)$, if both $E[p_{0, \mu, \Omega}]$ and $E[p_{1, \mu, \Omega}]$ converge to 0, then $MR(QDA)$ converges to 0, which means QDA is successful.


\subsection{Preparations and notations}\label{sec:lemmas}
To find the upper bounds, we should analyze the asymptotic distribution of the QDA score. In the analysis, we keep on using the quadratic terms of $X$, in the form of $X^\top A X + 2d^\top X$. The following lemma states the asymptotic distribution of such quadratic terms.
\begin{lemma}\label{lemma:quad}[Quadratic functional of normal distributions]
Consider $X \sim N(\mu, \Sigma)$ where $\Sigma$ is positive definite. Let $S = X^\top A X + 2d^\top X$ with a symmetric matrix $A$ and a vector $d$,  
\bea
&& E[S] = Tr(A\Sigma) + \mu^\top A \mu + 2d^\top \mu,\\
&& \var(S) = 2Tr((A\Sigma)^2) + 4(\mu^\top A\Sigma A\mu + \mu^\top A\Sigma d + d^\top \Sigma d).
\eea
\begin{itemize}
    \item[(a)] $\frac{\sum_{i=1}^p |\lambda_i|^3 (1 + |\tilde{\mu}(i)|^3)}{(\sum_{i=1}^p \lambda_i^2(1 + \tilde{\mu}_i^2))^{3/2}} \goto 0$; or 
    \item[(b)] $\var(S) = \sum_{i:\lambda_i = 0} \Tilde{d}^2(i)(1+o(1))$. 
\end{itemize}
\end{lemma}

\begin{lemma}\label{lemma:tweak}
Under current model and assumptions, for a given matrix $A$ with spectrum in $(1-o(1), 1+o(1))$, there exists a constant $C > 0$, so that with probability $1 - o(1)$, 
\[
\biggl|\left[\hat{\mu}_0^\top (I - A) \hat{\mu}_0 -\mu^\top(I-A)\mu\right] + \frac{1}{n_0}Tr(A - I)\biggr|   \leq 
 C\sqrt{\frac{\ln p}{n}}(\|A - I\|_F/\sqrt{n}+\|(I - A)\mu\|).
\]
\end{lemma}

In the analysis, we have to relate the terms $\|\Omega_i - I\|$,  $\|\Omega_i - I\|_F^2$ and $\|\mu\|^2$ to the constant parameters. The following two lemmas describe how these terms rely on the parameters.
\begin{lemma}\label{lemma1}[Bounds on the signals in precision matrix]
Under models (\ref{Omega0}) and (\ref{Param}), when $p\to\infty$, with probability $1-o(1)$, 
\ben
\|V^{(k)}\| \leq \eta_p b(p, \beta) = \left\{\begin{array}{ll}
3\eta_p\sqrt{p\nu} = 3\eta_pp^{(1-\beta)/2}, & 0<\beta < 1,\\
2\eta_p\sqrt{{\ln p}/{\ln\ln p}}, & \beta = 1,\\
{2\eta_p}/{(\beta - 1)}, & 1 < \beta \leq 2.
\end{array}
\right.
\een
\end{lemma}


%

The results are summarized in  in the following lemma.  
\begin{lemma}\label{lemmaeb}
Consider model (\ref{M12}) with the parameterizations (\ref{Param1}) and (\ref{Omega1})--(\ref{con1}). With probability $1 - o(1)$, we have $\|\Omega_k - I\|= o(1)$ and 
\be\ba{l}
\|V^{(k)}\|_F^2  =  \eta_p^2 p^2 \nu_p(1 + o(1)), \quad 
\|\Omega_k - I\|_F^2  =  p\xi_p^2+ \eta_p^2 p^2 
\nu_p(1 + o(1)),\\
\|\mu\|^2 = p\tau_p^2  \epsilon_p(1 + o(1)). 
 \ea
\ee
\end{lemma}




\subsection{Proof of Theorem \ref{thm:unknownmu}}\label{sec:unknownmu}

When $\mu$ is unknown, we propose two algorithms that work in different regions. When the non-zeros in $\mu$ are weak and relatively dense, then we apply QDAw which averages all the features; when the non-zeros in $\mu$ are relatively strong but sparse, we apply QDAfs to select features first. We find the upper bounds for both algorithms to prove Theorem \ref{thm:unknownmu}.

%
%

\subsubsection{Performance of QDAw}
\label{sec:weaksig} 
In QDAw, we estimate labels by $\hat{Y} = I(Q_w > 0)$. Here, $Q^w = S^w - T_S^w$, where  
\[
S^w = X^\top(I -{\Omega}) X + 2\hat\mu^\top (I + \Omega)X, \quad 
T_S^w = \hat{\mu}_0^\top (I - \Omega)\hat{\mu}_0 - \ln |{\Omega}| - \frac{1}{n_0}Tr({\Omega} - I).
\]
Here, $\hat{\mu}_0 = \frac{1}{\sum_i I\{Y_i = 0\}} \sum_{i:Y_i = 0}X_i$ as the average of training samples in Class 0, and $\hat\mu = a*{\bf 1}$, a vector with all the entries as $a$. In the algorithm, we take $a = p^{(c-1)/2}$.
The errors are $p_{i,\mu, \Omega} = P((-1)^i(S^w - T_S^w) > 0 )$. We want to find the region that both $p_{i, \mu, \Omega} \goto 0$.

Consider $S^w$, which is a quadratic term with $A = I - \Omega$ and $d = (I + \Omega)\hat{\mu}$. Apply Lemma \ref{lemma:quad} to $S^2$ with $\Sigma = I$ for the case $Y = 0$ and $\Sigma = \Omega^{-1}$ for the case $Y = 1$. There is 
\begin{eqnarray}\label{eqn:QDAwe}
E[S^w|Y = 0] & = & \mu^\top (I - \Omega)\mu +  Tr(I - \Omega) - 2a\mu^\top (I + \Omega){\bf 1},\\
E[S^w|Y = 1] & = & \mu^\top (I - \Omega)\mu +
  Tr(\Omega^{-1} - I) + 2a\mu^\top (I + \Omega){\bf 1},
\end{eqnarray}
and 
\begin{eqnarray}\label{eqn:QDAwvar}
\var(S^w|Y = i) = 2\|\Omega - I\|_F^2 + (16pa^2 + \|(\Omega - I)\mu\|^2)(1+o(1)), \,\, i = 0, 1.
\end{eqnarray}
Given $Y = i$, we define $Z_i = [S^w - E[S^w|Y=i]]/\sqrt{\var(S^w|Y=i)}$, then  
$\sup\nolimits_{-\infty < x < \infty}|F_{Z_i}(x) - \Phi(x)| \goto 0$. So the asymptotic distribution of $S^w$ is clear. When $Y = 0$ and $Y = 1$, the mean of $S^w$ differs in two parts, $\mu^\top (I + \Omega){\bf 1}$ and $\|\Omega - I\|^2_F$, with a shift that $\mu^\top (I - \Omega)\mu + \ln|\Omega|$. 

Compare $E[S^w|Y]$ with $T_S^w$, we can see $T_S^w$ mainly captures the shift. The difference is that $T_S^w$ uses $\hat{\mu}_0$ instead of the true parameter $\mu$. 
Consider the relative term $\hat{\mu}_0^\top (I - \Omega)\hat{\mu}_0$ in $T_S^w$. Apply Lemma \ref{lemma:tweak} to it with $A = \Omega$ and we have 
\begin{eqnarray}\label{eqn:QDAwt}
T_S^w = \mu^{\top}(I - \Omega) \mu - \ln|\Omega| + \Delta T,    
\end{eqnarray}
where
$|\Delta T| \leq C\sqrt{\ln p}(\|\Omega - I\|_F/n + \|(I - \Omega)\mu\|/\sqrt{n})$.

Introduce the results about $S^w$ and $T_S^w$ into $p_{i, \mu, \Omega} = P((-1)^i(S^w - T_S^w) > 0|Y = i)$. By the asymptotic normality of $S^w$, the error is $\Phi(\frac{(-1)^i*(E[S^w|Y=i] - T_S^w) }{\sqrt{\var(S^w|Y = i)}}) + o(1)$. 
Introduce in (\ref{eqn:QDAwe}), (\ref{eqn:QDAwvar}) and (\ref{eqn:QDAwt}) into $p_{0, \mu, \Omega}$, and we have 
\bea
p_{0, \mu, \Omega} 
& = & \Phi(\frac{ \ln|\Omega| + Tr(I - \Omega) - 2a\mu^\top(\Omega + I){\bf 1} + \Delta T }{\sqrt{2\|\Omega - I\|_F^2 + 16pa^2 + \|(\Omega - I)\mu\|^2}}) + o(1)\nonumber\\
& = & \Phi(\frac{ -\|\Omega - I\|_F^2/2 - 4a\|\mu\|_1(1+o(1)) + \Delta T }{\sqrt{2\|\Omega - I\|_F^2 + 16pa^2 + \|(\Omega - I)\mu\|^2}}) + o(1).\nonumber
\eea

Since $|\Delta T| \leq C\sqrt{\ln p}(\|\Omega - I\|_F/n + \|(I - \Omega)\mu\|/\sqrt{n}) \ll \sqrt{2\|\Omega - I\|^2_F + \|(\Omega - I)\mu\|^2}$ the denominator, so $\Delta T$ has negligible effects. 
Consider $\|(\Omega - I)\mu\|^2$. When $0 < \beta < 1$, then $\|(\Omega - I)\mu\| \ll \sqrt{p}\|\Omega - I\|\leq \|\Omega - I\|_F$ by Lemma \ref{lemma1}. 
When $1 \leq \beta < 2$, there are at most constant non-zeros in each row of $\Omega$. Hence, with probability $1 - o(1)$, $\|(\Omega -I)\mu\|^2 = \|V^{(1)}\mu + \xi \mu\|^2 \leq p\xi^2 + (\eta^2 p^2 \nu) \tau^2 \epsilon \ll \|\Omega - I\|^2_F$. In all, 
We only need to discuss 
\[
\frac{ -\|\Omega - I\|_F^2/2 - 4a\|\mu\|_1}{\sqrt{2\|\Omega - I\|_F^2 + 16pa^2}}
\leq \frac{ -\|\Omega - I\|_F^2/2 - 4a\|\mu\|_1}{2\max\{\sqrt{2}\|\Omega - I\|_F, 4a\sqrt{p}\}}.
\]

Now we discuss two cases: 
\begin{itemize}
\item Case 1. Suppose $\|\Omega - I\|_F^2 \gg a\sqrt{p} \goto \infty$. In this case, both $\|\Omega - I\|_F^2/\|\Omega - I\|_F$ and $\|\Omega - I\|_F^2/a\sqrt{p}$ go to infinity, and the term of interest goes to negative infinity.  
\item Case 2. Suppose $\sqrt{p}\tau \epsilon \goto \infty$. Then $\|\mu\|_1 \goto \infty$ with probability $1 - o(1)$. If $\|\Omega - I\|_F^2 \gg a\sqrt{p}$, then it comes to case 1 which is solved. If $\|\Omega - I\|_F^2 \ll a\sqrt{p}$, then the term of interest comes to $a\|\mu\|_1/4a\sqrt{p} = \sqrt{p}\tau\epsilon (1+o(1)) \goto \infty$. 
\end{itemize}
Therefore, $p_{i, \mu, \Omega} \goto 0$ with probability $1 - o(1)$, and $E[p_{i, \mu, \Omega}] \goto 0$. The same derivation holds for $p_{1, \mu, \Omega}$. 
As a conclusion, $MR(QDAw) \goto 0$ in this region. 

\subsubsection{Performance of QDAfs}\label{sec:qdafsmu} 
Now we consider the case $\tau \gg 1/\sqrt{n}$, i.e., $\theta < \delta/2$. The signals in $\mu$ are individually strong enough for successful recovery. Hence, we select features first, and then apply QDA on the post-selection data. 

The feature selection step is as follows. 
Define $d$ as 
\be\label{eqn:d}
d \ = \ \Omega \hat{\mu}_1 - \hat{\mu}_0 \ \sim \ N\bigl((I + \Omega)\mu, \ \frac{1}{n_0}I+\frac{1}{n_1}\Omega\bigr).
\ee 
When $\max_{1 \leq j \leq p} |d_i| > 2{\ln p}/\sqrt{n}$, we let $d^{(t)}_j = I(|d_j| \geq t)$ with the threshold $t= 2\sqrt{\ln p}/\sqrt{n}$. 
Define $\hat \mu_0^{(t)}=\hat{\mu}_0 \circ d^{(t)}$ and $\hat \mu_d^{(t)}=d\circ d^{(t)}$ as the post-selection estimators. 
Define $\Omega^{(d)}$ as the sub-matrix of $\Omega$ consisting of rows and columns that $d^{(t)} = 1$.
When $\theta < \delta/2$, this feature selection step happens with probability $1 - o(1)$. 
In supplementary materials \cite{QDAsupp}, it is shown that the signals can be exactly recovered with probability $1 - o(1)$. 
Hence, we only consider the event that $\{t = \sqrt{2\ln p/n}\}$ and all the signals are exactly recovered.  



In QDAfs, the criteria is updated as $Q^s = S^s - T_S^s$, where  
\[
S^s = X^\top(I -{\Omega}) X + 2\hat\mu_d^\top X, \quad 
T_S^s = (\hat \mu_0^{(t)})^\top (\Omega - I)\hat \mu_0^{(t)} - \ln |{\Omega}| - \frac{1}{n_0}Tr({\Omega^{(d)}} - I).
\]
Compare it with the ideal case that $\mu$ is known, the difference in the criteria is 
$\Delta Q = Q^s - Q(X, \mu, \Omega)$, where 
\be\label{Decomp2}
 \Delta Q =  \displaystyle 2(\hat{\mu}_d^{(t)} - (I + \Omega)\mu)^\top X+\bigl[(\hat{\mu}_0^{(t)})^\top(I-\Omega)\hat{\mu}_0^{(t)}-\mu^\top(I-\Omega)\mu + \frac{1}{n_0}Tr(\Omega^{(d)} - I)\bigr]\nonumber.
\ee

In Supplementary Materials \cite{QDAsupp}, we prove that, $\frac{Q(X, \mu, \Omega)}{\sqrt{2\|\Omega - I\|_F^2 + 16\|\mu\|^2}}$ is asymptotically normal distributed with mean $(-1)^{Y+1}\sqrt{\|\Omega - I\|_F^2/8 + \|\mu\|^2}$ and variance 1. Therefore, the mis-classification rate by $I\{Q(X, \mu, \Omega) > 0\}$ converges to 0 when the mean diverges.

When $\mu$ is unknown, the classification rule is $I\{Q^s = Q(X, \mu, \Omega) + \Delta Q > 0\}$. The error rate can be bounded by 
\bea\label{eqn:pthm3}
p_{i, \mu, \Omega} & = & P((-1)^i(Q(X, \mu, \Omega) + \Delta Q) > 0)\nonumber\\
& = & P(\frac{(-1)^i(\|\Omega - I\|_F^2/2 + 4\|\mu\|^2) + \Delta Q}{\sqrt{2\|\Omega - I\|_F^2 + 16\|\mu\|^2}} > 0), \quad i = 0, 1.
\eea
Therefore, $|\Delta Q| \leq \sqrt{2\|\Omega - I\|_F^2 + 16\|\mu\|^2}$ with probability $1 + o(1)$ suffices to show the success of QDAfs. 

\begin{lemma}\label{lemma:StrongDeltaQ}
Under the model assumptions and the definition of $\Delta Q$, with probability $1 - o(1)$, there is
\be\label{eqn:lemmaDeltaQ2}
 |\Delta Q|   \leq   O(\sqrt{p \epsilon_p (\xi_p^2 + p\epsilon_p \eta_p^2 \nu_p)}/n) + \sqrt{p\epsilon_p} \tau_p \ln p (1 + o(1)).
\ee
\end{lemma}
By Lemma \ref{lemma:StrongDeltaQ} about the magnitude of $\Delta Q$,  when $p\xi_p^2 + \eta_p^2p^2\nu_p \goto \infty$ or $\tau_p^2p\epsilon_p \goto \infty$, 
\[
|\Delta Q|
\ll 
\sqrt{p\xi_p^2/8 + \eta_p^2p^2\nu_p/8 + \tau_p^2p\epsilon_p}(1 + o(1)) = \sqrt{2\|\Omega - I\|_F^2 + 16\|\mu\|^2}.
\]
Therefore, in the region of possibility identified by part (ii) of Theorem \ref{thm:unknownmu}, $MR(QDAfs)$ converges to 0. \qed

\subsection{Proof of Theorem \ref{thm:PCS}}\label{sec:proofpcs}
To prove Theorem \ref{sec:proofpcs}, we start with the proof of Proposition \ref{prop:unknownomega} when $\mu$ is known and $\Omega$ is estimated by PCS in Section \ref{sec:proofprop}. The effects of estimated $\Omega$ can be found. Then we use the result to prove Theorem \ref{thm:PCS}.

\subsubsection{Proof of Proposition \ref{prop:unknownomega}}\label{sec:proofprop}
When $\mu$ is known and $\Omega$ is estimated by PCS, we classify by $\hat{Y} = I(Q(X, \mu, \hat{\Omega}) > 0)$, where 
\[
Q(X, \mu, \hat{\Omega}) = X^\top (I - \hat{\Omega})X + 2\mu^\top(I + \hat{\Omega})X + \mu^\top(I - \hat{\Omega})\mu + \ln|\hat{\Omega}|. 
\]
We do not need to consider QDAw or QDAfs, and the focus is on $\hat{\Omega}$ by PCS only. 

Let $Q(X, \mu, \hat{\Omega}) = S^{PCS} - T_S^{PCS}$, where $S^{PCS} = X^\top (I - \hat{\Omega})X + 2\mu^\top(I + \hat{\Omega})X$, $T_S = \mu^\top(\hat{\Omega} - I)\mu - \ln|\hat{\Omega}|$. 
Note that $X$ and $\hat{\Omega}$ are independent. Given $\hat{\Omega}$, we derive the asymptotic distribution of $S^{PCS}$ by Lemma \ref{lemma:quad}. 
In details, the expectations and variances are 
\begin{itemize}
    \item $E[S^{PCS}|Y = 0] = T_S^{PCS} -4\mu^\top \hat\Omega \mu  + \ln|\hat\Omega| + Tr(I - \hat\Omega)$;
    \item 
    $E[S^{PCS}|Y = 1] = T_S^{PCS} + 4\mu^\top  \mu  + \ln|\hat\Omega| + Tr(\Omega^{-1}(I-\hat{\Omega}))$; 
    \item $\var(S^{PCS}|Y = 0) = 2Tr((\hat\Omega - I)^2) + 16 \mu^\top \hat\Omega^2 \mu$;
    \item $\var(S^{PCS}|Y = 1) = 2Tr((\Omega^{-1} - \hat\Omega\Omega^{-1})^2) + 16 \mu^\top \Omega^{-1} \mu$. 
\end{itemize}
Define $Z_i = [S^{PCS} - E[S^{PCS}|Y = i]]/\sqrt{\var(S^{PCS}|Y = i)}$ and $F_{Z_i}(x) = P(Z_i \leq x)$, then $\sup\nolimits_{-\infty < x < \infty}|F_{Z_i}(x) - \Phi(x)| \goto 0$.

The mean and variance for $S^{PCS}|Y = 0$ are similar as those of the ideal case, except all $\Omega$ are replaced by $\hat\Omega$. 
With similar derivations, we have that 
\be\label{eqn:p0prop}
p_{0, \mu, \Omega} = \Phi(\frac{E[S^{PCS}|Y = 0] - T_S^{PCS}}{\sqrt{\var(S^{PCS}|Y = 0)}}) = \Phi(-\sqrt{\|\hat\Omega - I\|^2_F/8 + \|\mu\|^2}) + o(1).
\ee

The derivation for $p_{1, \mu, \Omega}$ is more complicated. Both the mean and variance of $S^{PCS}|Y = 1$ involves the term $Tr(\Omega^{-1}(I-\hat{\Omega}))$, which is related to both $\hat{\Omega}$ and $\Omega$. To bound it, we compare the term with $Tr(\hat\Omega^{-1}(I-\hat{\Omega}))$. The difference between them is $\Delta V  = Tr((\Omega^{-1} - \hat{\Omega}^{-1})(I-\hat{\Omega}))$. The goal is to bound $|\Delta V|$. 

For any square matrices $A$ and $B$ with ordered singular values as $\alpha_i$ and $\beta_i$, respectively. By Von Neuman's trace inequality in \cite{VonTrace}, $|Tr(AB)| \leq \sum \alpha_i \beta_i \leq \sqrt{\sum \alpha_i^2 \sum \beta_i^2} = \sqrt{Tr(A^\top A)}\sqrt{Tr(B^\top B)}$. 
Apply this result to $\Delta V$ and recall that both $\Omega$ and $\hat{\Omega}$ has eigenvalues at $1+o(1)$. Then we have 
\begin{eqnarray*}
|\Delta V| & \leq & \sqrt{Tr((\Omega^{-1} - \hat{\Omega}^{-1})^2)} \sqrt{Tr((I - \hat{\Omega})^2)}\\
&\leq & \sqrt{Tr(\Omega^{-2}(\Omega - \hat{\Omega})^2\hat{\Omega}^{-2})} \|I - \hat{\Omega}\|_F(1+o(1)).\\
& \leq & \|\hat{\Omega} - \Omega\|_F \|I - \hat{\Omega}\|_F(1+o(1)).
\end{eqnarray*}

Introduce the bound of $\Delta V$ into $p_{1,\mu,\Omega}$,  
\begin{eqnarray}\label{eqn:p1prop}
p_{1, \mu, \Omega} & = &\Phi(\frac{T_S^{PCS} - E[S^{PCS}|Y = 1]}{\sqrt{\var(S^{PCS}|Y = 1)}}) \nonumber\\
&\leq &\Phi(\frac{-4\|\mu\|^2 - \|\hat\Omega - I\|^2_F/2  + \|\hat{\Omega} - \Omega\|_F \|I - \hat{\Omega}\|_F}{\sqrt{2\|\hat{\Omega}-I\|_F^2 + 16\|\mu\|^2}}) + o(1).
\end{eqnarray}

For $p_{i,\mu,\Omega}$, now we only need to consider $\|\hat{\Omega} - I\|^2_F$ and $\|\hat{\Omega} - \Omega\|^2_F$. 
According to Theorem 2.3 in \cite{PCS}, when $1 - \delta/2 < \beta < 2$ and $\eta \gg 1/\sqrt{n}$, PCS recovers the exact support with probability $1 - o(1/p^2)$, and 
$\max_{i,j}|\Omega(i,j) - \hat{\Omega}(i,j)| \leq C\sqrt{\ln p/n}$. 
Since $|{\Omega}(i,j)| \gg \sqrt{\ln p/n}$ on the off-diagonals, the estimation error is at a smaller order than the off-diagonal signals in $\Omega$. On the diagonals, we have to consider several cases. 
\begin{itemize}
    \item Case 1. $\xi_p \gg 1/\sqrt{n}$. With probability $1 - o(1)$, $|\hat{\Omega}(i,i) - \Omega(i,i)| \leq \sqrt{\ln p/n}$, which is at a smaller order than $\xi_p = |\Omega(i,i) - 1|$, for all $i$. Therefore, $\|\hat\Omega - I\|_F^2 = \|\Omega - I\|_F^2(1+o(1))$ and $\Delta V$ is negligible compared to $\|\hat\Omega - I\|_F^2$. Since $\xi \gg 1/\sqrt{n}$, $p\xi_p^2 \goto \infty$, therefore $\|\Omega - I\|_F^2 \goto \infty$ and $p_{1, \mu, \Omega} \goto 0$.
\item Case 2. $\xi_p \ll \max\{\eta_p\sqrt{p\nu_p}, 1/\sqrt{p}\}$. 
When $\xi_p \ll 1/\sqrt{n}$, with probability $1-o(1)$, $|\hat\Omega(i,i)-1| \leq \ln p/\sqrt{n}$ for all $i$ and therefore the diagonals of $\hat{\Omega}$ will be updated to 1. Hence, $\|I-\hat{\Omega}\|_F^2 = \eta_p^2p^2\nu(1+o(1))$ and 
$\|\Omega-\hat{\Omega}\|_F^2 = p\xi_p^2 + p^2\nu_p\ln p/n(1+o(1))$. When $\xi_p \ll \max\{\eta_p\sqrt{p\nu_p}, 1/\sqrt{p}\}$, $\|\hat\Omega - I\|_F^2 = \|\Omega - I\|_F^2(1+o(1)) + o(1)$ and $\Delta V$ is either $o(1)$ or negligible compared to $\|\hat{\Omega} - I\|_F^2$.
\end{itemize}

Introduce these terms into (\ref{eqn:p0prop}) and (\ref{eqn:p1prop}), we can see $p_{i,\mu,\Omega} \goto 0$ when a) $\xi \gg 1/\sqrt{n}$, or b) $\xi \ll \max\{\eta\sqrt{p\nu}, 1/\sqrt{p}\}$ and $\|\Omega - I\|_F^2 \goto \infty$ or $\|\mu\|^2 \goto \infty$. Proposition \ref{prop:unknownomega} is proved. 

\subsubsection{Proof of Theorem \ref{thm:PCS}}
We examine the performance of QDAw with PCS for the region $\tau_p \ll 1/\sqrt{n}$ and that of QDAfs with PCS for the region $\tau_p \gg 1/\sqrt{n}$. 

We first consider the weak signal region that $\tau \ll 1/\sqrt{n}$. Here we use adjusted PCS to estimate $\Omega$ and a constant vector $\hat{\mu} = a*{\bf 1}$ to estimate the mean vector. We classify $X$ to be in class 0 if $Q(X, \hat{\mu}, \hat{\Omega}) < 0$, where 
\[
Q(X, \hat{\mu}, \hat{\Omega}) = 
X^\top(I - \hat{\Omega})X + 2\hat{\mu}^\top (I + \hat{\Omega})X - \hat{\mu}_0^\top (I - \hat{\Omega})\hat{\mu}_0 + \ln|\hat{\Omega}| + \frac{1}{n_0}Tr(\hat{\Omega} - I).
\]
We rewrite it as $Q(X, \hat{\mu}, \hat{\Omega}) = S^{w, pcs} - T_S^{w, pcs}$, where $S^{w, pcs} = X^\top(I - \hat{\Omega})X + 2\hat{\mu}^\top (I + \hat{\Omega})X$ and $T_S^{w, pcs} = \hat{\mu}_0^\top (I - \hat{\Omega})\hat{\mu}_0 - \ln|\hat{\Omega}| - \frac{1}{n_0}Tr(\hat{\Omega} - I)$. 

Apply Lemma \ref{lemma:quad} to $S^{w, pcs}$ and we can prove that, 
\begin{itemize}
\item the expectations are 
\begin{eqnarray*}
E[S^{w,pcs}|Y = 0] & = &\mu^{\top}(I - \hat\Omega)\mu +   Tr(I - \hat\Omega) - 2a\mu^\top (I + \hat\Omega){\bf 1},\\
E[S^{w,pcs}|Y = 1] & = & \mu^{\top}(I - \hat\Omega)\mu +  Tr(\Omega^{-1}(I - \hat\Omega)) + 2a\mu^\top (I + \hat\Omega){\bf 1}.  
\end{eqnarray*}
\item when $p \to \infty$, the asymptotic variances are
\[
\var(S^{w,pcs}|Y = i) = 2\|\hat\Omega - I\|_F^2 + (16pa^2 + \|(\hat\Omega - I)\mu\|^2)(1+o(1)), \qquad i = 0, 1.
\]
\end{itemize}
Further, $S^{w,pcs}|Y = i$ normalized by mean and variance  converges to normal distribution when $p \goto \infty$.

Therefore, the error rates $p_{i,\mu, \Omega}$ can be approximated by 
\bea
p_{0, \mu, \Omega} & = & \Phi(\frac{ E[S^{w,pcs}|Y = 0]  - T_S^{w, pcs}  }{\sqrt{\var(S^{w,pcs}|Y = 0)}}) + o(1)\nonumber\\
& \leq & \Phi(\frac{ -\|\hat\Omega - I\|_F^2/2 - 4a\|\mu\|_1(1+o(1)) + \Delta T }{\sqrt{2\|\hat\Omega - I\|_F^2 + 16pa^2 + \|(\hat\Omega - I)\mu\|^2}}) + o(1),\nonumber
\eea
where $\Delta T = \hat{\mu}_0^{\top} (I - \hat{\Omega}) \hat{\mu}_0 - \mu^{\top} ( I - \hat{\Omega}) \mu + \frac{1}{n_0}Tr(\hat{\Omega} - I)$. 

Apply Lemma \ref{lemma:tweak} to $\Delta T$ with $A = \hat\Omega$, $|\Delta T| \leq C\ln p(\|\hat\Omega - I\|_F/n + \|(I - \hat\Omega)\mu\|/\sqrt{n}) \ll \sqrt{2\|\hat\Omega - I\|^2_F + \|(\hat\Omega - I)\mu\|^2}$, so $\Delta T$ has negligible effects. In Section \ref{sec:weaksig}, we found  $\|(\hat\Omega - I)\mu\| \ll \|\hat\Omega - I\|_F$ holds with probability $1-o(1)$. 
Hence, we only need to discuss 
\[
\frac{ -\|\hat\Omega - I\|_F^2/2 - 4a\|\mu\|_1}{\sqrt{2\|\hat\Omega - I\|_F^2 + 16pa^2}}
\leq \frac{ -\|\hat\Omega - I\|_F^2/2 - 4a\|\mu\|_1}{2\max\{\sqrt{2}\|\hat\Omega - I\|_F, 4a\sqrt{p}\}}.
\]
In Section \ref{sec:proofprop}, we have found $\|\hat\Omega - I\|_F = \|\Omega - I\|_F(1+o(1)) + o(1)$ in the current region of interest. Hence, it comes back to the equation when $\Omega$ is known. In the region of possibility identified by part (i) of Theorem \ref{thm:PCS}, $MR(QDAw) \goto 0$.

Now we consider the case $\tau_p \gg 1/\sqrt{n}$, i.e., $\theta < \delta/2$. The signals in $\mu$ are individually strong enough for successful recovery. Hence, we estimate $\Omega$ by PCS, then threshold on $d = \hat{\Omega}\hat{\mu}_1 - \hat{\mu}_0$. QDA is applied to the post-selection data.

In \cite[Appendix C.1]{QDAsupp}, it is shown that the signals can be exactly recovered with probability $1 - o(1)$. 
Hence, we only consider the event that $\{t = \sqrt{2\ln p/n}\}$ and all the signals are exactly recovered.  

By Proposition \ref{prop:unknownomega}, we analyze the performance of $Q(X, \mu, \hat{\Omega}) = S^{PCS} - T_S^{PCS}$. In QDAfs, the criteria is updated as 
\be\label{eqn:deltaqsigmastrong}
Q(X, \hat{\mu}, \hat{\Omega})  =
Q(X, \mu, \hat{\Omega}) + \Delta Q, 
\ee
where $\Delta Q =  \displaystyle 2(\hat{\mu}_d^{(t)} - (I + \hat\Omega)\mu)^\top X+\bigl[(\hat{\mu}_0^{(t)})^\top(I-\hat\Omega)\hat{\mu}_0^{(t)}-\mu^\top(I-\hat\Omega)\mu + \frac{1}{n_0}Tr(\hat\Omega^{(d)} - I)\bigr]$.
The following lemma bounds $|\Delta Q|$. 
\begin{lemma}\label{lemma:unknownstrongQ}
Under the model assumptions and the definition of $\Delta Q$, there is
\be
 |\Delta Q|   \leq  \eta \tau \max\{p\epsilon \nu, 1\} \ln p +  O_p(\sqrt{4n^{-1}p\epsilon}).
\ee
\end{lemma}

Combining Lemma \ref{lemma:unknownstrongQ} with Section \ref{sec:proofprop} about $Q$, the errors are 
\bea
p_{i, \mu, \Omega} & = & P((-1)^i*(S^{PCS} - T_S^{PCS} + \Delta Q) > 0)\nonumber\\
& = & \Phi(\frac{(-1)^i*(T_S^{PCS} - E[S^{PCS}|Y = i]) }{\sqrt{\var(S^{PCS}|Y = i)}} + \frac{(-1)^i*\Delta Q}{\sqrt{\var(S^{PCS}|Y = i)}}) + o(1), \quad i = 0, 1.
\eea
The first term $\frac{(-1)^i*(T_S - E[S^{PCS}|Y = i]) }{\sqrt{\var(S^{PCS}|Y = i)}}  = -\sqrt{\|\Omega - I\|^2_F/8 + \|\mu\|^2}(1 + o(1)) + o(1)$ in Section \ref{sec:proofprop}. The second term can be bounded by
\[
\frac{|\Delta Q|}{\sqrt{\var(S^{PCS}|Y = i)}} 
\leq 
\frac{ \eta \tau \max\{p\epsilon \nu, 1\} \ln p +  O_p(\sqrt{4n^{-1}p\epsilon})}{\sqrt{p\xi^2/8 + \eta^2p^2\nu/8 + \tau^2p\epsilon}(1 + o(1))}.
\]
It goes to 0 in the region of possibility identified in part (ii) of Theorem \ref{thm:PCS}. 

Therefore, in the region of possibility identified by part (ii) of Theorem \ref{thm:PCS}, $MR(QDAfs)$ converges to 0. \qed

\section{Real Data Analysis}\label{sec:rats}
In this paper, we consider the rats dataset present in \cite{liver}. As we introduced in Section \ref{sec:quickdata}, this data set record the gene expressions of live rats in response to different drugs and toxicant. There are 181 samples and 8491 genes, where 61 samples are labeled as toxicant and the other 120 are labeled as other drugs. 
We compare QDA with LDA, where the latter one is shown to enjoy the best performance compared to classifiers such as SVM, RandomForest, GLasso and FoBa. 
The QDA with feature selection for the real data is discussed in Section \ref{sec:real-proc} and the implementation details and results are in Section \ref{sec:dataresults}. 


\subsection{Procedure for the real data}\label{sec:real-proc}
Here, we present a procedure for the classification based on QDA for the real data. For the real data, we have to estimate $\Omega_0$, $\Omega_1$, $\mu_0$ and $\mu_1$ separately. Further, we need to eliminate the effect of the feature variances. Hence, there is an additional scaling step in the following algorithm. 

\begin{table}[ht!]
\caption{Algorithm 2: Pseudocode for QDA with feature selection on real data} 
\scalebox{1}{ 
{\begin{tabular}{ll}\\ \hline
& \underline{Input}: data points $(X_i, Y_i)$, $1 \leq i \leq n$;  threshold $t > 0$; new data point $X$; tuning parameters: $C$, $t$.\\
& \underline{Output}: label $\hat{Y}$. \\
1. & Find $\hat{\Omega}_0$ and $\hat{\Omega}_1$ by PCS. Let $\hat{\mu}_0 = \frac{1}{n_0} \sum_{i: Y_i = 0} X_i$  and $\hat{\mu}_1 = \frac{1}{n_1} \sum_{i: Y_i = 1} X_i$, where $n_1 = \sum_{i =1}^n Y_i$ \\
&and $n_0 = n - n_1$.   \\ 
2. & Let $ \hat{\Omega}_ {\rm diff} = \hat{\Omega}_0 - \hat{\Omega}_1 - diag(\hat{\Omega}_0 - \hat{\Omega}_1)$.\\
3.  & Let $d_0=\hat{\Omega}_0\hat{\mu}_0/s_0$, $d_1= \hat{\Omega} \hat{\mu}_1/s_1$ and $d = d_1 - d_0=(d(1),\dots, d(p))^\top$.  Here $s_i$ are the standard \\
&deviation vector of the train data from class $i$, $i=0, 1$. The division  means element-wise division. 
\\
4.  & \underline{Thresholding}: Let $d^{(t)}$ denote the indicator vector of feature selection, i.e. $d^{(t)}(j) = 1\{|d(j)| \geq t\}$, \\
&for $j=1,\dots,p$. Let $\hat \mu_d^{(t)}$ be the hard-thresholded $\hat \mu_d^{(t)}=d\circ d^{(t)}$. \\
5.  & Scale $X$ as $x_j = [X_j - \bar{\mu}_j]/s_j$, where $\bar{\mu}=({\hat{\mu}_1+\hat{\mu}_0})/{2}$ and $s$ is the standard error of the pooled data \\
&$s_j   =    \sqrt{[(n_0-1) ((s_0)_j)^2+(n_1-1)((s_1)_j)^2] / (n_0+n_1-2)}$.
\\
6.  & \underline{QDA Score}: Calculate the QDA score $Q=x^\top \hat{\Omega}_ {\rm diff}x + 2( \hat\mu_d^{(t)})^\top X + C$. \\
7. &\underline{Prediction}: Predict $\hat Y = I\{Q > 0\}$.  \\
 \hline
\end{tabular}}
} 
\label{tab:alg2} 
\end{table}

Here are two tuning parameters, $t$ and $C$. In the implementations, we use a grid search to find the optimal values of them. Details in Section \ref{sec:dataresults}. 


\subsection{Implementation and Results}\label{sec:dataresults}
Following the setup of the data analysis in \cite{PCS}, we apply 4-fold data splitting to the sample. 
For each class, we randomly draw one fourth of the samples, and then combine them to be the test data while using the leftover to be the training data. We do the splitting for 15 times independently and record the error with QDA and LDA for each splitting. 
The data (sample indices) for the 15 splittings is available upon request. 

In the real data analysis section, we focus on comparing QDA and LDA. The LDA is implemented within the setting of QDA, where in Step (3)  of the algorithm in Section \ref{sec:real-proc} we use clipping thresholding instead of hard thresholding, and in Step (5) we set  $\hat{\Omega}_{\mbox{{\footnotesize diff}}}={0}$ for LDA. 
The clipping threshold is employed since it gives much more satisfactory results than hard thresholding for LDA; details in \cite{PCS}. For QDA, the two ways give similar results. 
Since the calculation of $\hat{d}$ involves the calculation of $\hat{\Omega}_0$ and $\hat{\Omega}_1$ and the thresholding, LDA algorithm has exactly the same tuning parameters with QDA. 
The procedure of determining these tuning parameters are the same for both algorithms, so that the results are comparable.

For PCS, there are four tuning parameters $(q_1, q_2, \delta, L)$. Here we use the same set of tuning parameters for the estimation of both $\hat{\Omega}_0$ and $\hat{\Omega}_1$, since the two classes are from the same data set and the performance of PCS is not sensitive to the choice of these parameters (\cite{PCS}). 
Following the setting in \cite{PCS}, we set $(\delta, L)=(.1, 30)$, and also tried $(\delta, L)=(.1, 50)$. For $(q_1, q_2)$, we consider $.1 \leq q_k \leq 1$, with an increment of .1, $k = 0, 1$. The selection is done by grid search. 

For Algorithm 2, there are two tuning parameters $t$ and $C$. We set the ranges $[t_{\mbox{{\footnotesize min}}}, t_{\mbox{{\footnotesize max}}}]=[0, \mbox{max}_{1\leq j \leq p} |d_j|]$ with an increment of .1 and $[C_{\mbox{{\footnotesize min}}}, C_{\mbox{{\footnotesize max}}}]=[-50, 50]$ or $[C_{\mbox{{\footnotesize min}}}, C_{\mbox{{\footnotesize max}}}]=[-100, 100]$ with an increment of 1.
The smallest error is obtained over a grid search of $t$, $C$, and $(q_1, q_2)$. 
This step is the same for both QDA and LDA to be fair. We compare the smallest error that LDA and QDA can achieve. 


\begin{figure}[htb!]
    \centering
   \subfigure[ $(C_{\mbox{{\footnotesize min}}}, C_{\mbox{{\footnotesize max}}})=(-50, 50)$ ]{ \includegraphics[width= 2.37 in]{Figure1.pdf}} \quad\quad\quad
		\subfigure[$(C_{\mbox{{\footnotesize min}}}, C_{\mbox{{\footnotesize max}}})=(-100, 100)$]{\includegraphics[width= 2.37 in]{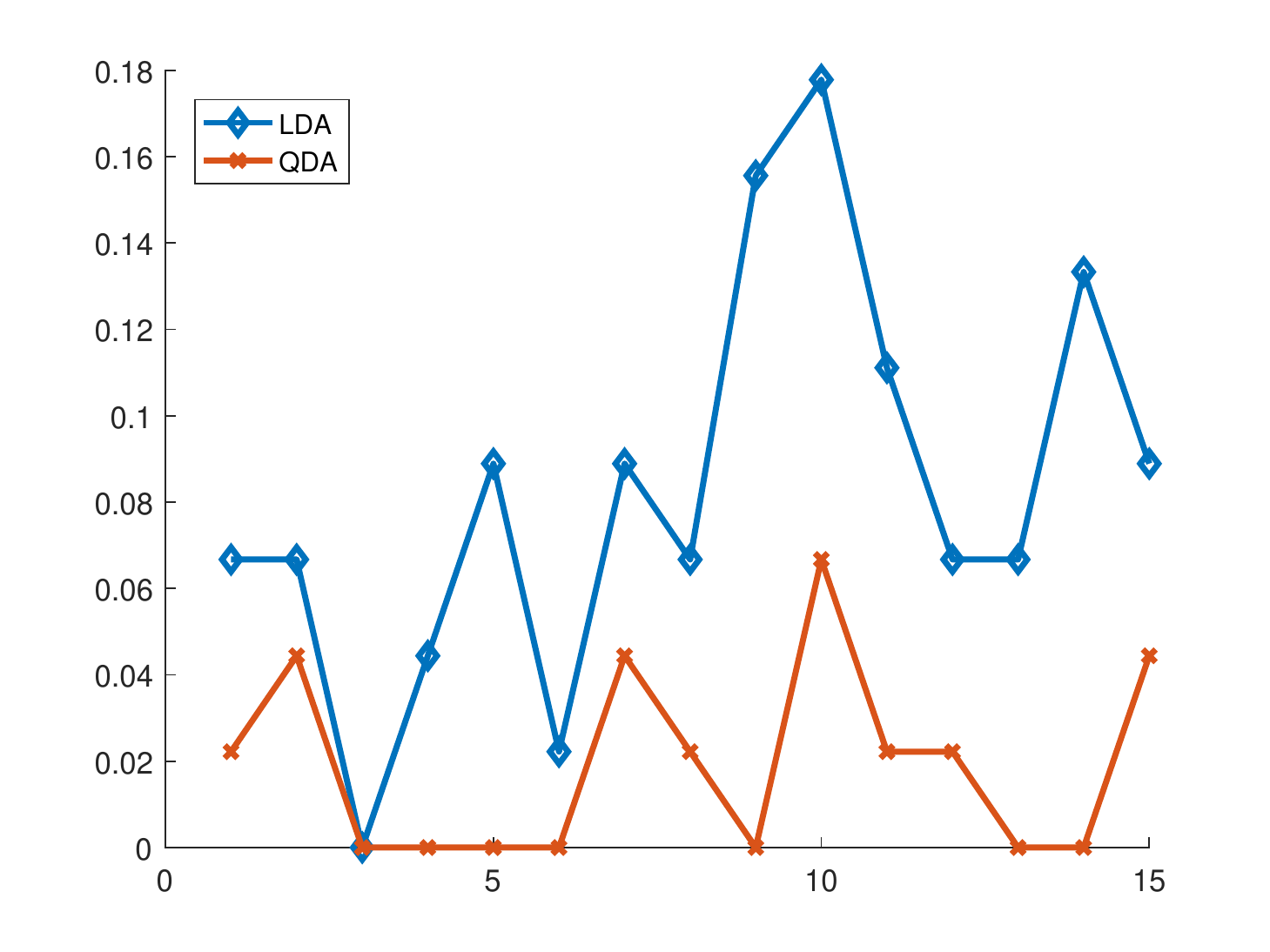}}
      \caption{Comparison of testing error rate (y-axis) of LDA and QDA for the rats data with $(\delta_1, L_1)=(\delta_2, L_2)=(.1, 30)$ and 15 data splittings.}
    \label{figure:errorrats-k30}
\end{figure}

Both the LDA test error (the best error) and the QDA test error (the best error) over all 15 data splittings are reported in Figure \ref{figure:errorrats-k30}.
In the left panel of Figure \ref{figure:errorrats-k30}, we can see that the error rates of LDA are all above QDA at every data splitting.  To better show the difference between them, we also plot the testing error rate in the right panel of Figure \ref{figure:errorrats-k30} for a wider grid-search range that $[C_{\mbox{{\footnotesize min}}}, C_{\mbox{{\footnotesize max}}}]=[-100, 100]$. 

\begin{figure}[htb!]
    \centering
    \subfigure[$L_1 = L_2 = 50$]{\includegraphics[width= 2.37 in]{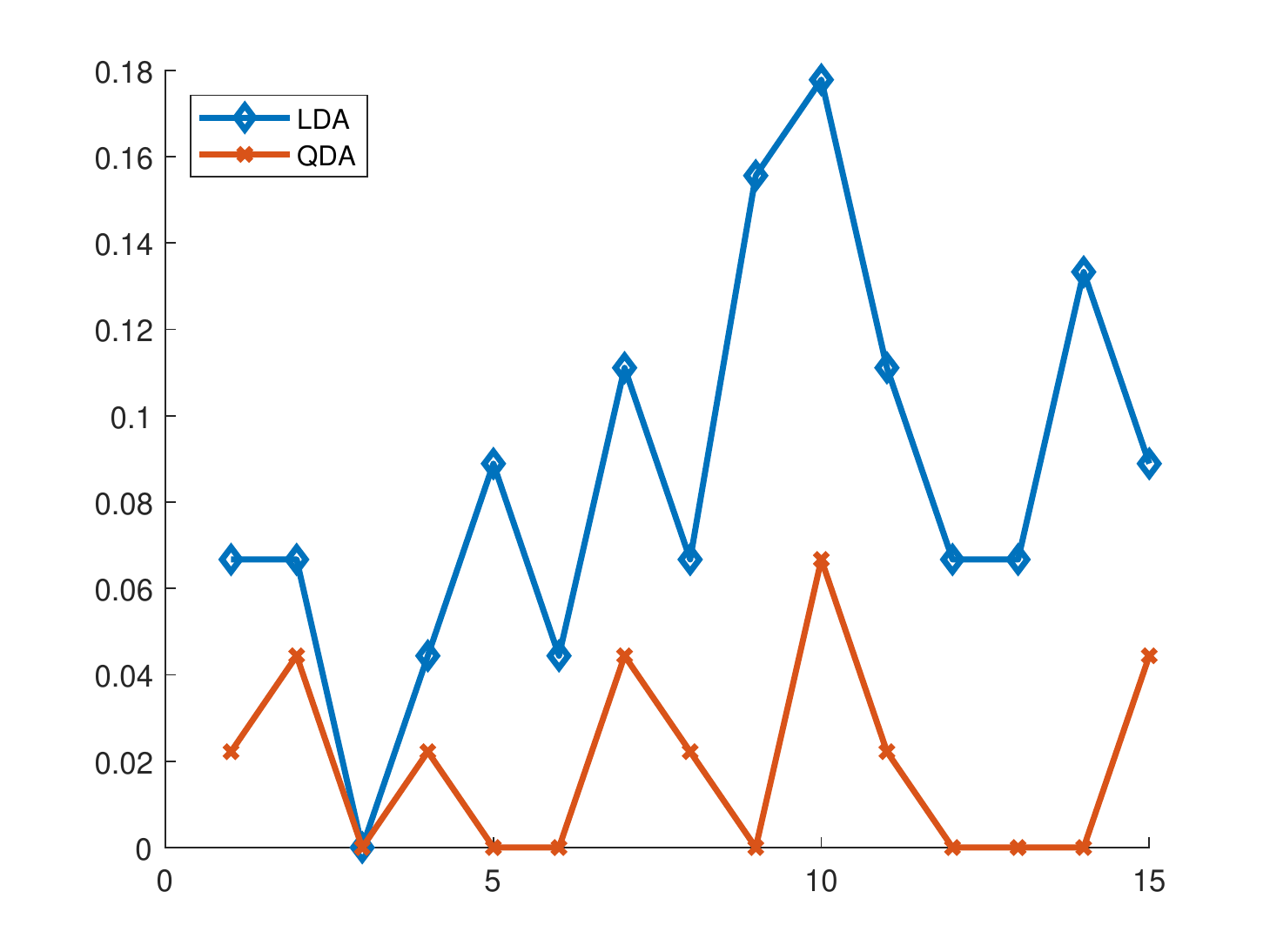}}\quad\quad\quad\quad
        \subfigure[$(q_1, q_2) \in (.1,1) \times (.1, 1)$]{ \includegraphics[width= 2.37 in]{Figure2}} \quad
      \caption{Comparison of testing error rate of LDA and QDA for the rats data on: (a) $(\delta_1, L_1)=(\delta_2, L_2)=(.1, 50)$ among 15 data splittings ($x$-axis); (b) varying choices of $(q_1, q_2) \in (.1,1) \times (.1, 1)$ for one splitting of 15 splittings in Figure 4(a). }
    \label{figure:errorrats-k50}
\end{figure}

%
The impact of the tuning parameters in PCS is presented in Figure \ref{figure:errorrats-k50}. 
When $L$ changes from 30 to 50, the results are summarized in subfigure (a), which is similar. This comparison clearly demonstrates the expected superiority of QDA over LDA. 
When $(q_1, q_2)$ changes, the results for one splitting are presented in subfigure (b). It suggests a proper choice of the tuning parameters will largely improve the QDA results, and overcome the LDA classifier. 

As a conclusion, the results suggest that, for rats data, QDA  outperforms LDA in terms of both best error rate and average error; with the results in \cite{PCS} for other methods, where the authors have shown that HCT-based LDA significantly outperforms all other HCT-based methods as well as SVM and RF, our findings also suggest that the QDA gives a better separation than the LDA  by taking into account the second order difference between the two classes.


\section{Discussion}\label{sec:dis}
This paper focuses on the classification problem associated with the use of QDA and feature selection for data of rare and weak signals. We derived the successful and unsuccessful classification regions, by using first the case of a known mean vector and covariance matrix, then the case of an unknown mean vector but known covariance matrix, and finally the case in which both mean vector and covariance matrix were unknown. We also proved that these regions were actually the possibility and impossibility regions under the same modeling, which indicates that QDA achieves the optimal classification results in this manner. In addition, we developed computing and classification algorithms that incorporated feature selection for rare and weak data. With these algorithms, our real data analysis showed that QDA had much-improved performance over LDA. 

Our theoretical results showed that the two sets of signal weakness and sparsity parameters, one set from the mean vector and the other set from the covariance matrix, influence the possibility/impossibility regions or QDA successful/unsuccessful regions almost independently (except for a $\max$ operator over the two sets of parameters) when the covariance matrix is known. When both the mean vector and covariance matrix are unknown, the two sets of parameters interact with each other as indicated in Theorem \ref{thm:PCS}. For the latter case, the analysis of the mis-classification rate is very complicated and we only obtained partial results for this most general case; further study is therefore warranted. Also, for the precision matrix $\Omega$ given in (\ref{Omega1}), we can introduce sparsity and weakness in the diagonal elements of $I-\Omega$, the difference in precision matrices, instead of using a constant $\xi=1-c$ for all diagonal elements.


\bibliographystyle{imsart-number}
\bibliography{Zhigang}

\newpage
\appendix

This file contains six sections. In Section~\ref{app:qda}, we present theoretical results  about the classical QDA in high-dimensional setting. We have the phase diagram in several scenarios, which cannot achieve the statistical lower bound. The proof of these results are deferred to Section~\ref{app:proof}. 
Section~\ref{app:ideal} proves Proposition \ref{thm:idealunequal}. 
Section~\ref{app:real} proves Theorem \ref{thm:real}, which contains the proof on screening accuracy. Section~\ref{app:seclemma} is to prove the asymptotic normality of the quadratic forms. 
Section~\ref{app:lemma} discusses the lemmas and useful properties in the proof of the theorems. 

\section{QDA on the weak and dense case}\label{app:qda}
In the main paper, we proposed QDAw for the case that $\mu$ has relatively dense and weak signals. The upper bound of QDAw matches the statistical lower bound for high dimensional classification problem. 
If we simply apply QDA without feature selection in this case, then the upper bound cannot match the lower bounds. 

In this section, we present our result about QDA without QDAw. It means, we apply QDA with feature selection when $\mu$ has strong signals and QDA without feature selection with estimated mean and precision matrix when $\mu$ has relatively weak signals. We still consider the following scenarios: 
\begin{enumerate}
\item $\mu$ is unknown, $\Omega_0 = I$ is known, $\Omega_1$ is known;
\item $\mu$ is unknown, $\Omega_1$ is unknown, $\Omega_0 = I$ is known.
\end{enumerate}
Remark. When both $\mu_i$ and $\Omega_i$'s are known, then it is the ideal case where we can apply the original QDA directly. It is the same with Proposition \ref{thm:idealunequal}. 

\subsection{Main results when precision matrices are known}\label{app:qdatheorem}
For the new QDA approach, the information in the quadratic term can be summarized by $\|\Omega - I\|_F^2$ and the information in the linear term by $\|\mu\|^2$. According to the model parameterizations and random matrix theory, 
$\|\Omega - I\|_F^2 \approx p^{1 - 2\gamma} + p^{2-2\alpha - \beta}$ and $\|\mu\|^2 \approx p^{1 - 2\theta - \zeta}$. When $\gamma < 1/2$, the information on the precision matrix diagonals ($p^{1 - 2\gamma}$) is sufficient for a satisfactory classification result (Proposition \ref{thm:idealunequal}). So we only consider the non-trivial case that $\gamma > 1/2$. 

When $\gamma > 1/2$,  $\|\Omega - I\|_F^2 \approx p^{2-2\alpha - \beta}$. Let $\kappa_1=2-2\alpha-\beta$ and $\kappa_2=1-2\theta-\zeta$. Here $\kappa_1$ and $\kappa_2$ are the synchronized indexes of signal weakness and sparsity in the mean difference and the covariance matrix difference, respectively. The total information can be represented as $p^{\kappa}$, where 
\be\label{kappa}
\kappa=\max\{\kappa_1, \kappa_2\}.
\ee
With $\kappa_1$, $\kappa_2$ and $\kappa$, we have the following two theorems. 

\begin{thm}\label{thm:qdaunknownmu}
Under model (\ref{M12}) and the parameterization (\ref{Param1}),  (\ref{Omega1})--(\ref{Param3}), and (\ref{eqn:gamma}) that $\gamma > 1/2$, 
\begin{itemize}
\item [(i)] When $\theta \geq \delta/2$, i.e., the signals are weak, 
\begin{itemize}
\item[(1)] If $\kappa > (1-\delta)/2$, then MR(QDA)$\goto 0$ as $p\to\infty$. 

\item[(2)] If $\kappa< (1-\delta)/2$, then $MR(QDA) \geq c > 0$ when $p \goto \infty$.
\end{itemize}

\item[(ii)] When $\theta < \delta/2$, i.e., the signals are strong, the results in Proposition \ref{thm:idealunequal} hold.
\begin{itemize}
\item[(1)] If $\kappa > 0$, then MR(QDAfs)$\goto 0$ as $p\to\infty$. 

\item[(2)] If $\kappa< 0$, then $MR(L) \geq c > 0$ when $p \goto \infty$ for any classifier $L$.
\end{itemize}

\end{itemize}
\end{thm}


Our results show that $\kappa$ is a key quantity in the phase transition.  
The sample size have two kinds of effects when the non-zeros in $\mu$ turns from weak to strong:
\begin{itemize}
\item[(i)] When $\theta < \delta/2$, the sample size is large enough so that the signals in the mean vector can be almost perfectly recovered. With the feature selection step, the QDA achieves an asymptotic misclassification rate of 0 when $\kappa > 0$, and $1/2$ when $\kappa<0$. In addition, the latter region is proven to be a failure region for all classifiers, which is referred to as the region of impossibility. Thus the boundary $\kappa=0$ partitions the phase space into the region of possibility and impossibility.

\item[(ii)] When $\theta > \delta/2$, the sample size is insufficient for the signal recovery, and the feature selection step is ineffective. The QDA misclassification rate converges to 0 when $\kappa > (1-\delta)/2$, and a positive value when $\kappa<(1-\delta)/2$. 
Thus the boundary $\kappa= (1 - \delta)/2$ separates the regions of success and failure for QDA. 
\end{itemize}

\begin{figure}[htb!]
    \centering
    \subfigure[Strong signal region ($\theta<\delta/2$)]{\includegraphics[scale=0.4]{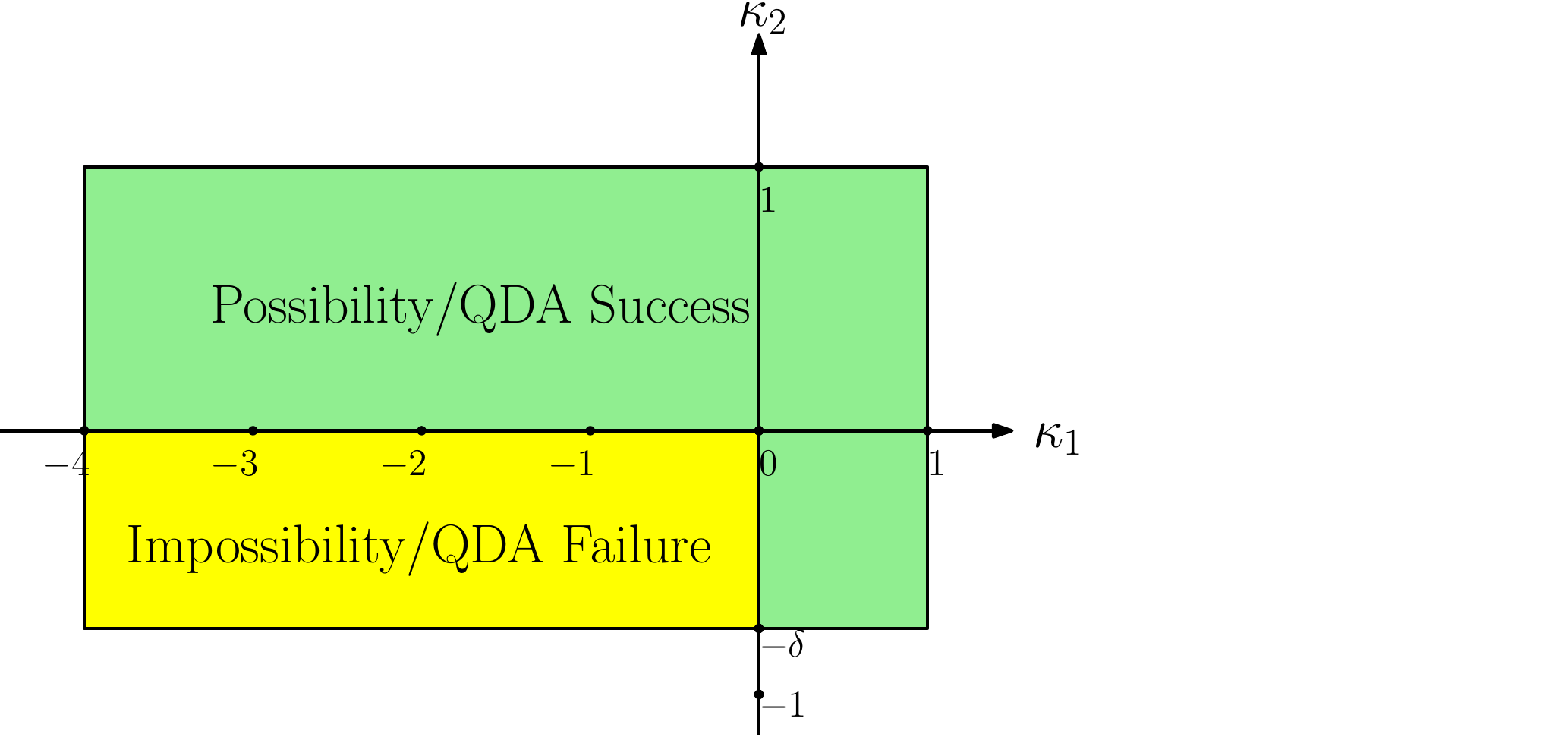}} \hspace{0.2cm} 
    \subfigure[Weak signal region ($\theta\geq \delta/2$)]{\includegraphics[scale=0.4]{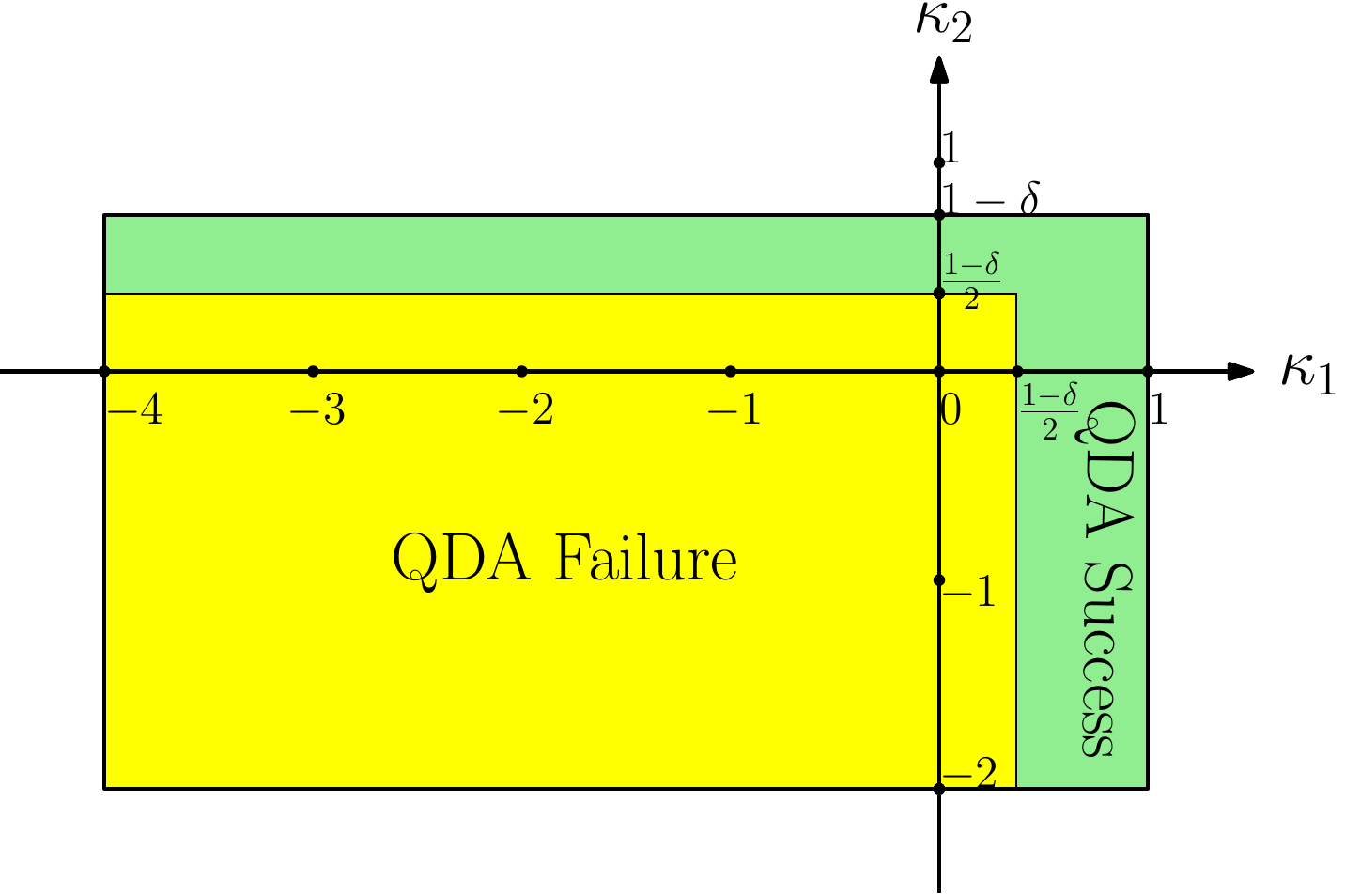}} 
      \caption{The possibility/impossibility regions and QDA successful/unsuccessful classification regions derived in Theorem \ref{thm:unknownmu} and defined in terms of $\kappa_1=2-2\alpha-\beta$ and $\kappa_2=1-2\theta-\zeta$ for fixed $\delta$ and for the two cases: (a) strong signal region ($\theta<\delta/2$) and (b) weak signal region ($\theta\geq \delta/2$).}
    \label{figure: phase}
\end{figure}

\begin{figure}[htb!]
    \centering
    \subfigure[$\theta<\delta/2$, $\kappa_2\leq 0$]{\includegraphics[scale=0.4]{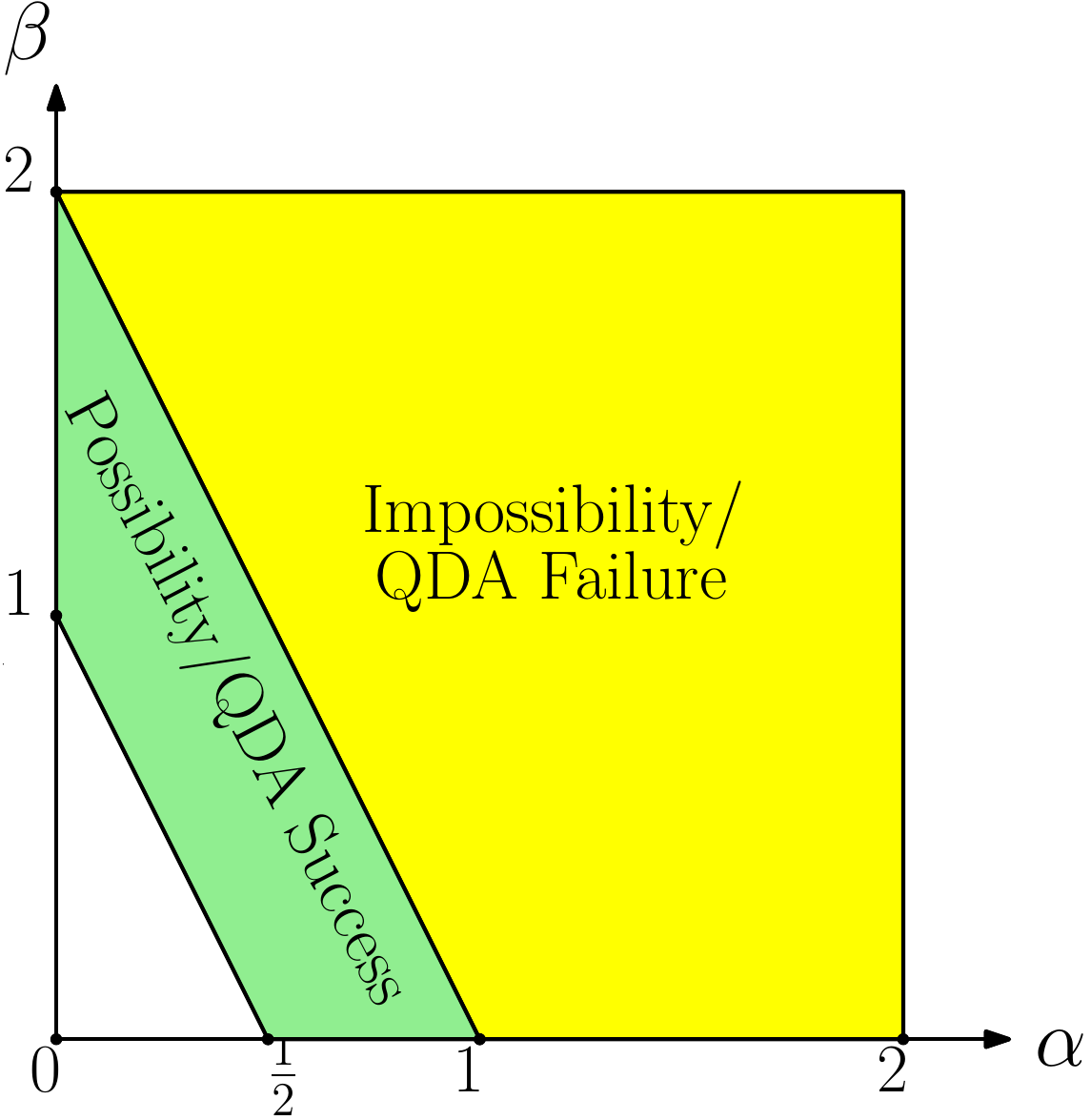}} 
               \hspace{0.2cm}
               \subfigure[$\theta\geq \delta/2$, $\kappa_2\leq (1-\delta)/2$]{\includegraphics[scale=0.4]{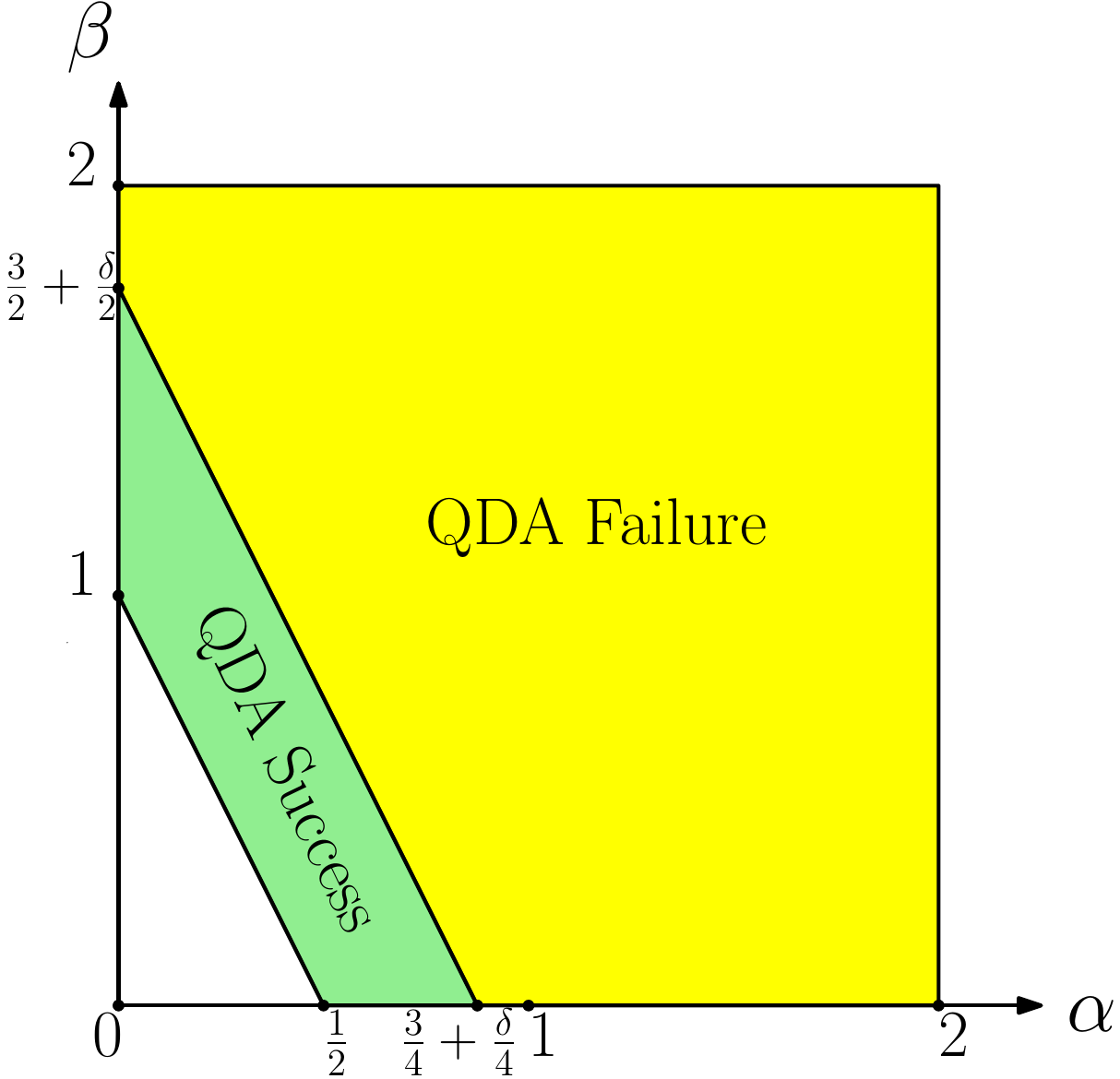}} 
    \subfigure[$0<\kappa_1\leq (1-\delta)/2$]{\includegraphics[scale=0.4]{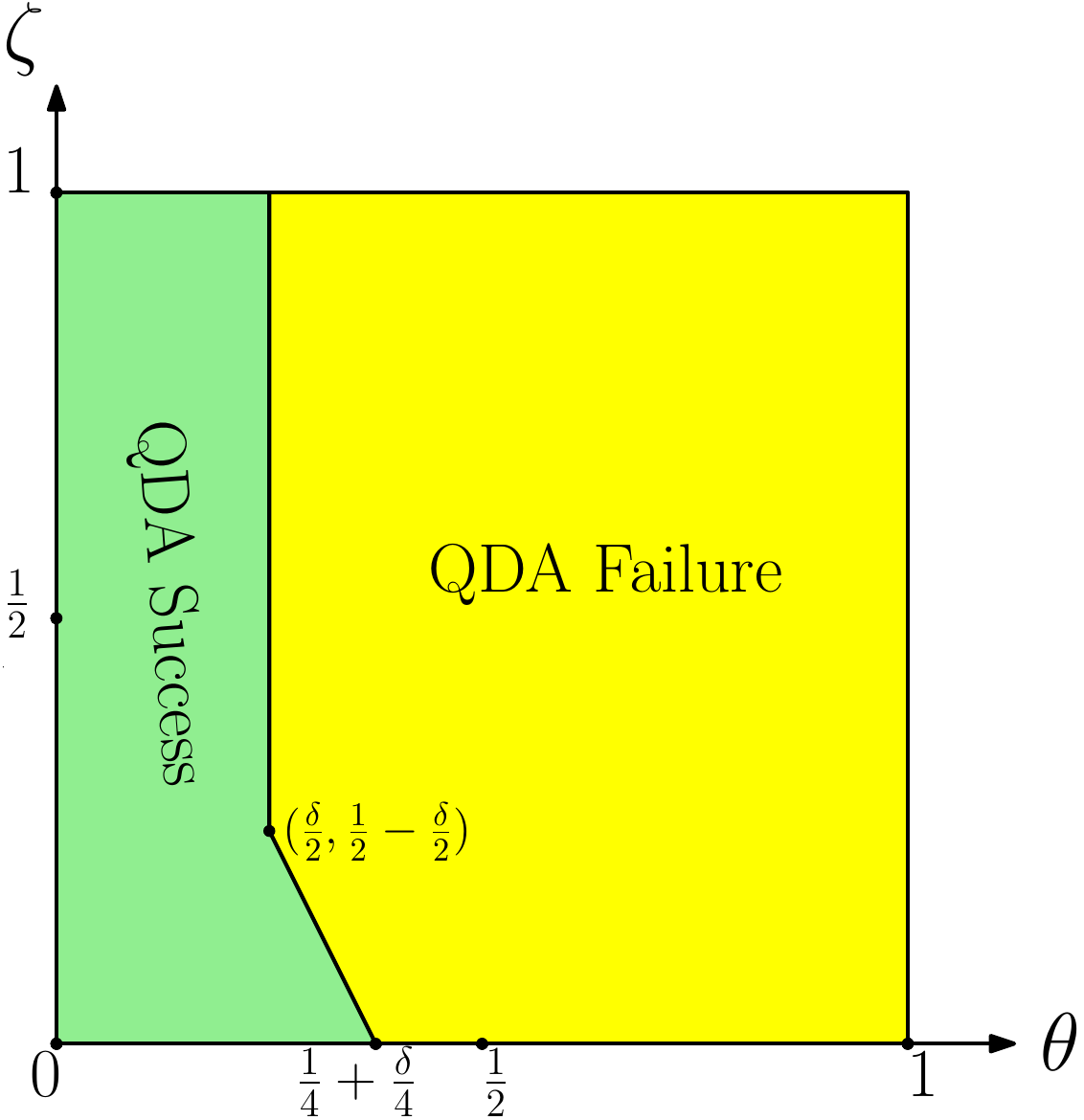}} 
                \hspace{0.7cm}
                \subfigure[$\kappa_1\leq 0$]{\includegraphics[scale=0.4]{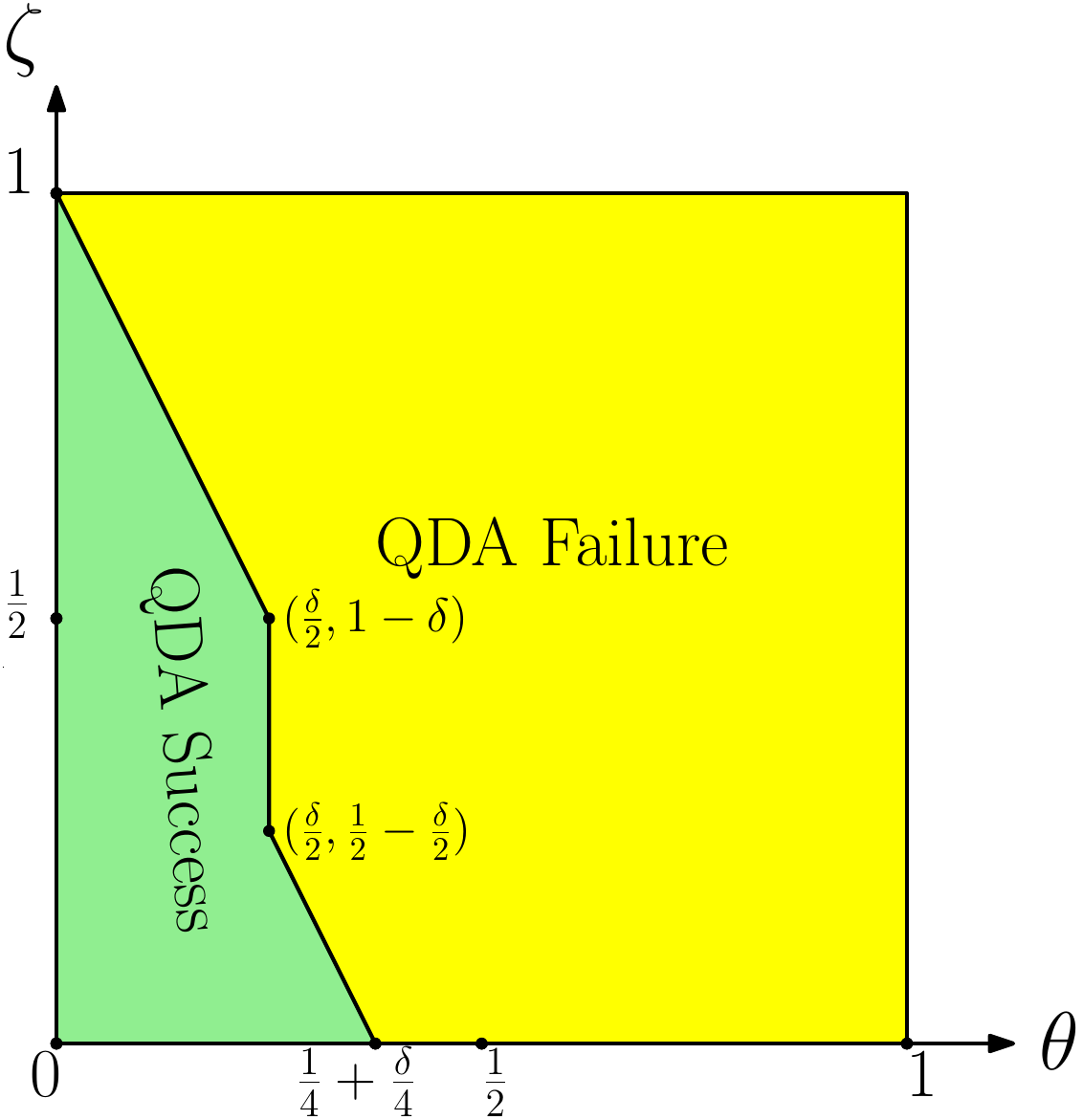}} 
      \caption{The possibility/impossibility regions and QDA successful/unsuccessful classification regions derived in Theorem \ref{thm:unknownmu} when $\delta$ and part of the rest parameters are fixed: (a) $\delta$, $\theta$ and $\zeta$ are fixed, $\theta<\delta/2$ and $1-2\theta-\zeta\leq 0$; (b) $\delta$, $\theta$ and $\zeta$ are fixed, $\theta\geq \delta/2$ and $1-2\theta-\zeta \leq (1-\delta)/2$; (c) $\delta$, $\alpha$ and $\beta$ are fixed, and $0<2-2\alpha-\beta \leq (1-\delta)/2$; (d) $\delta$, $\alpha$ and $\beta$ are fixed, and $2-2\alpha-\beta\leq 0$.}
    \label{figure: qdaparameter}
\end{figure}

Figures \ref{figure: phase} and \ref{figure: qdaparameter} below provide a visual representation of the above results.  Figure \ref{figure: phase} depicts those regions on the $\kappa_1$-$\kappa_2$ surface, with subfigure (a) for the strong signal case  and (b) for the weak signal case. In subfigure (a), the successful/failure regions of QDA are the possibility/impossibility regions for the classification problem, respectively.

Figure \ref{figure: qdaparameter} provides a sense of the relationship between the sparsity and weakness parameters of the mean and covariance matrix. 
Subfigures (a) and (b) are on the $\alpha$-$\beta$ plane about the precision matrix when $\kappa_2$ is fixed, while (c) and (d) are on the $\theta$-$\zeta$ plane about the mean vector when $\kappa_1$ is fixed. 
To better demonstrate the relationship between the parameters and the success/failure region, we only consider the cases $\kappa_2\leq 0$ in (a) and $\kappa_2\leq (1-\delta)/2$ in (b). Otherwise, the information in the mean vector is sufficient for successful classification. 
 Similarly, we do not consider the case $\kappa_1>(1-\delta)/2$ in (c) and (d). From (a) and (b) we can see that, when $\theta$ increases from less than $\delta/2$ to greater than $\delta/2$, the QDA successful region of $\alpha$ and $\beta$ decreases. As we can see from (c) and (d), the QDA success region of $\theta$ and $\zeta$ decreases when $\kappa_1$ decreases from a positive value to a negative value.

Compare the current results with Theorem 2 in the main paper. When the sample size is large to recover the signals in $\mu$, then we apply QDAfs and the bound is the same. When the sample size is small that the non-zeros in $\mu$ cannot be exactly recovered, QDA gives the bound as $\kappa > (1-\delta)/2$, which means the information is larger than $\sqrt{p/n}$. QDAw only needs $\kappa_1 > 0$ if the precision matrix $\Omega_1$ is informative. When $\mu$ is informative, QDAw requires $\theta + \zeta < 1/2$. In terms of $\epsilon$, $\tau$, $n$ and $p$, QDAw requires $\sqrt{p}\tau \epsilon \goto \infty$ while QDA requires $\tau^2 \epsilon\sqrt{np} = (\sqrt{p}\tau \epsilon)\sqrt{n}\tau \goto \infty$. 
Since $\sqrt{n}\tau \goto 0$ in the weak signal case, QDA has a stronger condition than QDAw.

\subsection{Main results when only one precision matrix is known}
Suppose $\Omega_1$ is unknown, hence we have to estimate $\Omega_1$ first by some precision matrix recovery method. 
It can be performed via any suitable approach. This has been discussed in numerous publications in the literature, such as \cite{CLIME, FJY, glasso, PCS}. Here, our goal is to develop the QDA approach with feature-selection step, instead of designing a new precision-matrix-estimation approach. In the following theorem, we consider a general precision-matrix-estimation approach, and let $\Delta_{\hat{\Omega}}$ denote the estimation error. We give the region that $MR(QDA) \goto 0$ based on $\Delta_{\hat{\Omega}}$. 



\begin{thm}\label{thm:qdaunknown}
Consider model (\ref{M12}) and the parameterization (\ref{Param1}),  (\ref{Omega1})--(\ref{con1}), (\ref{eqn:gamma}), and (\ref{Param3}). Assume $1 < \beta < 2$. For the employed precision-matrix-estimation approach, let $\Delta_{\hat\Omega} = \|\Omega - \hat{\Omega}\|$ be the spectral norm of the error. 
Suppose $\Delta_{\hat\Omega} \goto 0$ when $p \goto \infty$, and it satisfies that 
\[
\Delta_{\hat{\Omega}} (p\eta  + p \tau^2 \epsilon + \sqrt{p}\log p)+ p \Delta_{\hat{\Omega}}^2 +p/n \ll p\xi^2 + \eta^2 p^2 \nu + \tau^2 p \epsilon. 
\]
Then, QDA for the weak signal case or QDAfs for the strong signal case has a misclassification rate that converges to 0 as $p\to\infty$. 
\end{thm}

Here, we develop the general rule for the QDAfs. For the weak signal case that $\theta > \delta/2$, the condition can be relaxed by replacing $p/n$ to be $p\xi/n$. However, the term $p/n$ is not the dominating term when we apply the PCS method and the CLIME method as the precision matrix estimator. Therefore, we didn't differentiate the two cases. 

Compared to Proposition \ref{thm:idealunequal} and Theorem \ref{thm:unknownmu}, a big difference here is that the condition is an inequality that containing both the precision matrix parameters and the mean vector parameters. In the following two corollaries, we can see that the condition $\zeta < \alpha + \delta/2 - 2\theta$ indicates an intervention between the precision matrix weakness parameter $\alpha$, the mean vector sparsity and weakness parameters $\zeta$ and $\theta$, and the sample size parameter $\theta$. Under this situation, the dominating error term comes from $X^{\top}(\hat{\Omega}_1 - I)X$, which contains both $\Omega_1$ and $\mu$ (in $X$). 
 
We apply the PCS approach in \cite{PCS} and the CLIME approach in \cite{CLIME} to be the precision-matrix-estimation approach. The results can be found in the following corollaries. The boundaries they can achieve are the same. 
\begin{cor}\label{cor:PCS}
Under the conditions of Theorem \ref{thm:qdaunknown} and that $\alpha < \delta/2$, and PCS is employed for precision-matrix estimation.  Consider the conditions that 
\begin{itemize}
\item [(i)] $\beta < 1  - \alpha + \delta/2$; or
\item [(ii)] $\zeta < \alpha + \delta/2 - 2\theta$.
\end{itemize} 
If one of the above conditions is satisfied, then QDA for the weak signal case or QDAfs for the strong signal case has a misclassification rate that converges to 0 as $p\to\infty$. 
\end{cor}

\begin{cor}\label{cor:CLIME}
Under the conditions of Theorem \ref{thm:qdaunknown}, and that CLIME is employed for precision-matrix estimation.  Assume $\alpha < \delta/2$, and consider the conditions that 
\begin{itemize}
\item [(i)] $\beta < 1  - \alpha + \delta/2$;
\item [(ii)]  $\zeta < \alpha + \delta/2 - 2\theta$;
\end{itemize} 
If one of the above conditions is satisfied, then QDA for the weak signal case or QDAfs for the strong signal case has a misclassification rate that converges to 0 as $p\to\infty$. 
\end{cor}


\section{Proof of Proposition \ref{thm:idealunequal}}\label{app:ideal}
In Proposition \ref{thm:idealunequal}, we show the lower bound and upper bound of QDA given all the parameters. The proof for the lower bound is in Section \ref{sec:lowerproof} of the main paper, and here we only need to prove the upper bound by QDA. For short, we take $\Omega =\Omega_1$ without confusion in this section.

 To prove the upper bound, we want to show that in the region of possibility, $p_{i, \mu, \Omega} \goto 0$ with probability $1 - o(1)$, so that $MR(QDA) \goto 0$. 
For the ideal case, QDA estimates the label as $\hat{Y} = I\{Q > 0\}$, where $Q = Q(X, \mu, \Omega) = S - T_S$, that 
\bea\label{proofS}
S = - X^\top(\Omega - I) X + 2 \mu^{\top} (\Omega + I) X, \qquad 
T_S = \mu^\top (\Omega - I)\mu  - \log|\Omega|. 
\eea
The mis-classification rate $p_{0, \mu, \Omega} = P(S - T_S > 0)$ and $p_{1, \mu, \Omega} = P(S - T_S < 0)$. Hence, we want to find the distribution of $S$
and the magnitude of $T_S$. 

According to Lemma \ref{lemma:quad}, we have 
\begin{itemize}
\item $E[S|Y = 0] = Tr(I - \Omega) - \mu^\top(3\Omega + I) \mu$, $E[S|Y=1] = Tr(\Omega^{-1} - I) + \mu^\top(\Omega + 3I) \mu.$
\item $\var(S|Y=0) = 2Tr((\Omega - I)^2) + 16 \mu^\top \Omega^2 \mu$;
\item
$\var(S|Y = 1) = 2Tr((\Omega^{-1} - I)^2) + 16 \mu^\top \Omega^{-1} \mu$
\end{itemize}
Further, we can prove that $\frac{S - E[S|Y=i]}{\sqrt{\var(S|Y=i)}}|Y=i$ converges to standard normal distribution. 

Now, take $p_{0, \mu, \Omega}$ as an example and we have 
\bea
p_{0, \mu, \Omega} & = & P(S - T_S > 0)\nonumber\\
& = & P(\frac{Tr(I - \Omega) - \mu^\top(3\Omega + I) \mu - T_S}{\sqrt{2Tr((\Omega - I)^2) + 16 \mu^\top \Omega^2 \mu}} > 0) + o(1)\nonumber\\
& = & \Phi(\frac{4\mu^\top \Omega \mu - Tr(I - \Omega) + \log|\Omega|}{\sqrt{2Tr((\Omega - I)^2) + 16 \mu^\top \Omega^2 \mu}}) + o(1)\nonumber\\
& = & \Phi(\frac{ - 4\mu^\top \Omega \mu - \|\Omega - I\|_F^2/2}{\sqrt{2\|\Omega - I\|_F^2 + 16 \mu^\top \Omega^2 \mu}} > 0) + o(1)\nonumber\\
& = & \Phi(-\sqrt{\mu^\top \Omega \mu + \|\Omega - I\|_F^2/8}) + o(1).\nonumber
\eea
According to Lemma \ref{lemmaeb}, we have $\|\Omega - I\|_F^2 = (p^{1-2\gamma} + p^{2 -2\alpha -\beta})(1+o(1))$ and $\mu^\top \Omega \mu = p^{1 - 2\theta - \zeta}(1+o(1))$ with probability $1 - o(1)$. 

When $1 - 2\theta - \zeta > 0$ or $1 - 2\gamma > 0$ or $2 - 2\alpha - \beta > 0$, either $\|\Omega - I\|_F^2 \goto \infty$ or $\|\mu\|^2 \goto \infty$, which concludes that $p_{0,\mu,\Omega} \goto 0$. 
Similarly, we also have $p_{1,\mu, \Omega} \goto 0$ in this region.
Therefore, $MR(QDA) \goto 0$ in this region. The region of possibility is proved. 

\section{Proof of Theorem \ref{thm:real} }\label{app:real}

\subsection{Feature selection}\label{app:fs}
In this section, we want to prove that the feature selection step in QDAfs can successfully recover the signals with probability $1+o(1)$. To prove it, we start with the case that $\Omega_0 = I$, and then discuss the case that no information is given.

Consider the case that both are known. The definition and distribution of $d$ is that 
\be
d \ = \ \Omega \hat{\mu}_1 - \hat{\mu}_0 \ \sim \ N\bigl((I + \Omega)\mu, \ \frac{1}{n_0}I+\frac{1}{n_1}\Omega\bigr).
\ee 
Each entry $d_i \sim N(\mu_i + (\Omega \mu)_i, \frac{1}{n_0}+\frac{1 + \xi}{n_1})$. 
For the mean term, note that
\[
(\Omega\mu)_i = ((1+\xi) \mu + V\mu)_i = (1+\xi)\mu_i + \sum_{j \neq i, j = 1}^{p}  v_{ij}\mu_j.
\]
Under model (\ref{Param1}),  $| v_{ij}\mu_j| \leq  \eta\tau$ and $E[ v_{ij}\mu_j] = 0$. The magnitude $\eta$ is to measure the strength and not the number of non-zeros. So, let $\sigma^2 = \sum_{j \neq i} \frac{1}{\eta^2}E[ v_{ij}^2 \mu_j^2] = (p-1)  \tau^2 \epsilon \nu$; then, according to Bennett's inequality \citep{Bennett1962}, 
\[
P\biggl(\bigl|\sum_{j \neq i, j = 1}^p  v_{ij}\mu_j\bigr| > \eta t\biggr) \leq 2 \exp\left[-\frac{\sigma^2}{\tau^2} \biggl((1+\frac{\tau t}{\sigma^2})\log(1+\frac{\tau t}{\sigma^2}) - \frac{\tau t}{\sigma^2}\biggr)\right].
\]
So, with probability $1 - O(p^{-2})$, $|\sum_{j \neq i}  v_{ij}\mu_j| \leq 2 \eta\tau (\sqrt{p \epsilon \nu} + 1)\log p$. Under  (\ref{con1}), $2 \eta \tau (\sqrt{p \epsilon \nu} + 1)\log p \ll (1 + \xi)\tau = (1+\xi)\mu_i$ when $\mu_i \neq 0$. Hence, the feature-selection step would depend on $(1+\xi)\mu_i = \mu_i(1 + o(1))$.

We consider the weak signal case in which $\theta > \delta/2$ and the strong signal case in which $\theta < \delta/2$. In the former case, the signal strength $\tau = p^{-\theta} \ll L_p /\sqrt{n}$. Hence, $\max_{1 \leq j \leq p} |d_i| \leq 2{\log p}/\sqrt{n}$ with probability $1 - o(1)$, and the threshold will be 0. 
For the latter case, in which $\theta< \delta/2$, $\tau \gg 1/\sqrt n$. Hence, with probability $1-o(1)$, we have $t = 2\sqrt{\log p}/\sqrt{n}$, and the set $\{i: \mu_i \neq 0\}$ is recovered with zero error. 

\vspace{2em}
Now it comes to the case that $\Omega_0 = I$ and $\Omega_1 = \Omega$ is unknown. We estimate $\Omega$ with $\hat{\Omega}$ and $d$ by $\hat{d}$, that 
\be
\hat{d} = \hat{\Omega} \hat{\mu}_1 - \hat{\mu}_0 \ \sim \ N\bigl((I + \hat\Omega)\mu, \ \frac{1}{n_0}I+\frac{1}{n_1}\hat\Omega\Omega^{-1}\hat\Omega\bigr).
\ee 
Therefore, for each entry $\hat{d}_i$, it differs from $d_i$ by the mean and variance. Recall that to assure successful recovery of $\Omega$, it is required the number of non-zeros in each row is no larger than $\sqrt{n}$ and $\eta \gg 1/\sqrt{n}$. 
Under these two conditions, $\hat{\Omega}$ recovered all the non-zeros of $\Omega$ with probability $1-o(1/p)$ and $\max_{i,j}|\Omega(i,j) - \hat{\Omega}(i,j)| \leq \frac{\log p}{\sqrt{n}}$. Therefore, the mean differs at the order that 
\[
\|(\hat\Omega - \Omega)\mu\|_{\infty} 
\leq \|\hat\Omega - \Omega\|_1\|\mu\|_{\infty}
= o(\tau).
\]
For the variance, the maximal difference is that 
\[
\frac{1}{n_1}\|diag(\hat\Omega\Omega^{-1}\hat\Omega - \Omega)\|_{\infty} \leq \frac{1}{n_1}\|\hat\Omega\Omega^{-1}\hat\Omega - \Omega\| = o(1/n).
\]
Therefore, the difference between $d$ and $\hat{d}$ is a second order term compared to $d$. The results still hold. 

\vspace{2em}
When both $\Omega_0$ and $\Omega_1$ are unknown, the analysis is similar so we ignore it here. 

As a conclusion, the effectiveness of feature selection can be proved. 

\subsection{Proof of Theorem \ref{thm:real}}\label{app:realproof}
For the case that the precision matrices are unknown but there is information that the diagonals are around 1, we estimate the precision matrix $\tilde{\Omega}_k$ by finding the PCS estimate $\hat{\Omega}_k$ first and then adjust the diagonals to be close to 1. It reduces a large amout of noise. 

To prove the main theorem, we first consider the case that both $\Omega_0$ and $\Omega_1$ are estimated by PCS but $\mu$ is known. Then, we take into consideration that how the estimation of $\mu$ will affect the result. 

When $\mu$ is known and $\Omega_i$'s are estimated by PCS, we classify by $\hat{Y} = I(Q(X, \mu, \hat{\Omega}) > 0)$, where 
\[
Q(X, \mu, \tilde{\Omega}) = X^\top (\tilde{\Omega}_0 - \tilde{\Omega}_1)X + 2\mu^\top(\tilde{\Omega}_0 + \tilde{\Omega}_1)X - T_S.
\]
Here, $\tilde{\Omega}_k = \hat{\Omega}_k - D_k + I$, where $D_k$ is the diagonal matrix formed by the diagonals of $\hat{\Omega}_k$. Therefore, $\tilde{\Omega}_k$ is to force all the diagonals of $\hat{\Omega}_k$ to be 1. 
The constant term $T_S$ is given by 
\be 
T_S = \mu^\top(\tilde{\Omega}_1 - \tilde{\Omega}_0)\mu - \log|\tilde{\Omega}_1| + \log|\tilde{\Omega}_0|.
\ee

Let $Q(X, \mu, \hat{\Omega}) = S^{PCS} - T_S$, where $S^{PCS} = X^\top (\tilde{\Omega}_0 - \tilde{\Omega}_1)X + 2\mu^\top(\tilde{\Omega}_0 + \tilde{\Omega}_1)X$. 
Note that $X$ and $\tilde{\Omega}_i$'s are independent. Given $\tilde{\Omega}_1$ and $\tilde{\Omega}_0$, we derive the asymptotic distribution of $S^{PCS}$ by Lemma \ref{lemma:quad}. 
In details, the expectations and variances are 
\begin{itemize}
    \item $E[S^{PCS}|Y = 0] = T_S -4\mu^\top \tilde\Omega_1 \mu  -  \log|\tilde\Omega_0| + \log|\tilde\Omega_1| + Tr(\Omega_0^{-1}(\tilde\Omega_0 - \tilde\Omega_1))$;
    \item 
    $E[S^{PCS}|Y = 1] = T_S + 4\mu^\top \tilde\Omega_0 \mu  -  \log|\tilde\Omega_0| + \log|\tilde\Omega_1| + Tr(\Omega_1^{-1}(\tilde\Omega_0-\tilde{\Omega}_1))$; 
    \item $\var(S^{PCS}|Y = 0) = 2Tr((\tilde\Omega_0\Omega_0^{-1} - \tilde\Omega_1\Omega_0^{-1})^2) + 4 \mu^\top (3\tilde\Omega_1\Omega_0^{-1}\tilde\Omega_1 + \tilde\Omega_0\Omega_0^{-1}\tilde\Omega_0) \mu$;
    \item $\var(S^{PCS}|Y = 1) = 2Tr((\tilde\Omega_0\Omega_1^{-1} - \tilde\Omega_1\Omega_1^{-1})^2) + 4 \mu^\top (3\tilde\Omega_0\Omega_1^{-1}\tilde\Omega_0 + \tilde\Omega_1\Omega_1^{-1}\tilde\Omega_1) \mu$. 
\end{itemize}
Define $Z_i = [S^{PCS} - E[S^{PCS}|Y = i]]/\sqrt{\var(S^{PCS}|Y = i)}$ and $F_{Z_i}(x) = P(Z_i \leq x)$, then $\sup\nolimits_{-\infty < x < \infty}|F_{Z_i}(x) - \Phi(x)| \goto 0$.

We consider the case $Y = 0$, where we have to derive the asymptotic results for $p_{0,\mu,\Omega}$. According to the results above, we have 
\be\label{eqn:pf4p0}
p_{0,\mu,\Omega} = \Phi(\frac{-4\mu^\top \tilde\Omega_1 \mu  -  \log|\tilde\Omega_0| + \log|\tilde\Omega_1| + Tr(\Omega_0^{-1}(\tilde\Omega_0 - \tilde\Omega_1))}{\sqrt{2Tr((\tilde\Omega_0\Omega_0^{-1} - \tilde\Omega_1\Omega_0^{-1})^2) + 4 \mu^\top (3\tilde\Omega_1\Omega_0^{-1}\tilde\Omega_1 + \tilde\Omega_0\Omega_0^{-1}\tilde\Omega_0) \mu}}) + o(1).
\ee
To find the approximation of $p_{0,\mu,\Omega}$, we need proper approximation of the fraction inside. 

According to Theorem 2.3 in \cite{PCS}, when $1 - \delta/2 < \beta < 2$ and $\eta \gg 1/\sqrt{n}$, PCS recovers the exact support with probability $1 - o(1/p^2)$, and 
$\max_{i,j}|\Omega_1(i,j) - \hat{\Omega}_1(i,j)| \leq C\sqrt{\log p/n}$. In this case, with probability $1-o(1/p^2)$, the number of non-zero entries in each row of $\Omega_0$ and $\Omega_1$ is uniformly bounded by a constant. Therefore, we have the bound on the spectral norm of the estimator: 
\[
\|\Omega_0 - \hat{\Omega}_0\| \leq C\sqrt{\log p/n}, \quad
\|\Omega_1 - \hat{\Omega}_1\| \leq C\sqrt{\log p/n}.
\]
Because $\Omega_0$ and $\Omega_1$ are close to the identity matrix, so we have 
\[
\|\hat{\Omega}_0 - I\| = o(1), \quad 
\|\hat{\Omega}_1 - I\| = o(1).
\]
It means all the eigenvalues of $\hat{\Omega}_0$ are in the interval $1 - o(1) \leq \lambda_p(\hat{\Omega}_0) \leq \lambda_1(\hat{\Omega}_0) \leq 1 + o(1)$, and $1 - o(1) \leq \lambda_p(\hat{\Omega}_1) \leq \lambda_1(\hat{\Omega}_1) \leq 1 + o(1)$. 
A subsequent result is that 
\[
\|\hat{\Omega}_0^{-1/2}\hat{\Omega}_1\hat{\Omega}_0^{-1/2} - I\| = o(1).
\]
Finally, according to the definition of $\tilde{\Omega}_k$ and the condition that $\gamma > 1/2$, the error between on the diagonals of $\tilde{\Omega}_k - \Omega_k$ is $\xi_p$, which is smaller than that of the diagonals of $\hat{\Omega}_k - \Omega_k$ at $\sqrt{\log p/n}$. Therefore, all the above conclusions hold for $\tilde{\Omega}_k$, $k = 0, 1$. 

Now we derive the numerator. Since $\|\tilde{\Omega}_1 - I\| = 1+o(1)$, there is 
\[
-4\mu^\top \tilde\Omega_1 \mu = -4\|\mu\|^2(1+o(1)).
\]
We decompose 
\bea
Tr(\Omega_0^{-1}(\tilde{\Omega}_0 - \tilde{\Omega}_1)) & = &Tr(\tilde\Omega_0^{-1}(\tilde{\Omega}_0 - \tilde{\Omega}_1)) + Tr((\Omega_0^{-1} - \tilde{\Omega}_0^{-1})(\tilde{\Omega}_0 - \tilde{\Omega}_1)).
\eea
Then for the first term, we have 
\bea
&&Tr(\tilde\Omega_0^{-1}(\tilde{\Omega}_0 - \tilde{\Omega}_1)) -  \log|\Omega_0| + \log|\Omega_1|\nonumber\\
&=&Tr(I - \tilde{\Omega}_0^{-1/2}\tilde{\Omega}_1\tilde{\Omega}_0^{-1/2}) + \log|\tilde{\Omega}_0^{-1/2}\tilde{\Omega}_1\tilde{\Omega}_0^{-1/2}|\\
& = & -\frac12Tr((I - \tilde{\Omega}_0^{-1/2}\tilde{\Omega}_1\tilde{\Omega}_0^{-1/2})^2)(1+o(1))\\
& = & -\frac12Tr(\tilde{\Omega}_0^{-1}(\tilde{\Omega}_0 - \tilde{\Omega}_1)^2\tilde{\Omega}_0^{-1})(1+o(1))\\
& = & -\frac12Tr((\tilde{\Omega}_0 - \tilde{\Omega}_1)^2)(1+o(1)) = -\frac12\|\tilde{\Omega}_0 - \tilde{\Omega}_1\|_F^2(1+o(1)),
\eea
where the last equality comes from Lemma \ref{lemma:trace}. 
For the second term, the derivation is similar with the main paper. 
Let $\Delta V  = Tr((\Omega_0^{-1}-\tilde\Omega_0^{-1})(\tilde{\Omega}_0 - \tilde{\Omega}_1))$. 
For any square matrices $A$ and $B$ with ordered singular values as $\alpha_i$ and $\beta_i$, respectively. By Von Neuman's trace inequality, $|Tr(AB)| \leq \sum \alpha_i \beta_i \leq \sqrt{\sum \alpha_i^2 \sum \beta_i^2} = \sqrt{Tr(A^\top A)}\sqrt{Tr(B^\top B)}$. 
Apply this result to $\Delta V$ and recall that both $\Omega_0$ and $\tilde{\Omega}_0$ has eigenvalues at $1+o(1)$. Then we have 
\begin{eqnarray*}
|\Delta V| & \leq & \sqrt{Tr((\Omega_0^{-1} - \tilde{\Omega}_0^{-1})^2)} \sqrt{Tr((\tilde\Omega_0 - \tilde{\Omega}_1)^2)}\\
&\leq & \sqrt{Tr(\Omega_0^{-2}(\Omega_0 - \tilde{\Omega}_0)^2\tilde{\Omega}_0^{-2})} \|\tilde\Omega_0 - \tilde{\Omega}_1\|_F(1+o(1)).\\
& \leq & \|\tilde{\Omega}_0 - \Omega_0\|_F \|\tilde\Omega_0 - \tilde{\Omega}_1\|_F(1+o(1)).
\end{eqnarray*}
As a conclusion, the numerator is
\be\label{eqn:pfvarnum}
-4\|\mu\|^2(1+o(1)) -\frac12\|\tilde{\Omega}_0 - \tilde{\Omega}_1\|_F^2(1+o(1)) + O(\|\tilde{\Omega}_0 - \Omega_0\|_F \|\tilde\Omega_0 - \tilde{\Omega}_1\|_F).
\ee

Now we consider the denominator. 
Consider the first term, since $\|\Omega_0 - I\| = o(1)$, 
\[
Tr((\tilde\Omega_0\Omega_0^{-1} - \tilde\Omega_1\Omega_0^{-1})^2) = 
Tr((\tilde\Omega_0 - \tilde\Omega_1)^2)(1+o(1))
=\|\tilde\Omega_0 - \tilde\Omega_1\|_F^2(1+o(1))
\]
Consider the second term. Since all $\Omega_i$ and $\tilde\Omega_i$ have eigenvalues around 1, so we have 
\[
\mu^\top (3\tilde\Omega_1\Omega_0^{-1}\tilde\Omega_1 + \tilde\Omega_0\Omega_0^{-1}\tilde\Omega_0) \mu = 4\|\mu\|^2(1+o(1)). 
\]
As a summary, the denominator is
\be\label{eqn:pfvarden}
\sqrt{2\|\tilde\Omega_0 - \tilde\Omega_1\|_F^2(1+o(1)) + 16\|\mu\|^2(1+o(1))}.
\ee

Introduce (\ref{eqn:pfvarnum}) and (\ref{eqn:pfvarden}) into $p_{0,\mu,\Omega}$, there is 
\begin{eqnarray}\label{eqn:pfvarp0frac}
p_{0, \mu, \Omega} & = &\Phi(\frac{T_S^{PCS} - E[S^{PCS}|Y = 0]}{\sqrt{\var(S^{PCS}|Y = 0)}}) \nonumber\\
&= &\Phi(\frac{-4\|\mu\|^2 -\frac12\|\tilde{\Omega}_0 - \tilde{\Omega}_1\|_F^2 + O(\|\tilde{\Omega}_0 - \Omega_0\|_F \|\tilde\Omega_0 - \tilde{\Omega}_1\|_F)}{\sqrt{2\|\tilde\Omega_0 - \tilde\Omega_1\|_F^2 + 16\|\mu\|^2}}
) + o(1)\\
&= &\Phi(-\sqrt{\|\tilde\Omega_0 - \tilde\Omega_1\|_F^2/8 + \|\mu\|^2} + O(\frac{\|\tilde{\Omega}_0 - \Omega_0\|_F \|\tilde\Omega_0 - \tilde{\Omega}_1\|_F}{\sqrt{2\|\tilde\Omega_0 - \tilde\Omega_1\|_F^2 + 16\|\mu\|^2})}
) + o(1).
\end{eqnarray}

According to Theorem 2.3 in \cite{PCS}, under current conditions, PCS recovers the exact support with probability $1 - o(1/p^2)$, and 
$\max_{i,j}|\Omega_k(i,j) - \tilde{\Omega}_k(i,j)| \leq C\sqrt{\log p/n}$, $k = 0,1$. Further, we force all the diagonals to be 1 in $\tilde{\Omega}_k$. Therefore, we have 
\[
\|\Omega_k - \tilde{\Omega}_k\|_F^2 = o(\|\Omega_k - I\|_F^2 + 1).
\]
Since $\tilde{\Omega}_k$ both have identity diagonals and $\gamma > 1/2$, we have 
\[
\|\tilde\Omega_0 - \tilde{\Omega}_1\|_F^2 = \|\Omega_0 - {\Omega}_1\|_F^2(1+o(1)) + o(1). 
\]
Further, since $\gamma > 1/2$, $\|{\Omega}_0 - \Omega_1\|_F = \frac{1}{2}\|{\Omega}_0 - I\|_F(1+o(1))$, so $\|\tilde{\Omega}_0 - \Omega_0\|_F = o(\|\Omega_1 - \Omega_0\|_F + 1)$. 
Therefore, we have 
\be\label{eqn:pfvaromega01}
\|\tilde{\Omega}_0 - \Omega_0\|_F \|\tilde{\Omega}_0 - \tilde\Omega_1\|_F 
= o(\|{\Omega}_0 - \Omega_1\|^2_F ).
\ee

To make sure the fraction in (\ref{eqn:pfvarp0frac}) goes to negative infinity, we need 
\[
\|\tilde{\Omega}_0 - \Omega_0\|_F \|\tilde{\Omega}_0 - \tilde\Omega_1\|_F \ll 4\|\mu\|^2 + \|\tilde{\Omega}_0 - \tilde\Omega_1\|_F^2/2, \mbox{ and }
\|\mu\|^2 + \|\tilde{\Omega}_0 - \tilde\Omega_1\|_F^2/8 \goto \infty. 
\]
The first inequality always holds by (\ref{eqn:pfvaromega01}). For the second inequality, by Lemma \ref{lemmaeb}, we can see that when 
\[
\eta^2p^2\nu \gg 0, \mbox{ or } 
\tau^2p\epsilon \gg 0,
\]
there is $p_{0, \mu, \Omega} \goto 0$. 
 
The similar derivation works for $p_{1,\mu,\Omega}$. Hence, we can see $p_{i,\mu,\Omega} \goto 0$ when a) $2 - 2\alpha - \beta > 0$, or b) $1 - 2\theta - \zeta > 0$. 

Now we introduce in the randomness of $\mu$. Suppose $\tau_p \gg 1/\sqrt{n}$, therefore the signals in $\mu$ are individually strong enough for successful recovery. 
We estimate $\Omega_k$ by PCS, then threshold on $d = \hat{\Omega}_1\hat{\mu}_1 - \hat{\Omega}_0\hat{\mu}_0$. QDA is applied to the post-selection data. In Section \ref{app:fs}, it is shown that the signals can be exactly recovered with probability $1 - o(1)$. 
Hence, we only consider the event that $\{t = \sqrt{2\log p/n}\}$ and all the signals are exactly recovered.

In previous analysis, we analyze the performance of $Q(X, \mu, \hat{\Omega}) = S^{PCS} - T_S$. In QDAfs, the criteria is updated as 
\be
Q(X, \hat{\mu}, \tilde{\Omega})  =
Q(X, \mu, \tilde{\Omega}) + \Delta Q, 
\ee
where $\Delta Q =  \displaystyle 2(\hat{\mu}_d^{(t)} - (\tilde{\Omega}_0 + \tilde\Omega_1)\mu)^\top X+\bigl[(\hat{\mu}_0^{(t)})^\top(\tilde{\Omega}_0-\tilde\Omega_1)\hat{\mu}_0^{(t)}-\mu^\top(\tilde{\Omega}_0-\tilde\Omega_1)\mu\bigr]$.
The following lemma bounds $|\Delta Q|$. 
\begin{lemma}\label{lemma:unknownstrongQ2}
Under the model assumptions and the definition of $\Delta Q$, there is
\be
 |\Delta Q|   \leq  \eta \tau \max\{p\epsilon \nu, 1\} \log p +  O_p(\sqrt{4n^{-1}p\epsilon}).
\ee
\end{lemma}

Combining Lemma \ref{lemma:unknownstrongQ2} with (\ref{eqn:pfvarp0frac}), the errors are 
\bea
p_{i, \mu, \Omega} & = & P((-1)^i*(S^{PCS} - T_S + \Delta Q) > 0)\nonumber\\
& = & \Phi(\frac{(-1)^i*(T_S^{PCS} - E[S^{PCS}|Y = i]) }{\sqrt{\var(S^{PCS}|Y = i)}} + \frac{(-1)^i*\Delta Q}{\sqrt{\var(S^{PCS}|Y = i)}}) + o(1), \quad i = 0, 1.
\eea
The first term is identified in (\ref{eqn:pfvarp0frac}). 
The second term can be bounded by
\[
\frac{|\Delta Q|}{\sqrt{\var(S^{PCS}|Y = i)}} 
\leq 
\frac{ \eta \tau \max\{p\epsilon \nu, 1\} \log p +  O_p(\sqrt{4n^{-1}p\epsilon})}{\sqrt{\eta^2p^2\nu/8 + \tau^2p\epsilon}(1 + o(1))}.
\]

Therefore, in the region of possibility identified by Theorem \ref{thm:real}, MR(QDAfs-PCS) converges to 0. \qed

\subsection{Proof of Lemma \ref{lemma:unknownstrongQ2}}
\begin{lemma*}
Under the model assumptions and the definition of $\Delta Q$, there is
\be
 |\Delta Q|   \leq  \eta \tau \max\{p\epsilon \nu, 1\} \log p +  O_p(\sqrt{4n^{-1}p\epsilon}).
\ee
\end{lemma*}
{\it Proof}. 
Recall that $\hat{\mu}^{(t)} = \hat{\mu}_0 \circ d^{(t)}$ and $\hat{\mu}_d^{(t)} = d \circ d^{(t)}$. 
For simplicity, in this section, we use $\hat{\mu}_0$ and $d$ to denote $\hat{\mu}_0^{(t)}$ and $\hat{\mu}_d^{(t)}$, respectively. Since all the signals are exactly recovered, $\hat{\mu}_0$ and $d$ have zeros on the non-signal entries and non-zeros on the signals.

Let $k=\|\mu\|_0$ denote the number of non-zeros in $\mu$. 
Without loss of generality, we permute $\mu$ such that the first $k$ entries are the non-zeros and the rest are the zeros. Permute $W$, $\Omega$, $\hat\Omega$ and $X$ accordingly, and rewrite $W = V/\eta$ and $\Omega$ as $2\times 2 $ block matrices $W=\left(^{W_{11} \ W_{12}}_{W_{21} \ W_{22}}\right)$ and $\Omega=\left(^{\Omega_{11} \ \Omega_{12}}_{\Omega_{21} \ \Omega_{22}}\right)$, where $W_{11}$ and $\Omega_{11}$ are $k\times k$ sub-matrices of $W$ and $\Omega$,  respectively. Let $X^{(k)}$, $d^{(k)}$, $\mu^{(k)}$, and $\hat{\mu}_0^{(k)}$ denote, respectively, $X$, $d$, $\mu$, and $\hat{\mu}_0$ restricted on the first $k$ entries, and let $X^{(p-k)}$ denote $X$ restricted on the last $(p-k)$ entries. Then $\mu^{(k)}$ is a length $k$ vector with all elements as $\tau$. 

With all the notations, $\Delta Q$ is 
\be\label{Decomp2appunknown2}
\ba{rl}
 \Delta Q = & \displaystyle 2(d - (\tilde{\Omega}_0 + \tilde\Omega_1)\mu)^\top X+\left[\hat{\mu}_0^\top(\tilde{\Omega}_0-\tilde\Omega_1)\hat{\mu}_0-\mu^\top(\tilde{\Omega}_0-\tilde\Omega_1)\mu\right] \\
 = & 2I_k + II_k.
\ea\ee

Now we analyze $I_k$ and $II_k$. The discussion focuses on the case $Y = 0$, i.e. $X \sim N(-\mu, \Omega_0^{-1})$. The derivation for $Y = 1$, i.e., $X \sim N(\mu, \Omega_1^{-1})$ is similar and the results are at the same order.  The result will include $k$. Recall that $k$ is the number of non-zeros in $\mu$, where $k \sim Binomial(p, \epsilon)$. According to Bernstein's inequality, 
\begin{eqnarray*}
  P(|k - p\epsilon| \geq \sqrt{p\epsilon} \log p)  &\leq&
 2\exp\{-\frac{(\sqrt{p\epsilon} \log p)^2/2}{p\epsilon(1-\epsilon) + (\sqrt{p\epsilon} \log p)/3}\} = o(p^{-1}).
\end{eqnarray*}
Since $p\epsilon \goto \infty$, with probability $1-o(1)$, we have $k =  p\epsilon(1+o(1))$.

\begin{itemize}
\item We consider $I_k$ first. Since $d = \left(^{d^{(k)}}_{0_{p-k}}\right)$, $\mu = \left(^{\mu^{(k)}}_{0_{p-k}}\right)$, and $X = \left(^{X^{(k)}}_{X^{(p-k)}}\right)$, where $0_{p-k}$ is a zero vector with length $p - k$,
\ben
\ba{lll}
I_k & = & \left({d^{(k)}}^\top \,\, {0_{p-k}}^\top \right)  
\left(\ba{l}
X^{(k)}\\
X^{(p-k)}
\ea\right) - 
\left({\mu^{(k)}}^\top \,\, {0_{p-k}}^\top \right)  (\tilde\Omega_0 + \tilde\Omega_1)
\left(\ba{l}
X^{(k)}\\
X^{(p-k)} 
\ea \right) \\
& = & (d^{(k)})^\top X^{(k)} - (\mu^{(k)})^\top( (\tilde\Omega_0)_{11} +  (\tilde\Omega_1)_{11}) X^{(k)} - (\mu^{(k)})^\top ((\tilde\Omega_0)_{12} + (\tilde\Omega_1)_{12}) X^{(p-k)}\\
& = &  (d^{(k)} - ((\tilde\Omega_0)_{11} + (\tilde\Omega_1)_{11}) \mu^{(k)})^\top X^{(k)} - (\mu^{(k)})^\top ((\tilde\Omega_0 + \tilde\Omega_1)_{12}) X^{(p-k)}\\
& = & Ia + Ib. 
\ea
\een

Consider $Ia$ first. Let $\tilde{d} = d^{(k)} - ((\tilde\Omega_0)_{11} + (\tilde\Omega_1)_{11})\mu^{(k)}$, then $Ia = \tilde{d}^\top X^{(k)}$, and 
\[
\tilde{d} \sim N(0, \frac{1}{n_0}(\tilde\Omega_0\Omega_0^{-1}\tilde\Omega_0)_{11} + \frac{1}{n_1} (\tilde\Omega_1\Omega_1^{-1}\tilde\Omega_1)_{11}), \quad 
X^{(k)} \sim N(-\mu^{(k)}, (\Omega_0^{-1})_{11}). 
\]
$\tilde{d}$ is independent with $X^{(k)}$, so $E[Ia] = 0$. 
The variance can be obtained by the law of total variance, that 
\ben\ba{lll}
\var(Ia)  
& \leq &\displaystyle (\frac{1}{n_0} + \frac{1}{n_1})k(1+o(1)) + \frac{1}{n_0}\|(\tilde\Omega_0\Omega_0^{-1}\tilde\Omega_0)_{11}\|k\tau^2 + \frac{1}{n_1}\|(\tilde\Omega_1\Omega_1^{-1}\tilde\Omega_1)_{11}\| k\tau^2 \\
& \lesssim &\displaystyle   \ 4n^{-1}k = 4p\epsilon/n(1+o(1)),
\ea
\een
where the trace of $(\tilde\Omega_0\Omega_0^{-1}\tilde\Omega_0)_{11}$ is constrained by $k\|(\tilde\Omega_0\Omega_0^{-1}\tilde\Omega_0)_{11}\| = k(1+o(1))$, and the same for the case with $\Omega_1$. 
For the case in which $X \sim N(\mu, \Omega^{-1})$, the same result is obtained.

We also prove the aymptotic normality according to Lemma \ref{lemma2} and the Berry-Ess\'een theorem. Therefore, $\sup_x|F_{Ia/\sqrt{\var(Ia)}}(x)-\Phi(x)|\overset{\mathcal P}\to 0$, and, hence, $Ia = O_p(\sqrt{4n^{-1}p\epsilon})$. 

Next, consider $Ib$.  Recall that $X^{(p-k)} \sim N(0, I)$. Therefore, 
\[
Ib = (\mu^{(k)})^\top ((\tilde\Omega_0 + \tilde\Omega_1)_{12}) X^{(p-k)} \sim N(0, (\mu^{(k)})^\top ((\tilde\Omega_0 + \tilde\Omega_1)_{12})(\Omega_0^{-1})_{11}((\tilde\Omega_0 + \tilde\Omega_1)_{12})^\top \mu^{(k)}). 
\]
Consider the variance term. ${\Omega}_0$ and ${\Omega}_1$ are independent. Further, the non-zeros are very sparse that the probability that $\Omega_0$ and $\Omega_1$ have the same non-zero element is a relatively smaller order term. 
Hence, $((\Omega_0 + \Omega_1)_{12})$ can be seen as $\Omega_{12}$ that follows the same model 
where the sparsity parameter is $2\nu$. According to the property of PCS estimator, with high probability, $\tilde{\Omega} = \tilde\Omega_0 + \tilde\Omega_1$ have the same non-zero off-diagonals with $\Omega$. We use $\tilde{\Omega}$ for short. 

For the variance term, since $\tilde\Omega_{12} = \tilde\Omega_{21}^\top $, we have $(\mu^{(k)})^\top\tilde\Omega_{12} (\Omega^{-1}_0)_{11}\tilde\Omega_{21} \mu^{(k)} \leq \tau^2 \|\tilde\Omega_{21}\|_\infty^2(1+o(1))$. 
Currently we require there are $o(\sqrt{n})$ non-zero entries in each row of $\Omega$ and $\eta \gg 1/\sqrt{n}$. Further, the distribution on non-zeros in $\Omega$ are independent with the non-zeros in $\mu$. Hence, with probability $1 - o(1)$,  
$\|\tilde\Omega_{21}\|_\infty^2 \leq 2\eta^2 \max\{4k^2\nu^2, 1\}$. 
As a result, with probability $1-o(1)$, 
\be
(\mu^{(k)})^\top\tilde\Omega_{12} (\Omega^{-1}_0)_{11}\tilde\Omega_{21} \mu^{(k)} \leq  \tau^2 \|\tilde\Omega_{21}\|_\infty^2(1+o(1)) \leq 9\tau^2 \eta^2 \max\{k^2\nu^2, 1\}.
\ee
So, with probability $1 - o(1)$, 
\be\label{thm3case2part1unknown2}
|Ib| \leq C\eta \tau \max\{k\nu, 1\} \log p= C\eta \tau \max\{p\epsilon\nu, 1\} \log p.
\ee
For the case in which $X \sim N(\mu, \Omega_1^{-1})$, the analysis is similar. 

To conclude, we have 
\be\label{eqn:Ikunknown2}
|I_k| \leq |Ia| + |Ib| \lesssim \eta \tau \max\{p\epsilon\nu, 1\} \log p + O_p(\sqrt{4n^{-1}p\epsilon}). 
\ee

\item 
Next, we analyze $II_k$. Removing the zero part, we can find 
\[
II_k = -(\hat{\mu}_0^{(k)})^\top ((\tilde{\Omega}_0)_{11} - (\tilde{\Omega}_1)_{11}) \hat{\mu}_0^{(k)}+ (\mu^{(k)})^\top ((\tilde{\Omega}_0)_{11} - (\tilde{\Omega}_1)_{11}) \mu^{(k)}.
\]  
Let $R = \hat{\mu}_0^{(k)} + \mu^{(k)}$, then $R \sim N(0, \frac{1}{n_0} I_k)$. 
Rewrite $II_k$ as 
\be
II_k = [R^\top ((\tilde{\Omega}_0)_{11} - (\tilde{\Omega}_1)_{11}) R ]+ 2 (\mu^{(k)})^\top ((\tilde{\Omega}_0)_{11} - (\tilde{\Omega}_1)_{11}) R = IIa + 2IIb. 
\ee

We first consider $IIa = R^\top (I - \Omega_{11}) R$. This follows a non-central chi-square distribution. Since $\eta \gg 1/\sqrt{n}$ and $\tilde{\Omega}$ can recover exactly the non-zeros of $\Omega$, 
\ben
\ba{lll}
E[IIa] & = & 0, \\ 
\var(IIa) & = & \frac{2}{n_0^2} Tr(((\tilde{\Omega}_0)_{11} - (\tilde{\Omega}_1)_{11})^2).
\ea
\een
Furthermore, we can prove that $\sup_x|F_{IIa/\sqrt{\var(IIa)}}(x)-\Phi(x)|\overset{\mathcal P}\to 0$, and so 
\[
IIa =  O(n^{-1}\sqrt{Tr(((\tilde{\Omega}_0)_{11} - (\tilde{\Omega}_1)_{11})^2)}).
\]
If we introduce in the terms, then 
\[
\frac{2}{n_0^2}Tr(((\tilde{\Omega}_0)_{11} - (\tilde{\Omega}_1)_{11})^2) \leq Cn^{-2}(k\xi^2*1\{\xi \gg 1/\sqrt{n}\} + \eta^2\max\{k^2\nu, 1\})
\]
for some constant $C > 0$. 
And so 
\[
IIa = O\biggl(n^{-1}\sqrt{p \epsilon \xi^2*1\{\xi \gg 1/\sqrt{n}\} + \eta^2 \max\{p^2\epsilon^2 \nu, 1\}})\biggr).
\]

Then, we consider $IIb =  (\mu^{(k)})^\top ((\tilde{\Omega}_0)_{11} - (\tilde{\Omega}_1)_{11}) R$. Since $Z\sim N(0, \frac{1}{n_0} I)$, it is clear that $IIb \sim N(0, \frac{1}{n_0} (\mu^{(k)})^\top ((\tilde{\Omega}_0)_{11} - (\tilde{\Omega}_1)_{11})^2 \mu^{(k)})$. 
According to the definition of $\mu^{(k)}$, $\frac{1}{n_0} (\mu^{(k)})^\top ((\tilde{\Omega}_0)_{11} - (\tilde{\Omega}_1)_{11})^2 \mu^{(k)} = \frac{\tau^2}{n_0}\|(\tilde{\Omega}_0)_{11} - (\tilde{\Omega}_1)_{11}\|_{\infty}^2$. 
Therefore, with probability $1 - o(1/p)$, 
\[
|IIb| \leq n^{-1/2}\tau \|(\tilde{\Omega}_0)_{11} - (\tilde{\Omega}_1)_{11}\|_{\infty}.
\]
As a result, with $k = p\epsilon (1+o(1))$, 
\[
|IIb| \leq n^{-1/2}\tau(\xi*1\{\xi \gg 1/\sqrt{n}\} + 2\eta \max\{2p\epsilon\nu, 1\}). 
\]
Combining the results for $IIa$ and $IIb$, we have
\be\label{eqn:IIkunknown2}
|II_k| \leq |IIa| + |IIb| \lesssim O\biggl(n^{-1}\sqrt{p \epsilon \xi^2*1\{\xi \gg 1/\sqrt{n}\} + \eta^2 \max\{p^2\epsilon^2 \nu, 1\}})\biggr).
\ee
\end{itemize}

Combining the results for $I_k$ and $II_k$ in (\ref{eqn:Ikunknown2}) and (\ref{eqn:IIkunknown2}), 
\ben
\ba{lll}
\Delta Q & \leq & \eta \tau \max\{p\epsilon\nu, 1\} \log p + O_p(\sqrt{4n^{-1}p\epsilon}) + O\biggl(n^{-1}\sqrt{p \epsilon \xi^2*1\{\xi \gg 1/\sqrt{n}\} + \eta^2 \max\{p^2\epsilon^2 \nu, 1\}})\biggr) \\
& = &  \eta \tau \max\{p\epsilon \nu, 1\} \log p + O_p(\sqrt{4n^{-1}p\epsilon}). 
\ea
\een
The result is proved. \qed

\section{Proof of Theorem \ref{thm:qdaunknownmu} and \ref{thm:qdaunknown}}\label{app:proof}
We show the proof of Theorems \ref{thm:qdaunknownmu} and \ref{thm:qdaunknown}, followed by the two corollaries where we consider PCS and CLIME as the precision matrix estimators. 

In this section, we use $\Omega = \Omega_1$ for short when there is no confusion.  
\subsection{Proof of Theorem \ref{thm:qdaunknownmu}}\label{sec:proof2}

In this section, we focus on the algorithm for QDA with feature selection, when $\mu$ is unknown. To estimate $\mu$, we use $\mu^* = -\hat{\mu}_0$ for the quadratic part and $d$ for the linear part:
\be\label{eqn:d2}
d \ = \ \Omega \hat{\mu}_1 - \hat{\mu}_0 \ \sim \ N\bigl((I + \Omega)\mu, \ \frac{1}{n_0}I+\frac{1}{n_1}\Omega\bigr).
\ee 
For a threshold $t$, we let $\hat{d}_{j} = d_j*I(|d_j| \geq t)$. When $\max_{1 \leq j \leq p} |d_i| \leq 2{\log p}/\sqrt{n}$, we take $t = 0$ which means the original QDA; otherwise we take $t = 2\sqrt{\log p}/\sqrt{n}$, which means the QDAfs algorithm.

In Section \ref{app:fs}, it is shown that $t = 0$ happens with probability $1 - o(1)$ when $\tau \ll 1/\sqrt{n}$, i.e. in the weak signal region; and $t \neq 0$ happens with probability $1 - o(1)$ when $\tau \gg 1/\sqrt{n}$, i.e. in the strong signal region. For the latter case, the signals can be exactly recovered with probability $1 - o(1)$. Hence, we have original QDA for the weak signal case and QDAfs for the strong signal case. We will discuss them separately.

%
%

\subsubsection{The weak signal region}\label{app:weaksig}
Consider the event $\{t = 0\}$. It happens with probability $1 - o(1)$, so we focus on this event only. It means we apply the original QDA method with estimated $\hat{\mu}$ and $d$. 

For original QDA, the estimated label is $I\{Q(X, \hat{\mu}, \Omega) > 0\}$ where 
\bea
Q(X, \hat{\mu}, \Omega) & = & X^\top(I -{\Omega}) X + 2d^\top X + 
 \hat{\mu}_0^\top (I - {\Omega}) \hat{\mu}_0 + \log |{\Omega}|\\
 & = & Q(X, \mu, \Omega) + \Delta Q \nonumber = S - T_S + \Delta Q.
\eea
Here, $S$ and $T_S$ are given (\ref{proofS}) and  $\Delta Q = Q(X, \hat{\mu}, \Omega) - Q(X, \mu, \Omega)$ gives the difference between the two criteria. 

According to the analysis in Section~\ref{app:ideal}, $S|Y=i$ asymptotically converges to a normal distribution. 
So the error $p_{i, \mu, \Omega} = P((-1)^i*Q > 0)$ can be presented by $\Phi(\cdot)$ with $o(1)$ error, 
\bea\label{eqn:pthma1}
p_{i, \mu, \Omega} & = & P((-1)^i*(S - T_S + \Delta Q) > 0)\nonumber\\
& = & \Phi(\frac{(-1)^i*(T_S - E[S|Y=i]) }{\sqrt{\var(S|Y=i)}} + \frac{(-1)^i*\Delta Q}{\sqrt{\var(S|Y=i)}}) + o(1), \quad i = 0, 1.
\eea
We only need to check the value of this normal probability. 
By Section \ref{app:ideal}, the first term comes to $ -\sqrt{\|\Omega - I\|^2_F/8 + \|\mu\|^2}(1 + o(1)) + o(1)$. So we only need to check $\Delta Q/\sqrt{\var(S|Y=i)}$. 
Results about $\Delta Q$ are in the following lemma.
\begin{lemma}\label{lemma:DeltaQ}
Under the model assumptions, Under the model assumptions,  with probability $1+o(1)$, there is
\be
\left| \frac{\Delta Q}{\sqrt{\var(S|Y=k)}} \right| 
 \leq  \frac{\|\Omega - I\|_F\log p/n + \sqrt{p/n \log\log p}}{\sqrt{\|\Omega - I\|_F^2 + 8\|\mu\|^2}}, \quad k = 0, 1.
\ee
\end{lemma}

Combining (\ref{eqn:pthma1}) with Lemma \ref{lemma:DeltaQ}, the misclassification probability follows that 
\bea\label{eqn:pthm3phi}
p_{i, \mu, \Omega}& \leq & \Phi(-\sqrt{ \|\Omega - I\|_F^2/8+\|\mu\|^2} + c(\frac{\|\Omega - I\|_F\log p/n + \sqrt{p/n \log\log p}}{\sqrt{\|\Omega - I\|_F^2 + 8\|\mu\|^2}} + 1))
\nonumber\\
& \leq & \Phi(-\sqrt{(p\xi^2 + \eta^2 p^2 \nu)/8 + \tau^2 p \epsilon} + c(\frac{\sqrt{p\xi^2 + \eta^2 p^2 \nu}\log p/n + \sqrt{p/n \log\log p}}{\sqrt{(p\xi^2 + \eta^2 p^2 \nu)/8 + \tau^2 p \epsilon} } + 1)).
\eea
In the region of possibility, the main term $\sqrt{ \|\Omega - I\|_F^2/8+\|\mu\|^2}$ goes to infinity. If the second term is at a smaller rate, then the whole term still goes to $-\infty$, and the normal probability is 0. The requirement that the second term is at a smaller rate is equivalent as the successful region in Theorem \ref{thm:qdaunknownmu}:  $\frac{1-\delta}{2} < \max \{2 - 2\alpha - \beta, 1 - 2\theta - \zeta\} = \kappa$. 

Therefore, $MR(QDA) \goto 0$ in this region. 

\vspace{3em}
We then prove that QDA cannot work in the complement region. It means $p_{0, \mu, \Omega} > c$ for a positive constant $c$, and then the classification error cannot converge to 0. According to the definition of $p_{0, \mu, \Omega}$, it is to prove $P(Q(X,\hat{\mu}, \Omega) > 0|Y = 0) > c$. 

Given $X$, we have the following lemma about $Q(X,\hat{\mu}, \Omega)$. 
\begin{lemma}\label{lemma:weakfail}
With probability at least $\Phi(-C)/4$, there is 
\be\label{eqn:qimposs}
Q(X,\hat{\mu}, \Omega) \geq S + 2C\sqrt{X^\top \Omega X /n} + \mu^{\top}(I - \Omega)\mu+ \log|\Omega|+Tr(I - \Omega)/n_0, 
\ee
where $S$ is defined in (\ref{proofS}).
\end{lemma}
We further find the lower bound of RHS in the lemma. 
Consider $S$, $P(S \geq E[S|Y=0]|Y=0) = 1/2 + o(1)$ according to Lemma \ref{lemma:S}. 
The second term is $\sqrt{X^\top \Omega X/n}$. 
By Markov inequality, 
$P(X^\top \Omega X/n \geq E[X^\top \Omega X/n] - 2 \sqrt{\var(X^\top \Omega X/n)}) \geq 3/4$. 
The probability that both inequality holds is no smaller than $1/4$. 

According to Lemma \ref{lemma:quad}, $E[S|Y=0]$, $E[X^\top \Omega X/n]$, and $\sqrt{\var(X^\top \Omega X/n)}$ can be derived. Introducing them into (\ref{eqn:qimposs}) and combining with Lemma \ref{lemmaeb}, with probability at least $\Phi(-C)/16$,  
\begin{eqnarray}\label{eqn:qresult}
Q & \geq & -4\mu^\top \Omega \mu+ Tr(I - \Omega)+ 2C\sqrt{\mu^\top \Omega \mu/n + Tr(\Omega)/n - 2\sqrt{p/n^2}}\nonumber\\
&&  + \log|\Omega|+ Tr(I - \Omega)/n_0\nonumber\\
& \gtrsim & 2C\sqrt{p/n} - 4  \tau^2 p \epsilon - \eta^2p^2\nu.
\end{eqnarray} 

Therefore, when both terms $\tau^2 p \epsilon$ and $\eta^2p^2\nu$ are much smaller than $2C\sqrt{p/n} \goto \infty$, then $Q > 0$ with probability  $\Phi(-C)/16$. It is equivalent with the failure region defined in Theorem \ref{thm:qdaunknownmu}:  $\kappa = \max\{2-2\alpha-\beta, 1-2\theta-\zeta\}< (1-\delta)/2$.
Therefore, $MR(QDA) \geq c$ where $c > 0$ is a constant in this region. 

Combining the region that QDA will have MR converging to 0 and the region that QDA will have at least a constant MR, the result about the weak signal region that $\eta > \delta/2$ in Theorem \ref{thm:qdaunknownmu} is proved. \qed


\subsubsection{The strong signal region}\label{sec:thm3strong} 
Suppose the signals are strong. With probability $1 - o(1)$, the threshold $t = \sqrt{2\log p/n}$ and all the signals $\mu_i \neq 0$ are exactly recovered. Hence, we only consider the event that $\{t = \sqrt{2\log p/n}\}$ and all the signals are exactly recovered. Hence we have QDAfs in this region.  



\vspace{2em}
Similar as the derivation in Section \ref{app:weaksig}, we have 
\be\label{eqn:deltaqstrong}
Q(X, \hat{\mu}, \Omega)  =
\sqrt{\var(S|Y=i)} Z + E[S|Y=i] - \mu^\top (\Omega - I)\mu  + \log|\Omega| + \Delta Q, 
\ee
where $Z \sim N(0,1)$ and 
\[
 \Delta Q =  \displaystyle 2(\hat{\mu}_d^{(t)} - (I + \Omega)\mu)^\top X+\bigl[(\hat{\mu}_0^{(t)})^\top(I-\Omega)\hat{\mu}_0^{(t)}-\mu^\top(I-\Omega)\mu + \frac{1}{n_0}Tr(\Omega_1^{(d)} - I)\bigr] .
\]
According to the analysis in Section \ref{app:weaksig}, to show $MR(QDA) \goto 0$, we only need to show, with probability $1 - o(1)$, 
\be\label{eqn:deltaqmag} 
\Delta Q  = o({p\xi^2 + \eta^2 p^2 \nu + 8 \tau^2 p \epsilon}). 
\ee

To prove (\ref{eqn:deltaqmag}), recall that we have the magnitude about $\Delta Q$ in Lemma \ref{lemma:StrongDeltaQ}, 
\[
 |\Delta Q|   \leq   O(\sqrt{p\epsilon(\xi^2 + p\epsilon\eta^2 \nu)}/n) + \sqrt{p\epsilon} \tau \log p (1 + o(1)).
\]
Consider the three terms in the RHS. Since $\xi \ll p^{-1/2}$ and $1/\sqrt{n} \ll \tau $, $c p\epsilon \xi/n \ll \sqrt{p}\epsilon/n \ll \tau^2 \sqrt{p}\epsilon$. The last two terms can be combined so that 
\[
O(\sqrt{p\epsilon(\xi^2 + p\epsilon\eta^2 \nu)}/n) + \sqrt{p\epsilon} \tau \log p (1 + o(1)) 
\leq C\sqrt{\eta^2p^2\nu + \tau^2 p\epsilon}.
\]
So, it is of $o(\eta^2p^2\nu + \tau^2 p\epsilon)$ when $\eta^2p^2\nu + \tau^2 p\epsilon \goto \infty$. 

Therefore, in the region of possibility identified by part (ii) of Theorem \ref{thm:qdaunknownmu}, (\ref{eqn:deltaqmag}) always holds, and hence $MR$ converges to 0. 

Now we consider the region of impossibility. When the mean vector $\mu$ is unknown, the region of impossibility is no smaller than the case that the mean vector $\mu$ is already known (\cite{LeCam}). For the latter, the region of impossibility is depicted in Proposition \ref{thm:idealunequal}. 

When the signals are strong, the region of impossibility is the same with that in Proposition \ref{thm:idealunequal}. Hence, when the mean vector is unknown, any classifier $L$ has $MR(L) \goto 1/2$ in this region.  \qed

\subsection{Proof of Theorem \ref{thm:qdaunknown}}\label{sec:proof4}
Similar as the proof of Theorem \ref{thm:qdaunknownmu}, now we want to find the difference between $Q(X, \hat{\mu}, \hat{\Omega})$ and $Q(X, \mu, \Omega)$. Let $\Delta Q  = Q(X, \hat{\mu}, \hat{\Omega}) - Q(X, \mu, \Omega)$, then 
\[
Q(X, \hat{\mu}, \hat{\Omega}) = Q(X, \mu, \Omega) + \Delta Q. 
\]
Hence, 
\[
p_{i, \mu, \Omega} = \Phi(-T + (-1)^i*\frac{\Delta Q}{\sigma_S}) + o(1), \quad i = 0, 1, 
\]
where $\sigma^2_S = 2(p\xi^2 + \eta^2 p^2 \nu + 8\tau^2 p \epsilon)$, and $T = \frac{1}{2\sqrt{2}}\sqrt{p\xi^2 + \eta^2 p^2 \nu + 8\tau^2 p \epsilon}(1 + o(1))$. Therefore, we need that 
\be\label{omegadeltaq}
|\Delta Q| \ll \frac{1}{2}(p\xi^2 + \eta^2 p^2 \nu + 8\tau^2 p \epsilon). 
\ee
Once it holds, then $MR(QDA) \goto 0$ and the region that QDAfs succeeds can be figured out. 

Recall that in Theorem \ref{thm:qdaunknown}, there is $\Delta_{\hat{\Omega}} (p\eta  + p \tau^2 \epsilon + \sqrt{p}\log p)+ p \Delta_{\hat{\Omega}}^2 +p/n \ll p\xi^2 + \eta^2 p^2 \nu + \tau^2 p \epsilon$. So, to prove (\ref{omegadeltaq}), we only need to show 
\be\label{eqn:omegadelta2}
|\Delta Q| \leq C \Delta_{\hat{\Omega}} (p\eta  + p \tau^2 \epsilon + \sqrt{p}\log p)+ p \Delta_{\hat{\Omega}}^2 + \frac{2p}{n_1}.
\ee

The main term to use in the proof is the spectral norm of $\hat{\Omega} - \Omega$, denoted by $\Delta_{\hat{\Omega}}$. For different algorithms to estimate the precision matrix, the resultant $\Delta_{\hat{\Omega}}$ is different. 
Hence, in this part, we derive the result based on $\Delta_{\hat{\Omega}}$, and explain the results when we choose PCS or CLIME in the next two subsections. 

In the proof of Theorem \ref{thm:qdaunknownmu}, we introduce $Q(X, \hat{\mu}, \Omega)$. Now, when the precision matrix is also unknown, the classification criteria change to  
\begin{eqnarray}\label{eqn:diffprecision}
\Delta Q & = & Q(X, \hat{\mu}, \hat{\Omega}) - Q(X, \hat{\mu}, \Omega)\nonumber\\
& = & X^\top (\Omega - \hat{\Omega}) X + 2(\hat{\mu} \circ \hat{d})^{\top}(\hat{\Omega} - \Omega) X  + (\hat{\mu}_0 \circ \hat{d})^\top (\Omega - \hat{\Omega})(\hat{\mu}_0 \circ \hat{d}) \nonumber\\
&& + \frac{1}{n_1} Tr(\hat{\Omega} - I) + \log |\hat{\Omega}|- \log |\Omega|\nonumber\\
& = & (X - \hat{\mu}_1 \circ \hat{d})^{\top}  (\Omega - \hat{\Omega}) (X - \hat{\mu}_1 \circ \hat{d})\nonumber\\
&& + (\hat{\mu}_0 \circ \hat{d})^\top (\Omega - \hat{\Omega})(\hat{\mu}_0 \circ \hat{d}) - (\hat{\mu}_1 \circ \hat{d})^\top (\Omega - \hat{\Omega})(\hat{\mu}_1 \circ \hat{d}) \nonumber \\
&& + \frac{1}{n_1} Tr(\hat{\Omega} - I) + \log |\hat{\Omega}|- \log |\Omega|\nonumber\\
& = & A + B + C.
\end{eqnarray}
Hence, we have to find proper bounds for the terms $A$, $B$, and $C$ under the sparse case and the dense case. 

Before we discuss this, we have to introduce the following lemma and results as preparation.
\begin{lemma}\label{lemma:trace}
Consider two symmetric matrices $A$ and $B$, where the eigenvalues of $B$ are $\lambda_1$, $\cdots, \lambda_p$; then,
\[
Tr(AB) \leq \|A\| *\sum_{i=1}^p |\lambda_i|.
\]
Furthermore, if $B$ is a positive semi-definite matrix,
\[
Tr(AB) \leq Tr(B)\|A\|.
\]
\end{lemma}
According to Lemma \ref{lemma:trace} and Lemma \ref{lemmaeb}, we have the following conclusion:
\begin{equation}\label{eqn:trdiff}
Tr((\Omega - \hat{\Omega})(\Omega^{-1} - I)) \leq p \|\Omega^{-1} - I\| \Delta_{\hat{\Omega}}, 
\,\, 
Tr((\Omega - \hat{\Omega})\Omega^{-1} ) \lesssim p \Delta_{\hat{\Omega}}.
\ee
Let the eigenvalues of $\Omega - I$ be $m_1, \cdots, m_p$, and those of $\hat{\Omega} - I$ be $\hat{m}_1, \hat{m}_2, \cdots, \hat{m}_p$. According to Weyl's inequality (\cite{weyl1912asymptotische}),
\be\label{diffeigenvalue}
|m_i - \hat{m}_i| \leq \Delta_{\hat{\Omega}}.
\ee
With these preparations, we begin to consider the three terms. 

{\bf Weak signal region $\theta > \delta/2$. } With such conditions, the signals are too weak to discover. Hence, the thresholding vector $\hat{d}$ is not helpful at all. We choose $t = 0$.  
\begin{itemize}
\item Consider the part $A$. In this case, it is 
\[
A = (X - \hat{\mu}_1)^{\top} (\Omega - \hat{\Omega}) (X - \hat{\mu}_1). 
\]
Obviously, It depends on the distribution of $X$. 

When $Y = 0$, then $X \sim N(-\mu, I)$ and hence $X - \hat{\mu}_1 \sim N(-2\mu, I + \frac{1}{n_1}\Omega^{-1})$. 
Consider the randomness of $\hat{\mu}_1$ first. According to Lemma \ref{lemma:quad}, we have
\be
E[A]  =  Tr(\Omega - \hat{\Omega} + \frac{1}{n_1}(\Omega - \hat{\Omega})\Omega^{-1}) + 4\mu^{\top} (\Omega - \hat{\Omega}) \mu.
\ee
The variance of $A$ is given as 
\be
\var(A) = 2Tr(((\Omega - \hat{\Omega})(I + \frac{1}{n}\Omega^{-1}))^2) + 16\mu^{\top}(\Omega - \hat{\Omega})(I + \frac{1}{n}\Omega^{-1})(\Omega - \hat{\Omega})\mu. 
\ee
Further, the Berry-Ess\'een theorem gives $\sup_x|F_{(A - E[A])/\sqrt{\var(A)}}(x)-\Phi(x)|\overset{\mathcal P} \to 0$. 
So, with probability $1 - o(1/p)$, $|A - E[A]| \leq \sqrt{\var(A)} \log p$. 

According to (\ref{eqn:trdiff}), rearrange the expectation, we have 
\be\label{eqn:case1d}
|E[A] - \frac{1+n_1}{n_1}Tr(\Omega - \hat{\Omega})| \leq \frac{1}{n_1}p \Delta_{\hat{\Omega}}\|\Omega^{-1} - I\| + 4\|\mu\|^2 \Delta_{\hat{\Omega}}, 
\ee
and $\var(A) = 2p \Delta_{\hat{\Omega}}^2(1+o(1))$. 
Therefore, combining with the result that $|A - E[A]| \leq \sqrt{\var(A)} \log p$, with probability $1 - o(1/p)$,  
\be\label{unknownvard}
|A - \frac{1+n_1}{n_1}Tr(\Omega - \hat{\Omega})| \leq \frac{1}{n_1}p \Delta_{\hat{\Omega}}\|\Omega^{-1} - I\| + \|\mu\|^2 \Delta_{\hat{\Omega}} + \sqrt{2p}\Delta_{\hat{\Omega}} \log p.
\ee

When $Y = 1$, then $X\sim N(\mu, \Omega^{-1})$.  
Therefore, $X - \hat{\mu}_1 \sim N(0, \frac{1+n_1}{n_1}\hat{\Omega}^{-1})$. Again, according to Lemma \ref{lemma:quad}, we have
\begin{align*}
E[A] & =  \frac{1+n_1}{n_1}Tr((\Omega - \hat{\Omega})\Omega^{-1}) + 4\mu^{\top} (\Omega - \hat{\Omega}) \mu \\
& =   \frac{1+n_1}{n_1}Tr(\Omega - \hat{\Omega}) + \frac{1+n_1}{n_1}Tr((\Omega - \hat{\Omega})(\Omega^{-1} - I)) + 4\mu^{\top} (\Omega - \hat{\Omega}) \mu.
\end{align*}
According to (\ref{eqn:trdiff}), rearrange it, and we have
\be\label{eqn:case2A}
|E[A] - \frac{1+n_1}{n_1}Tr(\Omega - \hat{\Omega})| \leq p \Delta_{\hat{\Omega}}\|\Omega^{-1} - I\| + 4\|\mu\|^2 \Delta_{\hat{\Omega}}. 
\ee
The variance of $A$ is given as 
\be
2(\frac{1+n_1}{n_1})^2Tr(((\Omega - \hat{\Omega})\Omega^{-1})^2) + 16\frac{1+n_1}{n_1}\mu^{\top}(\Omega - \hat{\Omega})\Omega^{-1}(\Omega - \hat{\Omega})\mu \lesssim 2p \Delta_{\hat{\Omega}}^2. 
\ee
Therefore, we have that 
\be\label{unknownvarA}
|A - \frac{1+n_1}{n_1}Tr(\Omega - \hat{\Omega})| \leq p \Delta_{\hat{\Omega}}\|\Omega^{-1} - I\| + \|\mu\|^2 \Delta_{\hat{\Omega}} + \sqrt{2p}\Delta_{\hat{\Omega}} \log p.
\ee

\item Consider $B$. Let $B_1 = \hat{\mu}_0^\top (\Omega - \hat{\Omega})\hat{\mu}_0$ and $B_2 = \hat{\mu}_1^\top (\Omega - \hat{\Omega})\hat{\mu}_1$, then $B = B_1 - B_2$. Note that 
\[
\hat{\mu}_0 \sim N(-\mu, \frac{1}{n_0} I), \qquad \hat{\mu}_1 \sim N(\mu, \frac{1}{n_1}\Omega^{-1}).
\]
With the same analysis as in $A$, we have
\[
E[B_1] = Tr(\Omega - \hat{\Omega})/n_0 + \mu^{\top}(\Omega - \hat{\Omega})\mu, 
\quad
E[B_2] = Tr((\Omega - \hat{\Omega})\Omega^{-1})/n_1 + \mu^{\top}(\Omega - \hat{\Omega})\mu.
\]
The variance of them is 
\[
\var[B_1] = 2Tr((\Omega - \hat{\Omega})^2)/n_0^2 + 4\mu^{\top}(\Omega - \hat{\Omega})^2\mu/n_0 \leq 2p\Delta_{\hat{\Omega}}^2/n^2 + 4 \|\mu\|^2 \Delta_{\hat{\Omega}}^2/n.
\]
Note that $\tau^2 \ll 1/n$, so $\|\mu\|^2 \ll p/n$. Therefore, $\var(B_1)  \lesssim 2p\Delta_{\hat{\Omega}}^2/n^2$. Similarly, we have that $\var(B_2)  \lesssim 2p\Delta_{\hat{\Omega}}^2/n^2$. 

In all, for part $B$, we have, with probability $1 - o(1/p)$, 
\bea\label{unknownvarB}
|B| & = & |B_1 - B_2| \nonumber\\
&\lesssim & |E[B_1] - E[B_2]| + 2\sqrt{2p \Delta_{\hat{\Omega}}^2/n^2}\log p\nonumber\\
& \lesssim &  |Tr(\Omega - \hat{\Omega})/n_0 - Tr((\Omega - \hat{\Omega})\Omega^{-1})/n_1| + 2\sqrt{2p \Delta_{\hat{\Omega}}^2/n^2}\log p \nonumber\\
& \lesssim & p \Delta_{\hat{\Omega}}/n + 2\sqrt{2p \Delta_{\hat{\Omega}}^2/n^2}\log p . 
\eea

\item Finally, we consider $C$. Note that  
\begin{align}
C&  =  \frac{1+n_1}{n_1}Tr(\hat{\Omega} - {\Omega}) +  \frac{1}{n_1} Tr(\Omega - I) + Tr(\Omega - \hat{\Omega}) + \log |\hat\Omega|- \log|\Omega| \nonumber \\
& =  \frac{1+n_1}{n_1}Tr(\hat{\Omega} - {\Omega}) + I + II,
\end{align}
where 
\be\label{unknownvar1}
I =  \frac{1}{n_1} Tr(\Omega - I)  = \frac{p \xi}{n_1}
\ee
and $II = Tr(\Omega - \hat{\Omega}) + \log |\hat\Omega|- \log|\Omega| $.

For Part $II$, note that, when $ \Delta_{\hat{\Omega}} = o(1)$, 
\bea
\log |\hat\Omega|- \log|\Omega|  & = & \sum_{i=1}^p \log \frac{1+\hat{m}_i}{1+{m}_i}\nonumber\\
& = & \sum_{i = 1}^p \log (1 + \frac{\hat{m}_i - m_i}{1+m_i})\nonumber\\
& = & \sum_{i = 1}^p  \left[\frac{\hat{m}_i - m_i}{1+m_i} - ( \frac{\hat{m}_i - m_i}{1+m_i})^2/2 + ( \frac{\hat{m}_i - m_i}{1+m_i})^3/3 + \cdots \right]\nonumber\\
& = & \sum_{i = 1}^p  \frac{\hat{m}_i - m_i}{1+m_i}  - p \Delta_{\hat{\Omega}}^2/2(1+o(1))\nonumber\\
& = & \sum_{i=1}^p (\hat{m}_i - m_i)\sum_{l = 1}^{\infty} (-m_i)^l  - p \Delta_{\hat{\Omega}}^2/2(1+o(1))\nonumber\\
& = & \sum_{i=1}^p (\hat{m}_i - m_i) + \sum_{i = 1}^p  \frac{(\hat{m}_i - m_i)(-m_i)}{1+m_i} - p \Delta_{\hat{\Omega}}^2/2(1+o(1)).
\eea 
Noting $\sum_{i=1}^p (\hat{m}_i - m_i) = - Tr(\Omega - \hat{\Omega})$, therefore, 
\be\label{unknownvar2}
II = \sum_{i = 1}^p  \frac{(\hat{m}_i - m_i)(-m_i)}{1+m_i} - p \Delta_{\hat{\Omega}}^2/2(1+o(1)).
\ee
Combining (\ref{unknownvar1}) and (\ref{unknownvar2}), we have
\begin{align}\label{unknownvarC}
C & =  \frac{1+n_1}{n_1}Tr(\hat{\Omega} - {\Omega}) +  \frac{p \xi}{n_1} +  \sum_{i = 1}^p  \frac{(\hat{m}_i - m_i)(-m_i)}{1+m_i} - p \Delta_{\hat{\Omega}}^2/2(1+o(1)) \nonumber\\
& \leq  \frac{1+n_1}{n_1}Tr(\hat{\Omega} - {\Omega}) +  \frac{p \xi}{n_1} +  p\Delta_{\hat{\Omega}} \|\Omega^{-1} - I\| - p \Delta_{\hat{\Omega}}^2/2(1+o(1)), 
\end{align}
where the last part comes from that $\|\Omega^{-1} - I\| = \max_{i} |\frac{-m_i}{1+m_i}|$. 
\end{itemize}
Combining (\ref{unknownvard}), (\ref{unknownvarB}), (\ref{unknownvarC}), and (\ref{diffeigenvalue}), when $Y = 0$, we have
\begin{align}\label{case1}
|\Delta| & \leq (\frac{1}{n_1}p \Delta_{\hat{\Omega}}\|\Omega^{-1} - I\| + \|\mu\|^2 \Delta_{\hat{\Omega}} + \sqrt{2p}\Delta_{\hat{\Omega}} \log p) + (p \Delta_{\hat{\Omega}}/n + 2\sqrt{2p\Delta_{\hat{\Omega}}^2/n^2\log p}) \nonumber\\
& + (\frac{p \xi}{n_1} + p\Delta_{\hat{\Omega}} \|\Omega^{-1} - I\| + p \Delta_{\hat{\Omega}}^2)\nonumber\\
& \lesssim  \frac{p \xi}{n_1} + \Delta_{\hat{\Omega}} (p\|\Omega^{-1} - I\| + \|\mu\|^2 + \sqrt{2p}\log p)+ p \Delta_{\hat{\Omega}}^2.
\end{align}
Combining (\ref{unknownvarA}), (\ref{unknownvarB}), (\ref{unknownvarC}), and (\ref{diffeigenvalue}), and with $\|\Omega^{-1} - I\| = \max_{i} |\frac{-m_i}{1+m_i}|$, when $Y = 1$, we have 
\begin{align}\label{case2}
|\Delta| & \leq (p \Delta_{\hat{\Omega}}\|\Omega^{-1} - I\| + \|\mu\|^2 \Delta_{\hat{\Omega}} + \sqrt{2p}\Delta_{\hat{\Omega}} \log p) + (p \Delta_{\hat{\Omega}}/n + \sqrt{2p\Delta_{\hat{\Omega}}^2/n^2\log p}) \nonumber\\
& + (\frac{p \xi}{n_1} + p\Delta_{\hat{\Omega}} \|\Omega^{-1} - I\| + p \Delta_{\hat{\Omega}}^2)\nonumber\\
& \leq  \frac{p \xi}{n_1} + \Delta_{\hat{\Omega}} (p\|\Omega^{-1} - I\| + \|\mu\|^2 + \sqrt{2p}\log p)+ p \Delta_{\hat{\Omega}}^2.
\end{align}

Combine (\ref{case1}) and (\ref{case2}), and we can see that, when $\Delta_{\hat{\Omega}} = o(1)$, with probability $1 - o(1)$, 
\be\label{eqn:unknowndelta}
|\Delta Q| \leq  \frac{p \xi}{n_1} + \Delta_{\hat{\Omega}} (p\|\Omega^{-1} - I\| + \|\mu\|^2 + \sqrt{2p}\log p)+ p \Delta_{\hat{\Omega}}^2.
\ee

\vspace{2em}
{\bf Strong signal region $\theta < \delta/2$. }
We consider the case in which the signals are strong. With such a condition, all the signals can be recovered exactly (see the analysis in the proof for Theorem \ref{thm:qdaunknownmu}). The threshold $\hat{d}$ guarantees that the true signals are non-zero, and the noises are all zero. 

For simplicity, we rearrange the features so that all the features with $\mu(j) \neq 0$ rank first, and those with $\mu(j) = 0$ rank last. Hence, the top $k$ features have non-zero means. We use $\mu^{(k)}$ to denote the sub-vector of $\mu$ that contains the first $k$ elements only, and do the same for $\hat{\mu}_0$, $\hat{\mu}_1$. We also decompose the precision matrix $\Omega$, as follows:
\[
\Omega = \left(\begin{array}{ll}
\Omega_{11} & \Omega_{12}\\
\Omega_{12}^{\top} & \Omega_{22}
\end{array}
\right), 
\]
where $\Omega_{11} \in \mathcal{R}^{k \times k}$. We perform the same decomposition on $\hat{\Omega}$ and $\Omega^{-1}$, correspondingly. 

For the sub-matrices, we have that $\|\Omega_{11} - \hat{\Omega}_{11}\| \leq \Delta_{\hat\Omega}$. 
According to Lemma \ref{lemma:trace} and Lemma \ref{lemmaeb}, we have the following conclusion:
\begin{equation}\label{eqn:trdiffstrong}
Tr((\Omega_{11} - \hat{\Omega}_{11})(\Omega^{-1})_{11} ) \lesssim k \Delta_{\hat{\Omega}}.
\ee

With the new notations, the three terms of interest now become
\begin{eqnarray*}
A &= &(X - \hat{\mu}_1^{(k)})^{\top}  (\Omega - \hat{\Omega}) (X - \hat{\mu}_1^{(k)}), \\
B &= &(\hat{\mu}_0^{(k)})^\top (\Omega_{11} - \hat{\Omega}_{11})(\hat{\mu}_0^{(k)}) - (\hat{\mu}_1^{(k)})^\top (\Omega_{11} - \hat{\Omega}_{11})(\hat{\mu}_1^{(k)}),  \\
C &= &\frac{1}{n_1} Tr(\hat{\Omega} - I) + \log |\hat{\Omega}|- \log |\Omega|.
\end{eqnarray*}
Again, the analysis of $A$ depends on the distribution of the new input $X$, yet $B$ and $C$ do not rely on $X$. 
Now we analyze them one by one. 
\begin{itemize}
\item 
For Part $A$, we first consider the case note that $X \sim N(-\mu, I)$.
According to Lemma \ref{lemma:quad}, we have
\be
E[A]   =   Tr(\Omega - \hat{\Omega}) + \frac{1}{n_1}Tr((\Omega_{11} - \hat{\Omega}_{11})(\Omega^{-1})_{11}) + 4\mu^{\top} (\Omega - \hat{\Omega}) \mu.
\ee
Rearrange it, and we have
\be\label{eqn:sparseA}
|E[A] - Tr(\Omega - \hat{\Omega})| \leq \frac{2k}{n_1} \Delta_{\hat{\Omega}} + 4\|\mu\|^2 \Delta_{\hat{\Omega}}. 
\ee
The variance of $A$ is given as 
\be
2Tr(((\Omega - \hat{\Omega})\tilde{\Sigma})^2) + 16\mu^{\top}(\Omega - \hat{\Omega})\tilde{\Sigma}(\Omega - \hat{\Omega})\mu \approx 2p \Delta_{\hat{\Omega}}^2. 
\ee
Then, with probability $1 - o(1/p)$,  
\be\label{unknownsparseA}
|A - Tr(\Omega - \hat{\Omega})| \leq 2k \Delta_{\hat{\Omega}}/n_1 +4 \|\mu\|^2 \Delta_{\hat{\Omega}} + \sqrt{2p}\Delta_{\hat{\Omega}} \log p.
\ee
For the case $X \sim N(\mu, \Omega^{-1})$, with similar derivations, we have that 
\be\label{unknownsparseA2}
|A - Tr(\Omega - \hat{\Omega})| \leq p \Delta_{\hat\Omega} \|\Omega^{-1} - I\| + 2k \Delta_{\hat{\Omega}}/n_1 +4 \|\mu\|^2 \Delta_{\hat{\Omega}} + \sqrt{2p}\Delta_{\hat{\Omega}} \log p.
\ee

\item Now we consider Part $B$. 

With Lemma \ref{lemma:quad}, we have
\begin{align*}
E[(\hat{\mu}_0^{(k)})^\top (\Omega_{11} - \hat{\Omega}_{11})(\hat{\mu}_0^{(k)})] & = \frac{1}{n_0} Tr(\Omega_{11} - \hat{\Omega}_{11}) + (\mu^{(k)})^\top (\Omega_{11} - \hat{\Omega}_{11})({\mu}^{(k)}), \\
E[(\hat{\mu}_1^{(k)})^\top (\Omega_{11} - \hat{\Omega}_{11})(\hat{\mu}_1^{(k)})] & = \frac{1}{n_1} Tr((\Omega_{11} - \hat{\Omega}_{11})(\Omega^{-1})_{11}) + (\mu^{(k)})^\top (\Omega_{11} - \hat{\Omega}_{11})({\mu}^{(k)}).
\end{align*}
Therefore, the difference between the two terms is that
\be\label{eqn:unknownsparse1}
E[B] = \frac{1}{n_0} Tr(\Omega_{11} - \hat{\Omega}_{11}) - \frac{1}{n_1} Tr((\Omega_{11} - \hat{\Omega}_{11})(\Omega^{-1})_{11}) \leq (\frac{1}{n_0} + \frac{1}{n_1}) k \Delta_{\hat{\Omega}}. 
\ee
The variance is
\begin{eqnarray*}
\var[(\hat{\mu}_0^{(k)})^\top (\Omega_{11} - \hat{\Omega}_{11})(\hat{\mu}_0^{(k)})] & = & \frac{2}{n_0^2}Tr((\Omega_{11} - \hat{\Omega}_{11})^2) +  \frac{4}{n_0}(\mu^{(k)})^\top (\Omega_{11} - \hat{\Omega}_{11})^2({\mu}^{(k)})\\
& \lesssim & \frac{2}{n^2} k \Delta_{\hat{\Omega}}^2 + \frac{4}{n} \|\mu\|^2 \Delta_{\hat{\Omega}}^2\\
& \leq & \frac{6}{n} \|\mu\|^2 \Delta_{\hat{\Omega}}^2,
\end{eqnarray*}
where the last inequality comes from $\|\mu\|^2 = k\tau^2$ and $\tau^2 > 1/n$ in the sparse region. 

The variance for $(\hat{\mu}_1^{(k)})^\top (\Omega_{11} - \hat{\Omega}_{11})(\hat{\mu}_1^{(k)})$ can be calculated in the same way, and the result is the same. Therefore, we have
 \be\label{unknownsparseB}
|B| \lesssim (\frac{1}{n_0} + \frac{1}{n_1}) k \Delta_{\hat{\Omega}} + 8\sqrt{\frac{\log p}{n} } \|\mu\| \Delta_{\hat{\Omega}}.
\ee

\item For Part $C$, note that it does not change in the strong signal case. Therefore, (\ref{unknownvarC}) still works. Here we revise it as  
 \begin{align}\label{unknownsparseC}
C & = & Tr(\hat{\Omega} - {\Omega}) + \frac{1}{n_1}Tr(\hat{\Omega} - I)  +  \sum_{i = 1}^p  \frac{(\hat{m}_i - m_i)(-m_i)}{1+m_i} - p \Delta_{\hat{\Omega}}^2/2(1+o(1))\nonumber \\
& \leq & Tr(\hat{\Omega} - {\Omega}) + \frac{1}{n_1} Tr(\hat{\Omega} - I)  +  p \Delta_{\hat\Omega} \|\Omega^{-1} - I\| - p \Delta_{\hat{\Omega}}^2/2(1+o(1)).
\end{align}
\end{itemize}

Combining (\ref{unknownsparseA}), (\ref{unknownsparseA2}), (\ref{unknownsparseB}), and (\ref{unknownsparseC}), we have
\begin{align}\label{sparse}
|\Delta Q| & \leq (p \Delta_{\hat{\Omega}} \|\Omega^{-1} - I\| + 2k \Delta_{\hat{\Omega}}/n +4 \|\mu\|^2 \Delta_{\hat{\Omega}} + \sqrt{2p}\Delta_{\hat{\Omega}} \log p)\nonumber\\
& + (\frac{k}{n} \Delta_{\hat{\Omega}} + 8\sqrt{\frac{\log p}{n} } \|\mu\| \Delta_{\hat{\Omega}}) 
 + (\frac{2p}{n_1}  +  p \Delta_{\hat{\Omega}} \|\Omega^{-1} - I\| + p \Delta_{\hat{\Omega}}^2)\nonumber\\
& \leq  \Delta_{\hat{\Omega}} (7\|\mu\|^2 + \sqrt{2p}\log p + 2p\|\Omega^{-1} - I\| )+ \frac{2p}{n_1} + p \Delta_{\hat{\Omega}}^2. 
\end{align}

As a conclusion, combining (\ref{eqn:unknowndelta}) and (\ref{sparse}), we have that 
\be
|\Delta Q| \leq C \Delta_{\hat{\Omega}} (p\|\Omega^{-1} - I\|  + \|\mu\|^2 + \sqrt{p}\log p)+ p \Delta_{\hat{\Omega}}^2 + \frac{2p}{n_1}.
\ee
Further, according to Lemma \ref{lemma2}, that $\|\mu\|^2 = \tau^2 p\epsilon(1 + o(1))$, and $\|\Omega^{-1} - I\| \leq 2(\eta + \xi)$ when $1 < \beta < 2$, we have that 
\be
|\Delta Q| \leq C \Delta_{\hat{\Omega}} (p\eta  + p \tau^2 \epsilon + \sqrt{p}\log p)+ p \Delta_{\hat{\Omega}}^2 + \frac{2p}{n_1} .
\ee
Therefore, (\ref{eqn:omegadelta2}) is proved.  
Combining it with (\ref{omegadeltaq}), we have that when $C \Delta_{\hat{\Omega}} (p\eta  + p \tau^2 \epsilon + \sqrt{p}\log p)+ p \Delta_{\hat{\Omega}}^2 + \frac{2p}{n} \ll p\xi^2 + \eta^2 p^2 \nu + 8\tau^2 p \epsilon$, $MR(QDA) \goto 0$. Theorem \ref{thm:qdaunknown} is proved. 

\subsection{Proof of two corollaries}
Now we consider the methods when we use PCS as the precision matrix estimator and when we use CLIME as the precision matrix estimator. The results are presented in Corollaries \ref{cor:PCS} and \ref{cor:CLIME}. The proofs are as follows. We only need to check whether the condition in Theorem \ref{thm:qdaunknown} is satisfied or not. 

\subsubsection{PCS}
Consider the PCS algorithm in \cite{PCS}. Note that there is $1 < \beta < 2$ and $\alpha < \delta/2$. 
Under such conditions, PCS can recover the exact support with probability $1 - o(1/p)$, and 
\[
\|\Omega - \hat{\Omega}\|_{max} \leq \sqrt{\log p/n}. 
\]
Since $1 < \beta < 2$, with probability $1 - o(1/p)$, each row has at most two non-zero elements. Hence, we have
\[
\|\Omega - \hat{\Omega}\|_1 \leq 2\sqrt{\log p/n}, \quad \|\Omega - \hat{\Omega}\|_{\infty} \leq 2\sqrt{\log p/n}.
\]
Therefore, we have the spectral norm, 
\be
\Delta_{\hat{\Omega}} = \|\Omega - \hat{\Omega}\| \leq 2\sqrt{\log p/n}.
\ee

Therefore, introduce the term $\Delta_{\hat{\Omega}}$ into Theorem \ref{thm:qdaunknown}; we then have
\bea
&&C \Delta_{\hat{\Omega}} (p\eta  + p \tau^2 \epsilon + \sqrt{p}\log p)+ p \Delta_{\hat{\Omega}}^2 + \frac{2p}{n}\nonumber\\
& \lesssim & \frac{2p}{n} + C\frac{2\sqrt{\log p}}{\sqrt{n}} (p\eta + p\tau^2\epsilon + \sqrt{p} \log p) + \frac{4p}{n}\log p\nonumber\\
& \lesssim &  C\frac{2\sqrt{\log p}}{\sqrt{n}} (p\eta + p\tau^2\epsilon + \sqrt{p} \log p) + \frac{6p}{n}\log p. 
\eea
Recall that $\eta \gg \sqrt{\log p}/\sqrt{n}$ and $\tau \ll 1/\sqrt{n}$ under the current conditions. Hence, we further reduce the terms to $\frac{C\log p}{\sqrt{n}} p \eta (1 + o(1))$. 

According to Theorem \ref{thm:qdaunknown}, to make sure $MR(QDA) \goto 0$,  we need that, 
\[
\frac{C\log p}{\sqrt{n}} p \eta (1 + o(1)) \ll p\xi^2 + \eta^2 p^2 \nu + 8\tau^2 p \epsilon. 
\]
It is equivalently with one of the following conditions, that 
\[
\beta < 1  - \alpha + \delta/2,
\]
or
\[
\zeta < \alpha + \delta/2 - 2\theta.
\]
Hence, Corollary \ref{cor:PCS} is proved.

\subsubsection{CLIME}

In this section, we consider the Constrained $\ell_1$-Minimization for Inverse Matrix Estimation (CLIME) method in \cite{CLIME}. 

\cite{CLIME} has proved that, when $\lambda_1(\Omega)/\lambda_n(\Omega)$ is finite, and $s$, the maximum number of off-diagonal entries in each row, is bounded by the $o(\sqrt{n})$ term, the CLIME algorithm can achieve the following rate:
\[
\sup_{\mathcal{G}_0(s, M_{n,p})} E[\Delta_{\hat{\Omega}}^2] \asymp (1 - \xi + \eta s)^2 s^2 \frac{\log p}{n}.
\]

To make sure the results hold, we need $\beta > 1 - \delta/2$ and $\alpha < \delta/2$. 

\begin{itemize}
\item When $1 - \delta/2 < \beta < 1$, according to Chebyshev's inequality, we have
\be
P(\Delta_{\hat{\Omega}} > p^{(1 - \beta - \delta/2)+\epsilon'} \log p) \leq  \frac{(1 - \xi + \eta s)^2 s^2 \frac{\log p}{n}}{p^{2 - 2\beta - \delta + 2\epsilon'} \log^2 p} \leq  \frac{ 1}{p^{2\epsilon'} \log p} \rightarrow 0,
\ee
where the last inequality comes from $n = p^{\delta}$, $s \leq 2p^{1 - \beta}$, and $0 < \epsilon' < \beta - (1 - \delta/2)$ can be any constant. 

We introduce the result into Theorem \ref{thm:qdaunknown} and find that there is no successful region. 

\item When $1< \beta < 2$, according to Chebyshev's inequality, we have
 \be
P(\Delta_{\hat{\Omega}} > \log p/\sqrt{n}) \leq  \frac{(1 - \xi + \eta s)^2 s^2 \frac{\log p}{n}}{ \log^2 p/n} \leq  \frac{ 1}{ \log p} \rightarrow 0.
\ee
It is the same with the result in PCS, and hence the successful region is also the same. 
\end{itemize}

\subsection{Proof of Lemmas}
In this section, we prove the lemmas used in the proof of Theorems \ref{thm:qdaunknownmu} and \ref{thm:qdaunknown}. 

For the convenience of readers, we copy the lemmas to prove in each section. 
\subsubsection{Proof of Lemma \ref{lemma:DeltaQ}}
\begin{lemma*}
Under the model assumptions,  with probability $1+o(1)$, there is
\ben
\left| \frac{\Delta Q}{\sqrt{\var(S|Y=k)}} \right| 
 \leq  \frac{\|\Omega - I\|_F\log p/n + \sqrt{p/n \log\log p}}{\sqrt{\|\Omega - I\|_F^2 + 8\|\mu\|^2}}, \quad k = 0, 1.
\een
\end{lemma*}
Direct calculations show that  
\bea
\Delta Q & = &  2(d^\top -\mu^\top(I+\Omega))X+\left[\hat{\mu}_0^\top (I - {\Omega}) \hat{\mu}_0 -\mu^\top(I-\Omega)\mu\right] \\
&= &2 I_p + II_p, 
\eea
where $I_p = (d^\top -\mu^\top(I+\Omega))X$ and $II_p = \hat{\mu}_0^\top (I - {\Omega}) \hat{\mu}_0 -\mu^\top(I-\Omega)\mu$. Let $\tilde{\mu}=U^\top  \mu=(\tilde{\mu}_1, \dots, \tilde{\mu}_p)^\top$, $\tilde{d} = U^\top d$, and $\tilde{x} = U^\top X$, where $W = U \Lambda U^\top$ is the eigenvalue decomposition of $W$. 



\begin{itemize}
\item Consider $I_p = (d -(I+\Omega)\mu)^\top X$. 

Since $d - (I+\Omega)^\top \mu \sim N(0, \frac{1}{n_0}I + \frac{1}{n_1} \Omega)$ and independent with $X$, the expectation $E[I_p] = 0$. 

Next we show that 
\be
\var(I_p) = \frac{1}{n_0} E[X^\top X] + \frac{1}{n_1} E[X^\top \Omega X] \lesssim  4n^{-1}p. 
\ee
The equality comes the law of total variance. We need to check the inequality. 

When $Y = 0$, $X \sim N(-\mu, I)$. According to Lemma \ref{lemma:quad}, 
\ben \ba{lll}
\var(I_p)  & = &\displaystyle \frac{1}{n_0}(p + \|\mu\|^2) + \frac{1}{n_1}(p(1 + \xi) + \mu^\top \Omega \mu)\\
& \leq &\displaystyle (\frac{1}{n_0} + \frac{1}{n_1})p(1+o(1)) + \frac{1}{n_0}\|\Omega\| \|\mu\|^2  \lesssim \displaystyle   \ 4n^{-1}p,
\ea\een
where the last inequality comes from $\|\Omega\| = 1 + \max_{1 \leq i \leq p} m_i = 1+o(1)$ and $\|\mu\|^2 = \tau^2 p \epsilon (1 + o(1))$ with probability $1 - o(1)$. 

When $Y = 1$, $X \sim N(\mu, \Omega^{-1})$. Noting $Tr(\Omega^{-1}) = p(1+o(1))$, 
\ben 
\var(I_p)  =\displaystyle  \frac{1}{n_0}(Tr(\Omega^{-1}) + \|\mu\|^2) + \frac{1}{n_1}(p + \mu^\top \Omega \mu)  \lesssim 
  \displaystyle 4n^{-1}p.
\een
Therefore, for both cases, $\var(I_p) \lesssim 4n^{-1}p$ with probability $1 - o(1)$.

Finally, we give an asymptotic distribution for $I_p$. Element-wise,  $I_p = \sum_{i = 1}^p \tilde{x}_i (\tilde{d}_i - (2 + m_i)\tilde{\mu}_i)$. 
Hence, when $X \sim N(-\mu, I)$, 
\ben\ba{ll}
& \displaystyle [Var(I_p)]^{-3/2}\cdot \overset p{\underset{i=1}\sum} E|\tilde x_i|^3E|\tilde d_i-(2+m_i)\tilde \mu_i|^3 \\
\leq & \displaystyle [Var(I_p)]^{-3/2}\cdot 4 \overset p{\underset{i=1}\sum} (E|\tilde x_i+\tilde\mu_i |^3+|\tilde\mu_i|^3)E|\tilde d_i-(2+m_i)\tilde \mu_i|^3  \\
\leq & \displaystyle [Var(I_p)]^{-3/2}\cdot \frac{8\sqrt 2}{\sqrt \pi}\overset p{\underset{i=1}\sum}  (2\sqrt{2/\pi}+|\tilde \mu_i|^3)\left[ \frac{1}{n_0} + \frac{1 + m_i}{n_1} \right]^{3/2}  \\
\leq & \displaystyle  C \frac{p+\sum_{i=1}^p |\tilde\mu_i|^3}{[\sum_{i=1}^p(1+\tilde \mu_i^2)]^{3/2}} \ \ \ ({\rm for \ some \ constants \ }C>0) \\
\leq & \displaystyle C (p^{-1/2}+ R(\tilde{\mu})).
\ea\een
As a consequence, the Berry-Ess\'een theorem gives $\sup_x|F_{I_p/\sqrt{\var(I_p)}}(x)-\Phi(x)|\overset{\mathcal P} \leq C (p^{-1/2}+ R(\tilde{\mu}))\to 0$, by Lemma \ref{lemma2}. The case is the same when $X \sim N(\mu, \Omega^{-1})$. 
Hence, in both cases, with probability $1 - o(1)$, 
\be\label{eqn:Ip}
I_p \leq C\sqrt{\log\log(p)} \sqrt{4n^{-1}p}.
\ee

\item Consider $II_p = \hat{\mu}_0^\top (I - {\Omega}) \hat{\mu}_0 -\mu^\top(I-\Omega)\mu$. Let $R = \hat{\mu}_0 + \mu$, then $R \sim N(0, \frac{1}{n_0} I_p)$. 
 Therefore, $II_p$ can be rewritten as 
\be
\ba{lll}
II_p & = & (R - \mu)^\top (I - \Omega) (R - \mu) - \mu^\top (I - \Omega) \mu\\
&= & R^\top (I - \Omega) R - 2 \mu^\top (I - \Omega) R\\
& = & IIa + 2IIb. 
\ea
\ee

$IIa$ follows a non-central chi-square distribution. Hence, according to Lemma \ref{lemmaeb}, 
\ben
\ba{rl}
E[IIa] & =  \frac{1}{n_0} Tr(I - \Omega) = -\frac{1}{n_0}\sum_{i=1}^p m_i = -cn^{-1}p\xi, \\
\var(IIa) & =  \frac{2}{n_0^2}Tr((I - \Omega)^2) = \frac{2}{n_0^2}\sum_{i=1}^p m_i^2 = 2c^2n^{-2} (p\xi^2 + p^2 \eta^2 \nu),
\ea
\een
where $c = n/n_0$. 
Further, we can find $\sup_x|F_{IIa/\sqrt{\var(IIa)}}(x)-\Phi(x)|\overset{\mathcal P}\to 0$ by the Berry-Ess\'een theorem, so $|IIa| \leq cn^{-1}p\xi + \sqrt{\log p}(\sqrt{n^{-2} (p\xi^2 + p^2 \eta^2 \nu)})$. 

Now we consider $IIb \sim N(0, \frac{1}{n_0} \mu^\top (I - \Omega)^2 \mu)$.
Furthermore, $\mu^\top (I - \Omega)^2 \mu \leq \|I - \Omega\|^2\|\mu\|^2 \ll \|\mu\|^2 = p\tau^2 \epsilon$. 
Therefore, 
\[
|IIb| \leq \sqrt{\log p}\sqrt{p \tau^2 \epsilon/n}.
\]
Combining the results about $IIa$ and $IIb$, with probability $1 + o(1)$, 
\be\label{eqn:IIp}
|II_p| \leq |IIa| + |IIb| \lesssim c p\xi/n + \sqrt{\log p}(\sqrt{p\xi^2 + p^2 \eta^2 \nu}/n+\sqrt{p\tau^2\epsilon/n}).
\ee
\end{itemize}



Combining the results for $I_p$ and $II_p$ in (\ref{eqn:Ip}) and (\ref{eqn:IIp}), we have 
\ben
 \Delta Q  = 2 I_p + II_p \lesssim  c p\xi/n + \sqrt{\log p}(\sqrt{p\xi^2 + p^2 \eta^2 \nu}/n+\sqrt{p\tau^2\epsilon/n}) + \sqrt{\log \log p}\sqrt{p/n}.
\een
Recall that $\xi \leq p^{-1/2}$ ($\gamma \geq 1/2$ in Theorems \ref{thm:qdaunknownmu} and \ref{thm:qdaunknown}); then $cp\xi/n \leq \sqrt{p}/n = o(\sqrt{p/n})$. 
Since $\xi^2 = o(1)$, $\sqrt{p\xi^2}/n = o(\sqrt{p/n})$. 
Similarly, $\sqrt{p\tau^2\epsilon/n} = o(\sqrt{p/n})$. Therefore, in short, we have 
\be\label{eqn:deltaqpf}
 \Delta Q   \lesssim  \sqrt{\|\Omega - I\|_F^2\log p}/n + \sqrt{p/n \log\log p}.
\ee

The variances $\var(S|Y=k)$ can be found as $\var(S|Y=k) = 2(\|\Omega - I\|_F^2 + 8\|\mu\|^2)(1+o(1))$ according to Lemma \ref{lemma:quad}.

Introduce in (\ref{eqn:deltaqpf}) and $\var(S|Y = k)$, then 
\ben
\left| \frac{\Delta Q}{\sqrt{\var(S_i)}} \right| 
 \leq  \frac{\|\Omega - I\|_F\log p/n + \sqrt{p/n \log\log p}}{\sqrt{\|\Omega - I\|_F^2 + 8\|\mu\|^2}}.
\een
So, the result is proved. \qed

\subsubsection{Proof of Lemma \ref{lemma:weakfail}}
\begin{lemma*}
With probability at least $\Phi(-C)/4$, there is 
\ben
Q(X,\hat{\mu}, \Omega) \geq S + 2C\sqrt{X^\top \Omega X /n} + \mu^{\top}(I - \Omega)\mu+ \log|\Omega|+Tr(I - \Omega)/n_0, 
\een d
where $S$ is defined in (\ref{proofS}).
\end{lemma*}
According to the definition that $d = \Omega \hat{\mu}_1 - \hat{\mu}_0$, there is 
\be
\ba{lll}
Q & = & X^\top (I - \Omega) X + 2 X^{\top} (\Omega \hat{\mu}_1 - \hat{\mu}_0) + \hat{\mu}_0^\top (I - \Omega) \hat{\mu}_0 + \log |\Omega| \\
& = & 2 [X^{\top} \Omega \hat{\mu}_1] + [\hat{\mu}_0^\top (I - \Omega) \hat{\mu}_0 - 2X^{\top} \hat{\mu}_0]  + [X^\top (I - \Omega) X + \log |\Omega|]\\
& = & 2Q_1(\hat{\mu}_1, X) + Q_2(\hat{\mu}_0, X) + Q_3(X). 
\ea
\ee
Given $X$, then $Q_3(X)$ is a constant, and $Q_1(\hat{\mu}_1, X)$ is independent with $Q_2(\hat{\mu}_0, X)$. 

Consider $Q_1(\hat{\mu}_1, X)$, that $Q_1(\hat{\mu}_1, X)|X = X^{\top} \Omega \hat{\mu}_1 | X \sim N(\mu^\top \Omega X, \frac{X^\top \Omega X}{n_1})$. 
Therefore, given $X$, 
\be\label{eqn:q1result}
P(Q_1(\hat{\mu}_1, X) \geq \mu^\top \Omega X +C \sqrt{X^\top \Omega X /n_1}) \geq \Phi(-{C}). 
\ee

Next we consider $Q_2(\hat{\mu}_0, X)|X$. Rewrite it as 
\bea\label{eqn:q2imposs}
Q_2(\hat{\mu}_0, X)& = & (\hat{\mu}_0 - (I - \Omega)^{-1} X)^\top (I - \Omega) (\hat{\mu}_0 - (I - \Omega)^{-1} X) - X^\top (I - \Omega)^{-1} X \nonumber\\
& = & Q_0 - X^\top (I - \Omega)^{-1} X.
\eea
Consider $Q_0$. Note that $\hat{\mu}_0 - (I - \Omega)^{-1} X \sim N(-\mu - (I - \Omega)^{-1} X, \frac{1}{n_0} I)$. 
Let $\mu_X = \sqrt{n_0}[-\mu - (I - \Omega)^{-1} X]$ and $z = \sqrt{n_0}[\hat{\mu}_0 - (I - \Omega)^{-1} X - \mu_X] \sim N(0, I)$, then there is 
\[
Q_0 = Q_0(z) = (\mu_X + z)^\top (I - \Omega) (\mu_X + z)/n_0. 
\]
Simple calculations show that $Q_0(z) + Q_0(-z) = 2 \mu_X^\top (I - \Omega) \mu_X/n_0 + 2 z^\top (I - \Omega) z/n_0$, so $\max\{Q_0(z), Q_0(-z)\} \geq \mu_X^\top (I - \Omega) \mu_X/n_0 + z^\top (I - \Omega) z/n_0$ for any $z$. Since $z$ is symmetric about 0, with probability at least $1/2$, $Q_0(z) \geq \mu_X^\top (I - \Omega) \mu_X/n_0 + z^\top (I - \Omega) z/n_0$. 
 
Further, according to Property 3 of non-central chi-square distribution and Lemma \ref{lemma1}, $P(z^\top (I - \Omega) z/n_0 \geq \frac{1}{n_0} Tr(I - \Omega) ) = \frac{1}{2} + o(1)$ when $p \goto \infty$. 
Therefore, we have that 
\bea\label{eqn:q0imposs}
&&P(Q_0(z) \geq \mu_X^\top (I - \Omega) \mu_X/n_0 + \frac{1}{n_0} Tr(I - \Omega) \\
& = & P(Q_0(z) \geq \mu^\top (I - \Omega) \mu + 2 \mu^\top X + X^\top (I - \Omega)^{-1} X +  \frac{1}{n_0} Tr(I - \Omega))\nonumber\\
 &\geq& 1/4 + o(1). \nonumber
\eea
Combining the above equation (\ref{eqn:q0imposs}) with (\ref{eqn:q2imposs}), then we have, 
\be\label{eqn:q2result}
P(Q_2(\hat{\mu}_0, X) \geq \mu^\top (I - \Omega) \mu + 2 \mu^\top X +   \frac{1}{n_0} Tr(I - \Omega))) \geq 1/4 + o(1).
\ee

Since $Q_1(\hat{\mu}_1, X)$ is independent with $Q_2(\hat{\mu}_0, X)$ when $X$ is given. Hence, combining (\ref{eqn:q1result}) and (\ref{eqn:q2result}), with probability at least $\Phi(-{C})/4$, there is 
\bea
Q & = &  2Q_1(\hat{\mu}_1, X) + Q_2(\hat{\mu}_0, X) + Q_3(X)\nonumber\\
 & \geq&   2 \mu^\top \Omega X + 2C\sqrt{X^\top \Omega X /n_1}+\frac{1}{n_0} Tr(I - \Omega)\nonumber\\
&&  + \mu^{\top}(I - \Omega)\mu + 2\mu^{\top} X + X^\top (I - \Omega) X + \log|\Omega|\nonumber\\
& \geq & 2 \mu^\top (I + \Omega) X + X^\top (I - \Omega) X+ \mu^{\top}(I - \Omega)\mu+ 2C\sqrt{X^\top \Omega X /n}\nonumber\\
&& + \log|\Omega|+\frac{1}{n_0} Tr(I - \Omega)\nonumber\\
& = & S + 2C\sqrt{X^\top \Omega X /n} + \mu^{\top}(I - \Omega)\mu+ \log|\Omega|+\frac{1}{n_0} Tr(I - \Omega), 
\eea
where $S$ is defined in (\ref{proofS}). 

The result is proved. \qed


\subsubsection{Proof of Lemma \ref{lemma:trace}}
\begin{lemma*}
Consider two symmetric matrices $A$ and $B$, where the eigenvalues of $B$ are $\lambda_1$, $\cdots, \lambda_p$; then,
\[
Tr(AB) \leq \|A\| *\sum_{i=1}^p |\lambda_i|.
\]
Furthermore, if $B$ is a positive semi-definite matrix,
\[
Tr(AB) \leq Tr(B)\|A\|.
\]
\end{lemma*}

Let the eigen-value decomposition of $B = U\Lambda U^\top$. Then we have 
\[
Tr(AB) = Tr(AU\Lambda U^{\top}) = Tr(U^{\top}AU\Lambda) 
= Tr(\tilde{A}\Lambda),
\]
where $\tilde{A} = U^{\top} AU$. $U$ is an orthogonal matrix so $\|\tilde{A}\| = \|A\|$. 

Therefore, we have 
\[
Tr(\tilde{A}\Lambda) = 
\sum_{i=1}^p \tilde{A}(i,i)\lambda_i
\leq \sum_{i=1}^p \|A\||\lambda_i| 
= \|A\|\sum_{i=1}^p |\lambda_i|.
\]

When $B$ is positive semi-definite, all the eigenvalues of $B$ are non-negative and $|\lambda_i| = \lambda_i$, $i = 1, 2, \cdots, p$. Therefore, 
\[
Tr(AB) \leq \|A\| *\sum_{i=1}^p |\lambda_i|
 = \|A\| *\sum_{i=1}^p \lambda_i = \|A\| Tr(B).
\]



\section{The normality of quadratic forms}\label{app:seclemma}
In this section, we will discuss how the quadratic form satisfy the conditions in Lemma \ref{lemma:quad} so that it is asymptotically normal distributed. We first introduce two lemmas that is helpful and then show the proof. 

We first introduce some notations for this section. For a fixed $\Omega$, let the eigenvalue decomposition of $\Omega$ be 
$ \Omega=U\Lambda U^\top$,  
 where $U^\top U=U U^\top=I$ and $\Lambda=diag(\lambda_1, \lambda_2, \dots, \lambda_p)$. 
Recall $\Omega = D_{\Omega} + \eta W$. Let $m_i$ be the eigenvalues of $\Omega - I$. There is $m_i = \lambda_i - 1$, $1 \leq i \leq p$.  For any vector $a = (a_1, a_2, \cdots, a_p)$, define $R(a) = \sqrt{p}\sum_{i=1}^p |a_i|^3/\|a\|^3$.

\subsection{Proof of Lemma \ref{lemma1b}}
In this section, we will prove the following lemma about $R(W)$, which can later show the asymptotic normality of the quadratic function. 
\begin{lemma}\label{lemma1b}
Let $R(\Omega) = R((m_1,m_2, \cdots, m_p))$ be the function of the eigenvalues of $\Omega - I$. Under models (\ref{Omega1}) and (\ref{Param}), when $p\to\infty$, with probability $1-o(1)$, 
\ben
R(\Omega) /\sqrt{p} \longrightarrow 0. 
\een
\end{lemma}

\begin{itemize}
\item Not very sparse region ($0<\beta < 1$). 

In this region, almost surely, $|m_i| \leq \xi + 2\eta\sqrt{p \nu}(1+o(1)) \leq \xi + 3 \eta\sqrt{p \nu}$. 

Now we analyze $R(\Omega)$. Because $\sum m_i^2 = \|\Omega - I\|^2 \sim p\xi^2 + 2\eta^2 Binomial(p(p-1)/2, \nu)$. According to the Chernoff lower tail bound, 
\begin{equation}\label{eqn:dense2}
P(\|\Omega - I\|^2 \leq p\xi^2 + \eta^2 p(p-1)\nu/2)\leq e^{-p(p-1)\nu/16}. 
\end{equation}
Therefore, with high probability that at least $1 - e^{-p(p-1)\nu/16}$, $\|\Omega - I\|^2 > p\xi^2 + \eta^2 p(p-1)\nu/2$. 

Combining this with the result for $\|\Omega\|$, with probability tending to 1, we have
\begin{equation}\label{eqn:dense3}
R(\Omega) \leq \frac{(\xi+3\eta\sqrt{p\nu})^3}{[(p\xi^2+\eta^2 p(p-1)\nu/2)/p]^{3/2}} = \frac{54\sqrt 2 (\xi^2 + \eta^2p\nu)^{3/2}}{[(\xi^2+\eta^2 p(p-1)\nu/2)/p]^{3/2}} \leq 216.
\end{equation}
\item Sparse region ($1 \leq \beta < 2$). 

\begin{itemize}
\item When $\beta = 1$ ($\nu = p^{-1}$). 
Let $S = \sum_{i=1}^p m_i^2$.  Note that
\begin{equation}\label{eqn:sparse2}
R(\Omega) = \frac{\frac{1}{p}\sum_{i=1}^p |m_i|^3}{(\frac{1}{p}\sum_{i=1}^p m_i^2)^{3/2}} \leq \frac{54\sqrt 2 (\xi^2 + \eta^2\frac{\log p}{\log\log p})^{3/2}}{[(\xi^2+\eta^2 p(p-1)\nu/2)/p]^{3/2}},
\end{equation}
where $\|\Omega - I\|^2 > p\xi^2 + \eta^2 p(p-1)\nu/2$ with high probability that at least $1 - e^{-p(p-1)\nu/16}$, according to (\ref{eqn:dense2}).
Therefore, \begin{equation}\label{eqn:sparse3}
R(\Omega) \leq 4 \sqrt{\frac{\log p}{\log\log p}}.
\end{equation}

\item When $1 < \beta \leq 2$ ($p^{-2} \leq \nu < p^{-1}$). 

Recall that $\Omega = D_{\Omega} + V$. For the current case, there are very few non-zeros in $V$, which actually at the same order of non-zero eigenvalues in $V$. 
Let $S$ be the set of all the rows that has non-zero off-diagonal entries. Then $\Omega^{(S)}$, the sub-matrix of $\Omega$ restricted on $S \times S$, contains all the non-zeros off-diagonals. Similarly, if we define $\Omega^{(S^c)}$ where $\Omega$ is restricted on $S^c \times S^c$, then $\Omega^{(S^c)}$ is an diagonal matrix. The leftover sub-matrices are zero matrices. 

The eigenvalues of $\Omega$ can be decomposed as the union of two sets $\{m_i^{(S)}\}$ and $\{m_i^{(S^c)}\}$, where $m_i^{(S)}$ are the eigenvalues of $\Omega^{(S)}$. Hence,  
\[
|m_i^{(S)}|\leq \xi + \eta/(\beta - 1), \quad 1 \leq i \leq |S|;\qquad 
|m_i^{(S^c)}| = \xi, \quad 1 \leq i \leq |S^c|. 
\]

For each row, there at most two non-zeros because $\nu < p^{-1}$. Therefore, the number of all non-zero entries is $\leq 2|S|$. 
Similar to (\ref{eqn:sparse2}), 
\begin{eqnarray}\label{eqn:sparse5}
R(\Omega) & \leq & \frac{\xi^3 + \frac{8}{|\beta - 1|^3}\eta^3|S|/p}{(\frac{1}{p}\sum_{i=1}^p m_i^2)^{3/2}}\nonumber\\
& \leq & \frac{\xi^3 + \frac{8}{|\beta - 1|^3}\eta^3|S|/p}{(\xi^2 + 2\eta^2|S|/p)^{3/2}}\nonumber\\
& \leq & 1 + \frac{\frac{1}{p} \cdot |S| \cdot |\frac{2}{\beta-1}|^3}{(\frac{1}{p}|S|)^{3/2}}\nonumber\\
&=& 1 + |\frac{2}{\beta-1}|^3 (|S|/p)^{-1/2}.
\end{eqnarray}

Since, according to (\ref{eqn:dense2}), $S>p(p-1)\nu/2$ with high probability (larger than $1 - e^{-p(p-1)\nu/16}$), there is high probability that
\begin{equation}
R(\Omega) \leq 1 + |\frac{2}{\beta-1}|^3 ((p-1)\nu/2)^{-1/2} \leq 1 + C (p\nu) ^{-1/2}.
\end{equation}
\end{itemize}
Combining the three cases, we have that 
\ben
R(\Omega) \leq r(p, \beta) = \left\{\begin{array}{ll}
216, & 0<\beta < 1,\\
4\sqrt{{\log p}/{\log\log p}}, & \beta=1, \\
1 + Cp^{(\beta-1)/2}, & 1 < \beta \leq 2.
\end{array}
\right.
\een
Hence, $R(\Omega)/\sqrt{p} \goto 0$ with probability $1 - o(1)$. 
\end{itemize}

\subsection{Proof of Lemma \ref{lemma2}}\label{sec:lemma2}
Recall that $\tilde{\mu} = U^T \mu$ and $\tilde{\mu}_i = u_i^T\mu$, where $U$ is the orthonormal matrix from the eigendecomposition $\Omega = U\Lambda U^T$. 

\begin{lemma}\label{lemma2} 
As $p\to\infty$, under models (\ref{Param1}), (\ref{Param2}), (\ref{Param}), and (\ref{con1}), with probability $1 - o(1)$, 
\ben
R(\tilde{\mu})/\sqrt{p}  \ \longrightarrow \ 0.
\een
\end{lemma}
{\it Proof}. 
First, we figure out the bound for fixed $\mu$ and $U$. Note that the entries of $\mu$ are either 0 or $\tau$. Let $S(\mu) = \{i: \mu_i \neq 0\}$, and $|S(\mu)|$ be the cardinality of $S(\mu)$. Hence, 
\be\label{lemma2:mui}
\tilde{\mu}_i = u_i^T\mu = \tau \sum_{j \in S(\mu)} u_{ij}.
\ee
We introduce it in our target and rewrite it as
\begin{eqnarray}
\frac{1}{\sqrt p} \cdot \frac{\frac{1}{p}\sum_{i=1}^p |\tilde{\mu}_i|^3}{\left(\frac{1}{p}\sum_{i=1}^p |\tilde{\mu}_i|^2\right)^{3/2}} 
& = & 
 \frac{\sum_{i=1}^p |\tilde{\mu}_i|^3}{\left(\sum_{i=1}^p |\tilde{\mu}_i|^2\right)^{3/2}} 
= \frac{\tau^3 \sum_{i=1}^p |\sum_{j \in S(\mu)} u_{ij}|^3}{\left(\sum_{i=1}^p \tau^2 |\sum_{j \in S(\mu)} u_{ij}|^2\right)^{3/2}} 
\nonumber\\
& = &  \frac{\sum_{i=1}^p |\sum_{j \in S(\mu)} u_{ij}|^3}{\left(\sum_{i=1}^p |\sum_{j \in S(\mu)} u_{ij}|^2\right)^{3/2}}
\end{eqnarray}
Let $a_i = \sum_{j \in S(\mu)} u_{ij}$. Then, $\sum |a_i|^3 \leq \max_{1 \leq i \leq p} |a_i| \cdot \sum a_i^2 $. Hence, we can further revise the equation as follows:
\be\label{step1}
\frac{1}{\sqrt p} \cdot \frac{\frac{1}{p}\sum_{i=1}^p |\tilde{\mu}_i|^3}{\left(\frac{1}{p}\sum_{i=1}^p |\tilde{\mu}_i|^2\right)^{3/2}} 
 \leq  \frac{\max_{1 \leq i \leq p} |a_i|}{\left(\sum_{i=1}^p |\sum_{j \in S(\mu)} u_{ij}|^2\right)^{1/2}} = \frac{\max_{1 \leq i \leq p} |a_i|}{|S(\mu)|^{1/2}}.
\ee
That last equation (above) comes from
\be
\sum_{i=1}^p |\tilde{\mu}_i|^2 = \|\mu_-\|^2_2 = \|U^T \mu\|^2_2 = \|\mu\|^2_2 = \tau^2 |S(\mu)|.
\ee

Now, we apply the randomness of $\mu$ to control the upper bound in (\ref{step1}). Fix a constant $0 < c < \min\{\beta/4, \frac{1-\beta}{2}\}$. Since $\beta < 1$, a positive constant always exists. We define $0 < b_0 < b_1 < \cdots < b_K$, where 
\[
b_0 = \frac{1}{\sqrt{p\epsilon}} p^{-c} = p^{\frac{1-\beta}{2}-c}, \quad 
b_{i} = b_{i-1}\cdot p^{\beta - 2c}, \qquad 1 \leq i \leq K,
\]
and $K = \min\{k; b_k \geq 1\}$. Note that $b_i$ increases by $p^{\beta - 2c} \goto \infty$, such that $k$ always exists and does not change with respect to $p$. The increasing step also indicates that there is at most one $k$, so  $p^{-\beta/2} < b_k < 1$. If such a $k$ exists, $b_{k+1} = 1$. 

With the definition of vector $a$, we define a matrix $V = (v_{ij})$ as follows:
\be\label{defineV}
v_{ij} = \left\{\begin{array}{ll}
b_0, & |u_{ij}| \leq b_0;\\
b_i, & b_{i-1} < |u_{ij}| \leq b_i.\\
\end{array}
\right.
\ee
Since $0 \leq |u_{ij}| \leq 1$, $V$ is well defined. 
Clearly, $v_{ij} \geq u_{ij}$ holds for each entry, which indicates that
\be\label{utov}
\max_{1 \leq i \leq p}|a_i| = \max_{1 \leq i \leq p}|\sum_{j\in S(\mu)} u_{ij}| 
\leq  \max_{1 \leq i \leq p}\sum_{j\in S(\mu)} v_{ij}.
\ee
The bound in (\ref{step1}) is correspondingly replaced by $\max_{1 \leq i \leq p}\sum_{j\in S(\mu)} v_{ij}/|S(\mu)|^{1/2}$. 
According to the definitions of $V(i,j)$, the bound can be decomposed into $k+1$ parts. 
\begin{align}
&\max_{1 \leq i \leq p}\frac{\sum_{j\in S(\mu)} v_{ij}}{|S(\mu)|^{1/2} }\nonumber\\
 \leq & \frac{b_0 \max_i\sum_{j \in S(\mu)} 1\{|u_{ij}|\leq b_0\}}{|S(\mu)|^{1/2} } + \sum_{k = 1}^K \frac{b_k \max_i\sum_{j \in S(\mu)} 1\{b_{k-1} < |u_{ij}|\leq b_k\}}{|S(\mu)|^{1/2} }\nonumber\\
 =  & I + \sum_{k=1}^K IIk.\label{decompb}
\end{align}

Consider term $I$. Note that $\sum_{j \in S(\mu)}1\{|u_{ij}|\leq b_0\} \leq |S(\mu)|$ holds for any $1 \leq i \leq p$, and therefore
\be
I \leq \frac{b_0 |S(\mu)|}{|S(\mu)|^{1/2}} = b_0|S(\mu)|^{1/2} = p^{-\frac{1-\beta}{2} - c}|S(\mu)|^{1/2}. 
\ee
Let $\nu_i = 1\{\mu_i \neq 0\} - \epsilon$. Then, $E[\nu_i] = 0$, $Var(\nu_i) = \epsilon(1 - \epsilon)$, and $\max(|\nu_i|) < 1$. Hence, $E[\sum \nu_i] = |S(\mu)| - p\epsilon$. According to Bernstein's inequality (\cite{shorack2009empirical}), 
\begin{eqnarray*}
P(|S(\mu)| \geq p\epsilon p^{c})  & = & P(|S(\mu)| - p\epsilon \geq p\epsilon (p^{c}-1)) \leq \exp\{\frac{-p^2\epsilon^2(p^c - 1)^2}{2p\epsilon(1 - \epsilon) + 2/3p\epsilon(p^c - 1)}\}\\
& \leq & \exp\{-\frac{p\epsilon(p^c - 1)^2}{2(1 - \epsilon) + 2/3(p^c - 1)}\}\\
& \leq & \exp\{-3p\epsilon(p^c - 1)/4\} = o(p^{-1}).
\end{eqnarray*}
Hence, with probability $1 - o(p^{-1})$, 
\be\label{term1}
I \leq p^{-\frac{1-\beta}{2} - c}|p\epsilon p^{c}|^{1/2} = p^{-c/2}  \goto 0. 
\ee

Now we consider term $IIk$. Since $U$ is an orthonormal matrix, $\|u_i\|^2 = 1$, and, hence,
\[
b_{k-1}^2 \sum_{j=1}^p 1\{b_{k-1} \leq |u_{ij}| \leq b_k\}  \leq \sum_{j=1}^p u_{ij}^2 = 1 
\Longrightarrow  \sum_{j=1}^p 1\{b_{k-1} \leq |u_{ij}| \leq b_k\} \leq 1/b_{k-1}^2.
\]
Therefore, there are at most $1/b_{k-1}^2$ elements in the region $[b_{k-1}, b_k]$. Since $1\{b_{k-1} \leq |u_{ij}| \leq b_k\}$ is independent of $\mu$, each element in $[b_{k-1}, b_k]$ is selected independently with probability $\epsilon$. In addition, let $Y_i + \epsilon \sim Bernoulli(1, \epsilon)$, with mean $0$ and variance $\epsilon(1-\epsilon)$, and let $Y = \sum_{i = 1}^{1/b_{k-1}^2} Y_i$. 

We discuss the control of $IIk$ according to the magnitude of $b_{k-1}$. 
\begin{itemize}
\item {\it Case 1. $b_{k-1} < p^{-\beta/2} = \sqrt{\epsilon}$. } With Bernstein's inequality, 
\bea
&& P(\sum_{j \in S(\mu)} 1\{b_{k-1} \leq |u_{ij}| \leq b_k\} > \frac{\epsilon}{b_{k-1}^2}p^{c/4})\nonumber\\
& \leq & P(Y > \frac{\epsilon}{b_{k-1}^2}p^{c/4}) =  P(Y-\epsilon/b_{k-1}^2 > \frac{\epsilon}{b_{k-1}^2}(p^{c/4} -1)) \nonumber \\
& \leq & \exp\{\frac{-\epsilon^2(p^{c/4}-1)^2/b_{k-1}^4 }{2\epsilon(1-\epsilon)/b_{k-1}^2 + 2(p^{c/4}-1)\epsilon/(3b_{k-1}^2) }\}\nonumber\\
& \leq & \exp\{-\epsilon(p^{c/4}-1)/b_{k-1}^2\} \qquad (\mbox{since } p^{c/4} \goto \infty)\nonumber
\eea
If $b_{k-1} < p^{-\beta/2} = \sqrt{\epsilon}$, $\epsilon(p^{c/4}-1)/b_{k-1}^2 \geq p^{c/4}-1 \goto \infty$, and the result is reduced to 
\be\label{term2part1}
P(\sum_{j \in S(\mu)} 1\{b_{k-1} \leq |u_{ij}| \leq b_k\} > \frac{\epsilon}{b_{k-1}^2}p^{c/4}) 
\leq  \exp\{-\epsilon(p^{c/4}-1)/b_{k-1}^2\}= o(p^{-2}).
\ee
Therefore, 
\bea
&& P(\max_i\sum_{j \in S(\mu)} 1\{b_{k-1} \leq |u_{ij}| \leq b_k\} > \frac{\epsilon}{b_{k-1}^2}p^{c/4})\nonumber\\
& \leq & p P(\sum_{j \in S(\mu)} 1\{b_{k-1} \leq |u_{ij}| \leq b_k\} > \frac{\epsilon}{b_{k-1}^2}p^{c/4}) = o(p^{-1}).
\eea
Hence, with probability $1 -o(p^{-1})$, 
\bea\label{term2part1i}
b_{k}\max_i\sum_{j \in S(\mu)} 1\{b_{k-1} \leq |u_{ij}| \leq b_k\} & \leq & b_k \frac{\epsilon}{b_{k-1}^2}p^{c/4}\nonumber\\
& \leq & b_{k-1} \cdot p^{\beta - 2c} \frac{\epsilon}{b_{k-1}^2} p^{c/4}\nonumber\\
& = & \frac{p^{7c/4}}{b_{k-1}} \leq p^{-7c/4}/b_0 \nonumber\\
& = & \sqrt{p\epsilon} p^{-3c/4}.
\eea
\item {\it Case 2. $b_{k-1} > p^{-\beta/2} = \sqrt{\epsilon}$.} Note that, in this case, $\epsilon/b_{k-1}^2 \leq 1$ and $b_{k+1} = 1$. As with the derivative in case 1, using Bernstein's inequality, we have
\bea
&& P(\sum_{j \in S(\mu)} 1\{b_{k-1} \leq |u_{ij}| \leq b_k\} > \sqrt{p\epsilon p^{-2c}})\nonumber\\
& \leq &  P(Y-\epsilon/b_{k-1}^2 > \sqrt{p\epsilon p^{-2c}} - \frac{\epsilon}{b_{k-1}^2}) \nonumber \\
& \leq & \exp\{\frac{-(\sqrt{p\epsilon p^{-2c}} - \epsilon/b_{k-1}^2)^2 }{2\epsilon(1-\epsilon)/b_{k-1}^2 + 2(\sqrt{p\epsilon p^{-2c}} - \epsilon/b_{k-1}^2)/3 }\}\nonumber\\
& \leq & \exp\{-(\sqrt{p\epsilon p^{-2c}} - \epsilon/b_{k-1}^2)\} = o(p^{-2})\qquad (\mbox{since } p\epsilon p^{-2c} \goto \infty)\nonumber
\eea
Introducing the result into the following maximum question, we get
\be\label{term2part2}
P(\max_i\sum_{j \in S(\mu)} 1\{b_{k-1} \leq |u_{ij}| \leq b_k\} > \sqrt{p\epsilon p^{-2c}} ) = o(p^{-1}).
\ee
Hence, with probability $1- o(p^{-1})$, 
\be\label{term2part2i}
b_{k}\max_i\sum_{j \in S(\mu)} 1\{b_{k-1} \leq |u_{ij}| \leq b_k\} \leq \sqrt{p\epsilon p^{-2c}} \leq \sqrt{p\epsilon} p^{-3c/4}.
\ee

\end{itemize}

According to Bernstein's inequality, 
\begin{eqnarray}\label{term2denom}
P(|S(\mu)| \leq p\epsilon p^{-c}) & = & P(|S(\mu)| - p\epsilon \leq p\epsilon (p^{-c}-1) )\nonumber\\
&\leq& \exp\{\frac{-p^2\epsilon^2(p^{-c} - 1)^2}{2p\epsilon(1 - \epsilon) + 2/3p\epsilon(1 - p^{-c})}\}\nonumber\\
&\leq& \exp\{\frac{-p^2\epsilon^2(p^{-c} - 1)^2}{3p\epsilon}\} = \exp\{-p\epsilon(p^{-c} - 1)^2/3\} = o(p^{-1})
\end{eqnarray}
Introducing (\ref{term2part1i}), (\ref{term2part2i}), and (\ref{term2denom}) into $IIk$, with probability $1 - o(p^{-1})$, 
\be\label{term2}
IIk \leq \frac{ \sqrt{p\epsilon}p^{-3c/4} }{(p\epsilon p^{-c})^{1/2} } = p^{-c/4} \goto 0. 
\ee

Introducing (\ref{term1}) and (\ref{term2}) into (\ref{decompb}), with probability $1 - o(1)$,  
\[
\max_{1 \leq i \leq p}\frac{\sum_{j\in S(\mu)} v_{ij}}{|S(\mu)|^{1/2} }.
\]
Combining it with (\ref{step1}) and (\ref{utov}), Lemma \ref{lemma2} is proved. \qed

\subsection{The normality of quadratic forms}
We discussed a lot about the quadratic form. In this section, we prove the normality of 
\[
S = X^{\top} (\Omega - I) X + 2\mu^\top (\Omega + I) X. 
\]
\begin{lemma}\label{lemma:S}
Let $S_i = S|Y = i$, $i = 0, 1$. Define $Z_i = [S_i - E[S_i]]/\sqrt{\var(S_i)}$ and $F_{Z_i}(x) = P(Z_i \leq x)$, then 
\[
\sup\nolimits_{-\infty < x < \infty}|F_{Z_i}(x) - \Phi(x)| \goto 0.
\]
\end{lemma}
Let $m_i$ be the eigenvalues of $\Omega - I$ and $U$ be the matrix containing eigenvectors of $\Omega - I$. 
Let $\tilde{\mu}=U^\top \mu=(\tilde{\mu}_1, \dots, \tilde{\mu}_p)^\top$. Let $\tilde{d}=U^\top (\Omega + I)\mu=(\tilde{d}_1, \dots, \tilde{d}_p)^\top$. Then $\tilde{d}_i = (2+m_i)\tilde{\mu}_i$
We prove the normality for the case $Y = 0$ and $Y = 1$. For each case, we can apply Lemma \ref{lemma:quad} and check the bound of $|F_Z(x) - \Phi(x)|$. 

First, consider the case $Y = 0$.  In this case, $X \sim N(-\mu, I)$.
According to Lemma \ref{lemma:quad}, we only need to show
\[
\frac{\sum_{i=1}^p |m_i|^3(1 + |\tilde{\mu}_i|^3 + |\tilde{d}_i/m_i|^3)}{(\sum_{i=1}^p m_i^2((1 + |\tilde{\mu}_i|^3 + |\tilde{d}_i/m_i|^2)))^{3/2}} 
\leq
\frac{\sum_{i=1}^p |m_i|^3}{(\sum_{i=1}^p m_i^2)^{3/2}} +  \frac{\sum_{i=1}^p |2(1+m_i)\tilde{\mu}_i|^3}{[\sum_{i=1}^p 4(1+m_i)^2\tilde{\mu}_i^2]^{3/2}} \goto 0.
\]
For the first term, by Lemma \ref{lemma1b}, we have 
\[
\frac{\sum_{i=1}^p |m_i|^3}{(\sum_{i=1}^p m_i^2)^{3/2}}  = R(\Omega)/\sqrt{p} \goto 0, 
\]
 Let $\omega = \max_{1\leq i \leq p} |m_i| \goto 0$. Then, for the second term, by Lemma \ref{lemma2}, we have 
\be\label{w0normal}
\ba{lll}
\frac{\sum_{i=1}^p |(1+m_i)\tilde{\mu}_i|^3}{[\sum_{i=1}^p (1+m_i)^2\tilde{\mu}_i^2]^{3/2}}& \leq &  
\frac{\overset p {\underset {i=1}\sum}(1+m_i)^3|\tilde{\mu}_i|^3}{\left[\overset p {\underset {i=1}\sum} (1+m_i)^2 \tilde{\mu}_i^2\right]^{3/2}} \\
& \leq & \frac{(1+\omega)^3}{(1-\omega)^3}  \frac{\overset p {\underset {i=1}\sum}|\tilde{\mu}_i|^3}{\left(\overset p {\underset {i=1}\sum} \tilde{\mu}_i^2\right)^{3/2}} \\
& = &  \frac{(1+\omega)^3}{(1-\omega)^3} R(\tilde{\mu})/\sqrt{p} \goto 0.
\ea
\ee

For the case $Y = 1$, the analysis is similar. In this case, $X \sim N(\mu, \Omega^{-1})$. 
Hence, using Lemma \ref{lemma:quad}, we only need to show 
\be
\ba{l}
36^2\sqrt{2}\frac{\sum_{i=1}^n \frac{|m_i|^3}{|1+m_i|^3}}{[\sum_{i=1}^n \frac{m_i^2}{(1+m_i)^2}]^{3/2}} + \frac{36^2}{\sqrt{2\pi}} \frac{\sum_{i=1}^n |2\frac{\tilde{\mu}_i}{\sqrt{1+m_i}}|^3}{(\sum_{i=1}^n 4\frac{\tilde{\mu}_i^2}{1+m_i})^{3/2}}\\
\qquad\qquad \leq   \frac{(1+\omega)^3}{(1-\omega)^3}\frac{72^2\sqrt 2}{\sqrt p} \left[1+R(W)\right]+ \frac{36^2}{\sqrt{2\pi}}\frac{(1+\omega)^3}{(1-\omega)^3}  R(\tilde{\mu}).\\
\qquad\qquad \leq C(1/\sqrt{p} + R(W)/\sqrt{p} + R(\tilde{\mu})/\sqrt{p}). 
\ea
\ee
According to Lemma \ref{lemma1b} and Lemma \ref{lemma2}, $C((1 + R(W))/\sqrt{p} + R(\tilde\mu)) = o(1)$ with probability $1 - o(1)$. 
So, the result is proved. \qed

\section{Proof of Lemmas}\label{app:lemma}
In this section, we prove the lemmas appeared in the main paper. For the convenience of readers, we copy paste the lemmas to prove in each section. 

\subsection{Proof of Lemma \ref{lemma:quad}}
\begin{lemma*}
Consider $X \sim N(\mu, \Sigma)$ where $\Sigma$ is positive definite. Let $S = X^\top A X + 2d^\top X$ with a symmetric matrix $A$ and a vector $d$,  
\bea
&& E[S] = Tr(A\Sigma) + \mu^\top A \mu + 2d^\top \mu,\\
&& \var(S) = 2Tr((A\Sigma)^2) + 4(\mu^\top A\Sigma A\mu + \mu^\top A\Sigma d + d^\top \Sigma d).
\eea
Further, let $|\lambda_1| \geq |\lambda_2| \geq \cdots \geq |\lambda_p|$ be the eigenvalues of $\Sigma^{1/2}A\Sigma^{1/2}$ and $u_1, u_2, \cdots, u_p$ be the corresponding eigenvectors. Let $W = \frac{z - E[z]}{Var(z)}$ with cumulative density function $F_{W}(x)$. Let $\Phi(x)$ be the cumulative density function of standard normal distribution. Let $\tilde{\mu}(i) = u_i^\top (\Sigma^{-1/2}\mu - \Sigma^{1/2}d/\lambda_i)$ for $\lambda_i \neq 0$ and $\Tilde{d}(i) = u_i^\top \Sigma^{1/2}d$. Then $\sup_x |F_{W}(x) - \Phi(x)| = o(1)$ if one of the following conditions hold, 
\begin{itemize}
    \item[(a)] $\frac{\sum_{i=1}^p |\lambda_i|^3 (1 + |\tilde{\mu}(i)|^3)}{(\sum_{i=1}^p \lambda_i^2(1 + \tilde{\mu}_i^2))^{3/2}} \goto 0$; or 
    \item[(b)] $\var(S) = \sum_{i:\lambda_i = 0} \Tilde{d}^2(i)(1+o(1))$. 
\end{itemize}
\end{lemma*}

We are considering a quadratic form $S$ about $X$, where $X$ follows multivariate normal distribution. Therefore, according to the analysis in \cite{searle2016linear}, we have the mean and variance about the quadratic form, where 
\bea
&& E[S] = Tr(A\Sigma) + \mu^\top A \mu + 2d^\top \mu,\\
&& \var(S) = 2Tr((A\Sigma)^2) + 4(\mu^\top A\Sigma A\mu + \mu^\top A\Sigma d + d^\top \Sigma d).
\eea

Now we prove the asymptotic distribution. 
Let $U = (u_1, u_2, \cdots, u_p)$, then $\Sigma^{1/2}A\Sigma^{1/2} = U^\top \Lambda U$. Let $Y = \Sigma^{-1/2} X$, then $Y \sim N(\Sigma^{-1/2}\mu, I)$. 
If all the eigenvalues of $A$ are nonzero, we have
\bea
S & = & X^\top A X + 2d^\top X \nonumber\\
& = & (X - A^{-1}d)^{\top} A (X - A^{-1}d)\nonumber\\
& = &(Y - \Sigma^{-1/2}A^{-1}d)^{\top} \Sigma^{1/2} A \Sigma^{1/2} (Y - \Sigma^{-1/2}A^{-1}d)\nonumber\\
& = & (Y - \Sigma^{-1/2}A^{-1}d)^{\top} U^\top \Lambda U (Y - \Sigma^{-1/2}A^{-1}d). \nonumber
\eea
Because $Y \sim N(\Sigma^{-1/2}\mu, I)$, there is 
$U(Y - \Sigma^{-1/2}A^{-1}d) \sim N(U(\Sigma^{-1/2}\mu -  \Sigma^{-1/2}A^{-1}d), I)$. 
Now consider $\Sigma^{-1/2}A^{-1}d$. Recall that $\Sigma^{1/2}A\Sigma^{1/2} = U^\top \Lambda U$, so $\Sigma^{-1/2}A^{-1}\Sigma^{-1/2} = U^\top \Lambda^{-1} U$. Introduce it in and we have $\Sigma^{-1/2}A^{-1}d = U^\top \Lambda^{-1} U \Sigma^{1/2} d$. 
Therefore, the mean follows that 
\[
U(\Sigma^{-1/2}\mu -  \Sigma^{-1/2}A^{-1}d) 
= U\Sigma^{-1/2}\mu - \Lambda^{-1} U \Sigma^{1/2} d
= \tilde{\mu}.
\]
Hence, 
\[
S = (Y - \Sigma^{-1/2}A^{-1}d)^{\top} U^\top \Lambda U (Y - \Sigma^{-1/2}A^{-1}d) = \sum_{i=1}^p \lambda_i (U (Y - \Sigma^{-1/2}A^{-1}d))_i^2,
\]
where $U (Y - \Sigma^{-1/2}A^{-1}d) \sim N_p(\tilde{\mu}, I)$. Let $\tilde{Y} = U (Y - \Sigma^{-1/2}A^{-1}d) - \tilde{\mu}$. Then, $\tilde{Y} \sim N_p(0, I)$ and $z = \sum_{i=1}^n \lambda_i (\tilde{Y}_i + \tilde{\mu}_i)^2$. Since
\ben \ba{lll}
\displaystyle |\lambda_i|^3E\left|(\tilde{Y}_i + \tilde{\mu}_i)^2 - E[(\tilde{Y}_i + \tilde{\mu}_i)^2]\right|^3 & = &  |\lambda_i|^3E\left|\tilde{Y}_i^2+2\tilde{\mu}_i \tilde{Y}_i-1\right|^3\\
 & \leq &  |\lambda_i|^3\bigl[ 9E|\tilde{Y}_i|^6+72|\tilde{\mu}_i|^3E|\tilde{Y}_i|^3+9\bigr] \\
& = &  |\lambda_i|^3\bigl[9(15)+72|\tilde{\mu}_i|^3\left(2\sqrt{\frac{2}{\pi}}\right)+9\bigr] \\
& = & \displaystyle 144|\lambda_i|^3+144 \sqrt{\frac{2}{\pi}}|\lambda_i|^3|\tilde{\mu}_i|^3,
\ea
\een
according to the Berry-Ess\'een theorem, for given $\mu$ and $W$,
\[
\ba{lll}
\displaystyle \underset x \sup\left|F_{W}(x)-\Phi(x) \right|  & \leq & \displaystyle 36\frac{144\sum_{i=1}^p |\lambda_i|^3+144 \sqrt{\frac{2}{\pi}}\sum_{i=1}^p|\lambda_i|^3|\tilde{\mu}_i|^3}
{2\sqrt 2 \left\{\sum_{i=1}^p  \left(\lambda_i^2 +2\lambda_i^2\tilde{\mu}_i^2 \right) \right\}^{3/2}}\\
& \leq & \displaystyle 36^2\sqrt{2}\frac{\sum_{i=1}^p |\lambda_i|^3}{(\sum_{i=1}^p \lambda_i^2)^{3/2}} + \frac{36^2}{\sqrt{2\pi}} \frac{\sum_{i=1}^p |\lambda_i\tilde\mu_i |^3}{(\sum_{i=1}^p \lambda_i^2 \tilde{\mu}_i^2)^{3/2}}.
\ea
\]
The result is thus proved, and we have control over the distribution of $z$. 

Now we consider the case that there are multiple zeros in $\lambda_i$. Again, let $Y = \Sigma^{-1/2}X$ and $Z = UY$, then $Y \sim N(\Sigma^{-1/2}\mu, I)$ and $Z \sim N(U\Sigma^{-1/2}\mu, I)$.
The term $S$ then follows
\bea
S & = & Y^{\top} \Sigma^{1/2}A \Sigma^{1/2} Y + 2(\Sigma^{1/2}d)^{\top} Y\\
& = & Z^{\top} \Lambda Z + 2(U\Sigma^{1/2}d)^{\top} Z\\
& = & \sum_{i=1}^p (\lambda_i z_i^2 + 2\tilde{d}_i z_i). 
\eea
When there are multiple zeros in $\lambda_i$, then the terms to sum is $\sum_{i: \lambda_i = 0} 2\tilde{d}_i z_i \sim N(2\tilde{d}_iE[z_i], 4\sum_{i: \lambda_i = 0}\tilde{d}_i^2)$. It follows normal because every item is independent normal distribution. 

Hence, when $\var(S) = \sum_{i:\lambda_i = 0} \tilde{d}^2(i)(1+o(1))$, then $\frac{S - E[S]}{\var(S)}$ converges to normal distribution.  
\qed

\subsection{Proof of Lemma \ref{lemma:tweak}}
\begin{lemma*}
Under current model and assumptions, for a given symmetric matrix $A$, there exists a constant $C > 0$, so that with probability $1 - o(1)$, 
\be
\left|\left[\hat{\mu}_0^\top (I - A) \hat{\mu}_0 -\mu^\top(I-A)\mu\right] + \frac{1}{n_0}(A - I)\right|  \leq C(\sqrt{\log p}(\|A - I\|_F/n)+\sqrt{\log p}\|(I - A)\mu\|/\sqrt{n}).
\ee
\end{lemma*}

 Let $R = \hat{\mu}_0 + \mu$, then $R \sim N(0, \frac{1}{n_0} I_p)$. 
 Therefore, the term of interest can be rewritten as 
\be
\ba{lll}
& & (R - \mu)^\top (I - A) (R - \mu) - \mu^\top (I - A) \mu + \frac{1}{n_0}(A - I)\\
&= & R^\top (I - A) R + \frac{1}{n_0}(A - I) - 2 \mu^\top (I - A) R\\
& = & Ia + 2Ib, 
\ea
\ee
where $Ia = R^\top (I - A) R + \frac{1}{n_0}(A - I)$ and $Ib = \mu^\top (I - A) R$. 

$Ia$ is a quadratic form of $R$. Hence, according to Lemma \ref{lemma:quad}, 
\ben
\ba{rl}
E[Ia] & =  \frac{1}{n_0} Tr(I - A) + \frac{1}{n_0}(A - I)= 0, \\
\var(Ia) & =  \frac{2}{n_0^2}Tr((I - A)^2)  = 2c^2n^{-2} \|A - I\|_F^2,
\ea
\een
where $c = n/n_0$. 
Further, we can find $\sup_x|F_{Ia/\sqrt{\var(IIa)}}(x)-\Phi(x)|\overset{\mathcal P}\to 0$ by the Berry-Ess\'een theorem, so $|Ia| \leq \sqrt{\log p}(\|A - I\|_F/n)$. 

Now we consider $Ib \sim N(0, \frac{1}{n_0} \mu^\top (I - A)^2 \mu)$.
Furthermore, $\mu^\top (I - A)^2 \mu \leq \|I - A\|^2\|\mu\|^2 $. 
Therefore, 
\[
|Ib| \leq \sqrt{\log p}\|(I - A)\mu\|/\sqrt{n}.
\]
Combining the results about $Ia$ and $Ib$, with probability $1 + o(1)$, 
\be\label{eqn:Tp}
|\Delta T| \leq |Ia| + |Ib| \lesssim \sqrt{\log p}(\|A - I\|_F/n)+\sqrt{\log p}\|(I - A)\mu\|/\sqrt{n}.
\ee

So, the result is proved. \qed

\subsection{Proof of Lemma \ref{lemma1}}
\begin{lemma*}
Under models (\ref{Omega0}) and (\ref{Param}), when $p\to\infty$, with probability $1-o(1)$, 
\ben
\|V^{(k)}\| \leq \eta_p b(p, \beta) = \left\{\begin{array}{ll}
3\eta_p\sqrt{p\nu} = 3\eta_pp^{(1-\beta)/2}, & 0<\beta < 1,\\
2\eta_p\sqrt{{\log p}/{\log\log p}}, & \beta = 1,\\
{2\eta_p}/{(\beta - 1)}, & 1 < \beta \leq 2.
\end{array}
\right.
\een
\end{lemma*}

{\it Proof}. In the proof, we use $V$ as $V^{(k)}$ for short. Further, since all the entries in $V$ has the same magnitude $\eta$, we define 
\[
W = \frac{1}{\eta} V. 
\]
Then we only need to prove the spectral norm about $W$. 
We prove it in two cases, the matrix is relatively sparse and the matrix is extremely sparse. 
\begin{itemize}
\item Not very sparse region ($0<\beta < 1$). 

Refer to Theorem 1.5 in \cite{vu2005spectral}. Note that, in our setting, the off-diagonal entries $|w_{ij}| \leq 1$, $E[w_{ij}] = 0$, and $Var(w_{ij}) = \nu \geq p^{-1} \log^4 p$. Then, almost surely, 
\[
\|W\| \leq 2\sqrt{p \nu} + C (p \nu)^{1/4}\log p.
\]
In this region, $p \nu \rightarrow \infty$, so $\sqrt{p\nu} \gg (p\nu)^{1/4}\log p$ when $p \rightarrow \infty$. Therefore, almost surely, 
\begin{equation}\label{eqn:dense1}
\|W\| \leq 2\sqrt{p \nu}(1+o(1)) \leq 3 \sqrt{p \nu}.
\end{equation}

\item Sparse region ($1 \leq \beta < 2$). 

In this region, we apply theorems in random graph theory to prove our results. 

Let $U$ be the matrix where the sign is removed for the entries of $W$, i.e., $u_{ij} = |w_{ij}|$, $1 \leq i, j \leq p$. Therefore, $U$ is also a Wigner matrix, with diagonals 0, and off-diagonals
\[
u_{ij} = u_{ji} \stackrel{i.i.d}{\sim} Bernoulli(\nu).
\]
Clearly, $\lambda_1(W) \leq \lambda_1(U)$, and $\lambda_n(W) \geq -\lambda_1(U)$. Therefore, $\|W\| \leq \|U\|$. To control $\|W\|$, we only need to control $\|U\|$. 

Recall the Erdos–Renyi random-graph model in graph theory. For a graph $G = (V, E)$, where $V$ is the vertex set and $E$ is the edge set, the Erdos–Renyi undirected random-graph model $G = G(p, \nu)$ denotes a graph with $V = \{1, 2, \dots, p\}$, and $(i, j) \in E$ with probability $\nu$. Therefore, the adjacency matrix of $G$ is a Wigner matrix with the same distribution of $U$. The results for the largest eigenvalue of $G$ also apply to $\|U\|$. 

\begin{itemize}
\item When $\beta = 1$ ($\nu = p^{-1}$). 

With Corollary 1.2 in \cite{krivelevich2003largest}, almost surely, 
\[
\lambda_1(G(p, 1/p)) = (1 + o(1))\sqrt{\frac{\log p}{\log\log p}}.
\]
Therefore, almost surely,
\begin{equation}\label{eqn:sparse1}
\|U\| \leq 2\sqrt{\frac{\log p}{\log\log p}}.
\end{equation}

\item When $1 < \beta \leq 2$ ($p^{-2} \leq \nu < p^{-1}$). 

Let $\|U\|_1 = \max_{1\leq i \leq p} \sum_{j=1}^p |u_{ij}|$ denote the 1-norm of $U$, and $\|U\|_{\infty} = \max_{1\leq j \leq p} \sum_{i=1}^p |u_{ij}|$ denote the infinity norm of $U$. Since $U$ is symmetric, $\|U\|_1 = \|U\|_{\infty}$. According to the relationships between matrix norms, we have
\[
\|U\|  \leq \sqrt{\|U\|_1\|U\|_{\infty}} = \|U\|_1.
\]
Thus we only need to control $\|U\|_1$, the maximal degree of nodes.

According to Lemma 2.2 in \cite{krivelevich2003largest}, the maximum degree of the graph almost surely satisfies the following:
\[
\mbox{maximum degree} = (1+o(1)) \max\{k: p{p-1 \choose k} \nu^k (1 - \nu)^{p-k} \geq 1\}.
\]
Now we want to have an upper bound for $k$. Note that
\begin{eqnarray*}
p{p-1 \choose k} \nu^k (1 - \nu)^{p-k}  \leq  p*p^k \nu^k \leq p(p\nu)^k.
\end{eqnarray*}
When $k > 1/(\beta - 1) $, since $p\nu < 1$, $p(p\nu)^k < p(p\nu)^{1/(\beta-1)} = 1$. Therefore, $p{p-1 \choose k} \nu^k (1 - \nu)^{p-k} < 1$, and so $\max\{k: p{p-1 \choose k} \nu^k (1 - \nu)^{p-k} \geq 1\} \leq 1/(\beta - 1)$.
The derivation shows that, almost surely, 
\begin{equation}\label{eqn:sparse4}
\|U\| \leq 2/(\beta - 1).
\end{equation}

\end{itemize}
Combine (\ref{eqn:sparse1}) and (\ref{eqn:sparse4}) about $\|U\|$, we can see that $\|W\| \leq \|U\|$ follows the same bound. 
\end{itemize}
Combing the results in two regions, the results are proved. \qed

\subsection{Proof of Lemma \ref{lemmaeb}}
\begin{lemma*}
Consider model (\ref{M12}) with the parameterizations (\ref{Param1}) and (\ref{Omega1})--(\ref{con1}). With probability $1 - o(1)$, we have $\|\Omega_k - I\|= o(1)$ and 
\be\ba{l}
\|V^{(k)}\|_F^2  =  \eta_p^2 p^2 \nu_p(1 + o(1)), \quad 
\|\Omega_k - I\|_F^2  =  p\xi_p^2+ \eta_p^2 p^2 \nu_p(1 + o(1)),\\
\|\mu\|^2 = p\tau_p^2  \epsilon_p(1 + o(1)). 
 \ea
\ee
\end{lemma*}

{\it Proof}. 
Let $\mathcal{B}_p \sim Binomial(\frac{p(p-1)}{2}, \nu)$ and $\mathcal{B}'_p \sim Binomial(p, \epsilon)$. According to the model, we have 
\[
\|V^{(k)}\|_F^2 = \eta_p^2 \mathcal{B}_p, \qquad 
\|\mu\|^2 = \tau^2 \mathcal{B}'_p, \qquad 
\|\Omega_k - I\|^2 = p\xi_p^2 + \eta_p^2 \mathcal{B}_p.
\]
We want to analyze three terms: $p\xi^2$, $\eta^2 \mathcal B_p$, and $\tau^2 \mathcal B'_p$. 
\begin{itemize}
\item The term $p\xi^2 = p^{1-2\gamma}$. Hence, under the condition that $\gamma \geq 1/2$, clearly we have
\[
p\xi^2 \goto \left\{
\ba{lll}
1, & \gamma = 1/2,\\
0, & \gamma > 1/2. 
\ea
\right. 
\]
\item The term $\eta^2 \mathcal B_p$, where $\mathcal B_p \sim Binomial(\frac{p(p-1)}{2}, \nu)$. 
According to Bernstein’s inequality, 
\begin{eqnarray*}
 P(|\mathcal B_p -  \frac{p(p-1)}{2}\nu| \geq \sqrt{p^2\nu} \log p)  \leq
 2\exp\{-\frac{(\sqrt{p^2\nu} \log p)^2/2}{\frac{p(p-1)}{2}\nu(1 - \nu) + (\sqrt{p^2\nu} \log p)/3}\} = o(p^{-1}).
\end{eqnarray*}
Note that $p^2 \nu \goto \infty$. 
Hence, with probability $1-o(1)$, $\mathcal B_p = p^2\nu/2(1+o(1))$, and
\[
\eta^2 \mathcal B_p = \eta^2 p^2\nu/2(1+o(1)). 
\]
Therefore, under the condition $\beta<2-2\alpha$, $f_p=1$, or the condition $\beta=2-2\alpha$, $f_p=L_p$ (a $\log p$ term), we have $\eta^2 p^2 \nu \goto \infty$, which means that $\eta^2 \mathcal B_p \goto \infty$. Under the condition that $\beta>2-2\alpha$, we have $\eta^2 p^2 \nu \goto 0$, which indicates that $\eta^2 \mathcal B_p \goto 0$. 

\item The term $\tau^2 \mathcal B'_p$, where $\mathcal B'_p \sim Binomial(p, \epsilon)$. 
According to Bernstein's inequality, 
\begin{eqnarray*}
  P(|\mathcal B'_p - p\epsilon| \geq \sqrt{p\epsilon} \log p)  &\leq&
 2\exp\{-\frac{(\sqrt{p\epsilon} \log p)^2/2}{p\epsilon(1-\epsilon) + (\sqrt{p\epsilon} \log p)/3}\} = o(p^{-1}).
\end{eqnarray*}
Since $p\epsilon \goto \infty$, we get, with probability $1-o(1)$, 
\[
\tau^2 \mathcal B'_p = \tau^2 p\epsilon(1+o(1)). 
\]
Under the condition $0<\zeta<1-2\theta$, $g_p=1$, or the condition $\zeta=1-2\theta$, $g_p=L_p$, we have $\tau^2 p \epsilon \goto \infty$; hence, $\tau^2 \mathcal B'_p \goto \infty$. Under the condition $\zeta>1-2\theta$, $\tau^2 p \epsilon \goto 0$, we have $\tau^2 \mathcal B'_p \goto 0$. 
\end{itemize}
Combining the results, when any one of the following four conditions holds, 
\begin{itemize}
\item [(1)]  $\beta<2-2\alpha$;
\item [(2)] $0<\zeta<1-2\theta$;
\end{itemize}
$\eta^2 \mathcal B_p \goto \infty$ or $\tau^2 \mathcal B'_p \goto \infty$ with probability $1-o(1)$. Therefore, we have 
\[
p\xi^2 + \eta^2 \mathcal B_p + c\tau^2 \mathcal B'_p \goto \infty, \quad c > 0.
\]
If $\beta>2-2\alpha$ and  $\zeta>1-2\theta$, $\eta^2 \mathcal B_p \goto 0$ or $\tau^2 \mathcal B'_p \goto 0$. Hence the summation of them also converges to 0. 

Combining the results, the lemma is proved. \qed

\subsection{Proof of Lemma \ref{lemma:StrongDeltaQ}}
\begin{lemma*}
Under the model assumptions and the definition of $\Delta Q$ in (\ref{Decomp2}), there is
\ben
 |\Delta Q|   \leq   O(\sqrt{p \epsilon (\xi^2 + p\epsilon \eta^2 \nu)}/n) + \sqrt{p\epsilon} \tau \log p (1 + o(1)).
\een
\end{lemma*}

{\it Proof}. 
Recall that $\hat{\mu}^{(t)} = \hat{\mu}_0 \circ d^{(t)}$ and $\hat{\mu}_d^{(t)} = d \circ d^{(t)}$. 
For simplicity, in this section, we use $\hat{\mu}_0$ and $d$ to denote $\hat{\mu}_0 \circ d^{(t)}$ and $\hat{\mu}_d^{(t)}$, respectively. Note that since all the signals are exactly recovered, $\hat{\mu}_0$ has zeros on the non-signal entries and non-zeros on the signals, and the same for $d$. We use $V$ to denote $V^{(k)}$ and $D_{\Omega}$ to denote $D_{\Omega}^{(k)}$. Further, we define $W = \frac{1}{\eta} V$ to focus on the location of non-zeros in $V$.

Let $k=\|\mu\|_0$ denote the number of non-zeros in $\mu$. 
Without loss of generality, we permute $\mu$ such that the first $k$ entries are the non-zeros and the rest are the zeros. We also permute $V$, $\Omega$, and $X$ accordingly, and rewrite $W$ and $\Omega$ as $2\times 2 $ block matrices $W=\left(^{W_{11} \ W_{12}}_{W_{21} \ W_{22}}\right)$ and $\Omega=\left(^{\Omega_{11} \ \Omega_{12}}_{\Omega_{21} \ \Omega_{22}}\right)$ respectively, where $W_{11}$ and $\Omega_{11}$ are $k\times k$ sub-matrices of $W$ and $\Omega$,  respectively. Let $X^{(k)}$, $d^{(k)}$, $\mu^{(k)}$, and $\hat{\mu}_0^{(k)}$ denote, respectively, $X$, $d$, $\mu$, and $\hat{\mu}_0$ restricted on the first $k$ entries, and let $X^{(p-k)}$ denote $X$ restricted on the last $(p-k)$ entries. Then $\mu^{(k)}$ is a length $k$ vector with all elements as $\tau$. 

With all the notations, $\Delta Q$ is 
\be\label{Decomp2app}
\ba{rl}
 \Delta Q = & \displaystyle 2(d - (I + \Omega)\mu)^\top X+\left[\hat{\mu}_0^\top(I-\Omega)\hat{\mu}_0-\mu^\top(I-\Omega)\mu+ \frac{1}{n_0}Tr(\Omega^{(d)} - I)\right] \\
 = & 2I_k + II_k.
\ea\ee

Now we analyze $I_k$ and $II_k$; the result will include $k$. Recall that $k$ is the number of non-zeros in $\mu$, where $k \sim Binomial(p, \epsilon)$. According to Bernstein's inequality, 
\begin{eqnarray*}
  P(|k - p\epsilon| \geq \sqrt{p\epsilon} \log p)  &\leq&
 2\exp\{-\frac{(\sqrt{p\epsilon} \log p)^2/2}{p\epsilon(1-\epsilon) + (\sqrt{p\epsilon} \log p)/3}\} = o(p^{-1}).
\end{eqnarray*}
Since $p\epsilon \goto \infty$, with probability $1-o(1)$, we have $k =  p\epsilon(1+o(1))$.

\begin{itemize}
\item We consider $I_k$ first. Since $d = \left(^{d^{(k)}}_{0_{p-k}}\right)$, $\mu = \left(^{\mu^{(k)}}_{0_{p-k}}\right)$, and $X = \left(^{X^{(k)}}_{X^{(p-k)}}\right)$, where $0_{p-k}$ is a zero vector with length $p - k$,
\ben
\ba{lll}
I_k & = & \left({d^{(k)}}^\top \,\, {0_{p-k}}^\top \right)  
\left(\ba{l}
X^{(k)}\\
X^{(p-k)}
\ea\right) - 
\left({\mu^{(k)}}^\top \,\, {0_{p-k}}^\top \right)  (I + \Omega)
\left(\ba{l}
X^{(k)}\\
X^{(p-k)} 
\ea \right) \\
& = & (d^{(k)})^\top X^{(k)} - (\mu^{(k)})^\top( I +  \Omega_{11}) X^{(k)} - (\mu^{(k)})^\top \Omega_{12} X^{(p-k)}\\
& = &  (d^{(k)} - (I + \Omega_{11}) \mu^{(k)})^\top X^{(k)} - (\mu^{(k)})^\top \Omega_{12} X^{(p-k)}\\
& = & Ia + Ib. 
\ea
\een

Consider $Ia$ first. Let $\tilde{d} = d^{(k)} - (I + \Omega_{11})\mu^{(k)}$, then $Ia = \tilde{d}^\top X^{(k)}$, and 
\[
\tilde{d} \sim N(0, \frac{1}{n_0}I + \frac{1}{n_1} \Omega_{11}), \quad 
X^{(k)} \sim N(-\mu^{(k)}, I). 
\]
Noting that $\tilde{d}$ is independent with $X^{(k)}$, there is $E[Ia] = 0$. 
Again by Lemma \ref{lemma:quad}, 
\ben\ba{lll}
\var(Ia)  & = &\displaystyle \frac{1}{n_0}(k + \|\mu^{(k)}\|^2) + \frac{1}{n_1}(k(1 + \xi) + (\mu^{(k)})^\top \Omega_{11} \mu^{(k)})\\
& \leq &\displaystyle (\frac{1}{n_0} + \frac{1}{n_1})k(1+o(1)) + \frac{1}{n_0}\|\Omega_{11}\| k\tau^2 \\
& \lesssim &\displaystyle   \ 4n^{-1}k = 4p\epsilon/n(1+o(1)),
\ea
\een
For the case in which $X \sim N(\mu, \Omega^{-1})$, the same result is obtained.

We also prove the aymptotic normality according to Lemma \ref{lemma2} and the Berry-Ess\'een theorem. Therefore, $\sup_x|F_{Ia/\sqrt{\var(Ia)}}(x)-\Phi(x)|\overset{\mathcal P}\to 0$, and, hence, $Ia = O_p(\sqrt{4n^{-1}p\epsilon})$. 

Next, consider $Ib$.  Recall that $X^{(p-k)} \sim N(0, I)$. Therefore, 
\[
Ib = (\mu^{(k)})^\top  \Omega_{12}X^{(p-k)} \sim N(0, (\mu^{(k)})^\top \Omega_{12} \Omega_{12}^\top \mu^{(k)}). 
\]
For the variance term, since $\Omega_{12} = \Omega_{21}^\top $, we have $(\mu^{(k)})^\top\Omega_{12} \Omega_{21} \mu^{(k)} \leq \|\mu^{(k)}\|^2 \| \Omega_{12} \Omega_{21}\| =  \|\Omega_{12}\|^2*k\tau^2$. 
By $\Omega = D_{\Omega} + V$,  we have $\Omega_{12}=V_{12}$, and hence $\|\Omega_{12}\| = \| V_{12}\| \leq  \|V\|$. 
According to Lemma \ref{lemma1} and the condition $\eta^2 p \nu \goto 0$, with probability $1 - o(1)$, $\|V\|^2 \leq  \eta^2 \log p \goto 0$ when $\beta\geq 1$ and $\|V\|^2 \leq C\eta^2p\nu \goto 0$ when $\beta < 1$. As a result, with probability $1-o(1)$, 
\be
(\mu^{(k)})^\top\Omega_{12} \Omega_{21} \mu^{(k)} \leq  k \tau^2 \| \Omega_{12}\|^2 \leq k\tau^2  \|V\|^2 = o(k\tau^2).
\ee
So, with probability $1 - o(1)$, 
\be\label{thm3case2part1}
|Ib| \leq \sqrt{k} \tau \log p = \sqrt{p\epsilon} \tau\log p.
\ee
For the case in which $X \sim N(\mu, \Omega^{-1})$, the analysis is similar. 

To conclude, and noting that $\tau \gg 1/\sqrt{n}$ in this region, we have 
\be\label{eqn:Ik}
|I_k| \leq |Ia| + |Ib| \lesssim \sqrt{p\epsilon} \tau\log p + O_p(\sqrt{4n^{-1}p\epsilon}) = \sqrt{p\epsilon} \tau\log p (1+o(1)). 
\ee

\item 
Next, we analyze $II_k$. Removing the zero part, we can find 
\[
II_k = -(\hat{\mu}_0^{(k)})^\top (I - \Omega_{11}) \hat{\mu}_0^{(k)}+ (\mu^{(k)})^\top (I - \Omega_{11}) \mu^{(k)}+ \frac{1}{n_0}Tr(\Omega_{11} - I).
\]  
Let $R = \hat{\mu}_0^{(k)} + \mu^{(k)}$, then $R \sim N(0, \frac{1}{n_0} I_k)$. 
Rewrite $II_k$ as 
\be
II_k = R^\top (I - \Omega_{11}) R + 2 (\mu^{(k)})^\top (I - \Omega_{11}) R = IIa + 2IIb. 
\ee

We first consider $IIa = R^\top (I - \Omega_{11}) R - \frac{1}{n_0} Tr(I - \Omega_{11})$. This is a quadratic form of $R$. Hence, by Lemma \ref{lemma:quad}, 
\ben
\ba{lll}
E[IIa] & = & 0, \\ 
\var(IIa) & = & \frac{2}{n_0^2} Tr((I - \Omega_{11})^2) = 2c^2(k\xi^2 + k^2 \eta^2 \nu)/n^2\\
& = & 2c^2 p \epsilon (\xi^2 + p\epsilon \eta^2 \nu)/n^2(1+o(1)),
\ea
\een
where $c = n/n_0$. 
Furthermore, we can prove that $\sup_x|F_{IIa/\sqrt{\var(IIa)}}(x)-\Phi(x)|\overset{\mathcal P}\to 0$, and so 
\[
IIa =  O(n^{-1}\sqrt{p \epsilon (\xi^2 + p\epsilon \eta^2 \nu)}). 
\]

Then, we consider $IIb =  (\mu^{(k)})^\top (I - \Omega_{11}) R$. Since $Z\sim N(0, \frac{1}{n_0} I)$, it is clear that $IIb \sim N(0, \frac{1}{n_0} (\mu^{(k)})^\top (I - \Omega_{11})^2 \mu^{(k)})$. 
According to Lemma \ref{lemmaeb}, $\|I - \Omega_{11}\| \leq \|I - \Omega\| = o(1)$.
Therefore, $(\mu^{(k)})^\top (I - \Omega_{11})^2 \mu^{(k)}) \leq \|\mu^{(k)}\|^2 \|I - \Omega\|^2 = o(k\tau^2)$. 
As a result, with $k = p\epsilon (1+o(1))$, 
\[
|IIb| \leq \sqrt{\frac{1}{n_0} (\mu^{(k)})^\top (I - \Omega_{11})^2 \mu^{(k)}} \log p = o(\sqrt{n^{-1} p\epsilon\tau^2}). 
\]
Combining the results for $IIa$ and $IIb$, we have
\be\label{eqn:IIk}
|II_k| \leq |IIa| + |IIb| \lesssim O(n^{-1}\sqrt{p \epsilon (\xi^2 + p\epsilon \eta^2 \nu)})+ o(\sqrt{p\epsilon\tau^2/n}).
\ee
\end{itemize}

Combining the results for $I_k$ and $II_k$ in (\ref{eqn:Ik}) and (\ref{eqn:IIk}), 
\ben
\ba{lll}
\Delta Q & \leq & \sqrt{p\epsilon} \tau\log p (1+o(1)) + O(n^{-1}\sqrt{p \epsilon (\xi^2 + p\epsilon \eta^2 \nu)})+ o(\sqrt{p\epsilon\tau^2/n}) \\
& = &  O(n^{-1}\sqrt{p \epsilon (\xi^2 + p\epsilon \eta^2 \nu)}) + \sqrt{p\epsilon} \tau \log p (1 + o(1)). 
\ea
\een
The result is proved. \qed

\subsection{Proof of Lemma \ref{lemma:unknownstrongQ}}
\begin{lemma*}
Under the model assumptions and the definition of $\Delta Q$, there is
\be
 |\Delta Q|   \leq  \eta \tau \max\{p\epsilon \nu, 1\} \log p +  O_p(\sqrt{4n^{-1}p\epsilon}).
\ee
\end{lemma*}
{\it Proof}. 
Recall that $\hat{\mu}^{(t)} = \hat{\mu}_0 \circ d^{(t)}$ and $\hat{\mu}_d^{(t)} = d \circ d^{(t)}$. 
For simplicity, in this section, we use $\hat{\mu}_0$ and $d$ to denote $\hat{\mu}_0^{(t)}$ and $\hat{\mu}_d^{(t)}$, respectively. Since all the signals are exactly recovered, $\hat{\mu}_0$ and $d$ have zeros on the non-signal entries and non-zeros on the signals.

Let $k=\|\mu\|_0$ denote the number of non-zeros in $\mu$. 
Without loss of generality, we permute $\mu$ such that the first $k$ entries are the non-zeros and the rest are the zeros. Permute $W$, $\Omega$, $\hat\Omega$ and $X$ accordingly, and rewrite $W$ and $\Omega$ as $2\times 2 $ block matrices $W=\left(^{W_{11} \ W_{12}}_{W_{21} \ W_{22}}\right)$ and $\Omega=\left(^{\Omega_{11} \ \Omega_{12}}_{\Omega_{21} \ \Omega_{22}}\right)$, where $W_{11}$ and $\Omega_{11}$ are $k\times k$ sub-matrices of $W$ and $\Omega$,  respectively. Let $X^{(k)}$, $d^{(k)}$, $\mu^{(k)}$, and $\hat{\mu}_0^{(k)}$ denote, respectively, $X$, $d$, $\mu$, and $\hat{\mu}_0$ restricted on the first $k$ entries, and let $X^{(p-k)}$ denote $X$ restricted on the last $(p-k)$ entries. Then $\mu^{(k)}$ is a length $k$ vector with all elements as $\tau$. 

With all the notations, $\Delta Q$ is 
\be\label{Decomp2appunknown}
\ba{rl}
 \Delta Q = & \displaystyle 2(d - (I + \hat\Omega)\mu)^\top X+\left[\hat{\mu}_0^\top(I-\hat\Omega)\hat{\mu}_0-\mu^\top(I-\hat\Omega)\mu\right] \\
 = & 2I_k + II_k.
\ea\ee

Now we analyze $I_k$ and $II_k$. The discussion focuses on the case $Y = 0$, i.e. $X \sim N(-\mu, I)$. The derivation for $Y = 1$, i.e., $X \sim N(\mu, \Omega^{-1})$ is similar and the results are at the same order.  The result will include $k$. Recall that $k$ is the number of non-zeros in $\mu$, where $k \sim Binomial(p, \epsilon)$. According to Bernstein's inequality, 
\begin{eqnarray*}
  P(|k - p\epsilon| \geq \sqrt{p\epsilon} \log p)  &\leq&
 2\exp\{-\frac{(\sqrt{p\epsilon} \log p)^2/2}{p\epsilon(1-\epsilon) + (\sqrt{p\epsilon} \log p)/3}\} = o(p^{-1}).
\end{eqnarray*}
Since $p\epsilon \goto \infty$, with probability $1-o(1)$, we have $k =  p\epsilon(1+o(1))$.

\begin{itemize}
\item We consider $I_k$ first. Since $d = \left(^{d^{(k)}}_{0_{p-k}}\right)$, $\mu = \left(^{\mu^{(k)}}_{0_{p-k}}\right)$, and $X = \left(^{X^{(k)}}_{X^{(p-k)}}\right)$, where $0_{p-k}$ is a zero vector with length $p - k$,
\ben
\ba{lll}
I_k & = & \left({d^{(k)}}^\top \,\, {0_{p-k}}^\top \right)  
\left(\ba{l}
X^{(k)}\\
X^{(p-k)}
\ea\right) - 
\left({\mu^{(k)}}^\top \,\, {0_{p-k}}^\top \right)  (I + \hat\Omega)
\left(\ba{l}
X^{(k)}\\
X^{(p-k)} 
\ea \right) \\
& = & (d^{(k)})^\top X^{(k)} - (\mu^{(k)})^\top( I +  \hat\Omega_{11}) X^{(k)} - (\mu^{(k)})^\top \hat\Omega_{12} X^{(p-k)}\\
& = &  (d^{(k)} - (I + \hat\Omega_{11}) \mu^{(k)})^\top X^{(k)} - (\mu^{(k)})^\top \hat\Omega_{12} X^{(p-k)}\\
& = & Ia + Ib. 
\ea
\een

Consider $Ia$ first. Let $\tilde{d} = d^{(k)} - (I + \hat\Omega_{11})\mu^{(k)}$, then $Ia = \tilde{d}^\top X^{(k)}$, and 
\[
\tilde{d} \sim N(0, \frac{1}{n_0}I + \frac{1}{n_1} (\hat\Omega\Omega^{-1}\hat\Omega)_{11}), \quad 
X^{(k)} \sim N(-\mu^{(k)}, I). 
\]
$\tilde{d}$ is independent with $X^{(k)}$, so $E[Ia] = 0$. 
The variance can be obtained by the law of total variance, that 
\ben\ba{lll}
\var(Ia)  & = &\displaystyle \frac{1}{n_0}(k + \|\mu^{(k)}\|^2) + \frac{1}{n_1}(Tr((\hat\Omega\Omega^{-1}\hat\Omega)_{11}) + (\mu^{(k)})^\top (\hat\Omega\Omega^{-1}\hat\Omega)_{11} \mu^{(k)})\\
& \leq &\displaystyle (\frac{1}{n_0} + \frac{1}{n_1})k(1+o(1)) + \frac{1}{n_0}k\tau^2 + \frac{1}{n_1}\|(\hat\Omega\Omega^{-1}\hat\Omega)_{11}\| k\tau^2 \\
& \lesssim &\displaystyle   \ 4n^{-1}k = 4p\epsilon/n(1+o(1)),
\ea
\een
where the trace of $(\hat\Omega\Omega^{-1}\hat\Omega)_{11}$ is constraied by $k\|(\hat\Omega\Omega^{-1}\hat\Omega)_{11})\| = k(1+o(1))$.
For the case in which $X \sim N(\mu, \Omega^{-1})$, the same result is obtained.

We also prove the aymptotic normality according to Lemma \ref{lemma2} and the Berry-Ess\'een theorem. Therefore, $\sup_x|F_{Ia/\sqrt{\var(Ia)}}(x)-\Phi(x)|\overset{\mathcal P}\to 0$, and, hence, $Ia = O_p(\sqrt{4n^{-1}p\epsilon})$. 

Next, consider $Ib$.  Recall that $X^{(p-k)} \sim N(0, I)$. Therefore, 
\[
Ib = (\mu^{(k)})^\top  \hat\Omega_{12}X^{(p-k)} \sim N(0, (\mu^{(k)})^\top \hat\Omega_{12} \hat\Omega_{12}^\top \mu^{(k)}). 
\]
For the variance term, since $\hat\Omega_{12} = \hat\Omega_{21}^\top $, we have $(\mu^{(k)})^\top\hat\Omega_{12} \hat\Omega_{21} \mu^{(k)} \leq \tau^2 \|\hat\Omega_{21}\|_\infty^2$. 
By $\Omega = c I + \eta W$,  we have $\Omega_{12}=\eta W_{12}$. Currently we require there are $o(\sqrt{n})$ non-zero entries in each row of $W$ and $\eta \gg 1/\sqrt{n}$. Further, the distribution on non-zeros in $W$ are independent with the non-zeros in $\mu$. Hence, with probability $1 - o(1)$,  
$\|\hat\Omega_{21}\|_\infty^2 \leq 2\eta^2 \max\{k^2\nu^2, 1\}$. 
As a result, with probability $1-o(1)$, 
\be
(\mu^{(k)})^\top\Omega_{12} \Omega_{21} \mu^{(k)} \leq  \tau^2 \|\hat\Omega_{21}\|_\infty^2 \leq 2\tau^2 \eta^2 \max\{k^2\nu^2, 1\}.
\ee
So, with probability $1 - o(1)$, 
\be\label{thm3case2part1unknown}
|Ib| \leq C\eta \tau \max\{k\nu, 1\} \log p= C\eta \tau \max\{p\epsilon\nu, 1\} \log p.
\ee
For the case in which $X \sim N(\mu, \Omega^{-1})$, the analysis is similar. 

To conclude, we have 
\be\label{eqn:Ikunknown}
|I_k| \leq |Ia| + |Ib| \lesssim \eta \tau \max\{p\epsilon\nu, 1\} \log p + O_p(\sqrt{4n^{-1}p\epsilon}). 
\ee

\item 
Next, we analyze $II_k$. Removing the zero part, we can find 
\[
II_k = -(\hat{\mu}_0^{(k)})^\top (I - \hat\Omega_{11}) \hat{\mu}_0^{(k)}+ (\mu^{(k)})^\top (I - \hat\Omega_{11}) \mu^{(k)} + Tr(\hat\Omega_{11} - I).
\]  
Let $R = \hat{\mu}_0^{(k)} + \mu^{(k)}$, then $R \sim N(0, \frac{1}{n_0} I_k)$. 
Rewrite $II_k$ as 
\be
II_k = [R^\top (I - \hat\Omega_{11}) R  + Tr(\hat\Omega_{11} - I)]+ 2 (\mu^{(k)})^\top (I - \hat\Omega_{11}) R = IIa + 2IIb. 
\ee

We first consider $IIa = R^\top (I - \Omega_{11}) R$. This follows a non-central chi-square distribution. Since $\eta \gg 1/\sqrt{n}$ and $\hat{\Omega}$ can recover exactly the non-zeros of $\Omega$, 
\ben
\ba{lll}
E[IIa] & = & 0, \\ 
\var(IIa) & = & \frac{2}{n_0^2} Tr((I - \hat\Omega_{11})^2).
\ea
\een
Furthermore, we can prove that $\sup_x|F_{IIa/\sqrt{\var(IIa)}}(x)-\Phi(x)|\overset{\mathcal P}\to 0$, and so 
\[
IIa =  O(n^{-1}\sqrt{Tr((I - \hat\Omega_{11})^2)}).
\]
If we introduce in the terms, then 
\[
\frac{2}{n_0^2}Tr((I - \hat{\Omega}_{11})^2) \leq Cn^{-2}(k\xi^2*1\{\xi \gg 1/\sqrt{n}\} + \eta^2\max\{k^2\nu, 1\})
\]
for some constant $C > 0$. 
And so 
\[
IIa = O\biggl(n^{-1}\sqrt{p \epsilon \xi^2*1\{\xi \gg 1/\sqrt{n}\} + \eta^2 \max\{p^2\epsilon^2 \nu, 1\}})\biggr).
\]

Then, we consider $IIb =  (\mu^{(k)})^\top (I - \hat\Omega_{11}) R$. Since $Z\sim N(0, \frac{1}{n_0} I)$, it is clear that $IIb \sim N(0, \frac{1}{n_0} (\mu^{(k)})^\top (I - \hat\Omega_{11})^2 \mu^{(k)})$. 
According to the definition of $\mu^{(k)}$, $\frac{1}{n_0} (\mu^{(k)})^\top (I - \hat\Omega_{11})^2 \mu^{(k)} = \frac{\tau^2}{n_0}\|I - \hat\Omega_{11}\|_{\infty}^2$. 
Therefore, with probability $1 - o(1/p)$, 
\[
|IIb| \leq n^{-1/2}\tau \|I - \hat{\Omega}_{11}\|_{\infty}.
\]
As a result, with $k = p\epsilon (1+o(1))$, 
\[
|IIb| \leq n^{-1/2}\tau(\xi*1\{\xi \gg 1/\sqrt{n}\} + 2\eta \max\{p\epsilon\nu, 1\}). 
\]
Combining the results for $IIa$ and $IIb$, we have
\be\label{eqn:IIkunknown}
|II_k| \leq |IIa| + |IIb| \lesssim O\biggl(n^{-1}\sqrt{p \epsilon \xi^2*1\{\xi \gg 1/\sqrt{n}\} + \eta^2 \max\{p^2\epsilon^2 \nu, 1\}})\biggr).
\ee
\end{itemize}

Combining the results for $I_k$ and $II_k$ in (\ref{eqn:Ikunknown}) and (\ref{eqn:IIkunknown}), 
\ben
\ba{lll}
\Delta Q & \leq & \eta \tau \max\{p\epsilon\nu, 1\} \log p + O_p(\sqrt{4n^{-1}p\epsilon}) + O\biggl(n^{-1}\sqrt{p \epsilon \xi^2*1\{\xi \gg 1/\sqrt{n}\} + \eta^2 \max\{p^2\epsilon^2 \nu, 1\}})\biggr) \\
& = &  \eta \tau \max\{p\epsilon \nu, 1\} \log p + O_p(\sqrt{4n^{-1}p\epsilon}). 
\ea
\een
The result is proved. \qed

\end{document}